\documentclass[11pt,a4paper]{ansarticle}      
\usepackage{setspace,graphicx}     

\usepackage{enumerate}  

\usepackage{amssymb}
\usepackage{mystyle}    
\usepackage{mythesisStyle}   
\usepackage{setspace}
\usepackage[numbers]{natbib}
\usepackage{subfigure}
\usepackage[toc,page]{appendix}
\usepackage{esint}
\usepackage{comment}
\usepackage{upgreek}
\DeclareGraphicsExtensions{.pdf,.png}

\newcommand{\remarkA}[1]{{\color{red}\large\bf Andreas: #1}}
\newcommand{\remarkC}[1]{{\color{blue}\large\bf Chinedu: #1}}

\title{A two-layer approach for Coupling 1D/2D Shallow Water Flow Models}
\author{Andreas Dedner \thanks{a.s.dedner@warwick.ac.uk}}
\author{Chinedu Nwaigwe \thanks{corresponding author, c.nwaigwe@warwick.ac.uk}}
\affil{Centre for Scientific Computing and Warwick Mathematics Institute, University of Warwick, United Kingdom}

\begin{document}
\maketitle

\section*{Abstract}
In this paper, we propose a novel approach for coupling 2D/1D shallow water
flows based on a two layer model in the channel,
a 1D lower layer model and a 2D upper laayer model. The upper layer is only
used in regions where flooding occurs otherwise the model reduces to
a standard 1D channel model. To switch between the one
layer and the two layer models the user prescribes
an elevation above which the channel is considered to be full, i.e., a
flooding event may be taking place.
In the case of flooding the 2D upper layer model make it strightforward to
couple the channel flow solver to a 2D shallow water solver used in the floodplain.
We show that the resulting method (i) is well-balanced (ii) preserves a
\emph{no-numerical flooding} property (iii) preserves conservation
properties of the underlying 1D and 2D finite volume schemes used for the
flow in the channel and the floodplain.  Numerical tests show that the method
performs well compared to two horizontal coupling methods found in the literature.
The results show that the
method recovers the 2D flow structure in the channel in flooding regions,
retains 1D flow structure in non-flooding regions while maintains good
efficiency.

\section{Introduction}
One dimensional Saint Venant models are often used to simulate open channel flows
but they become inadequate once the channel overflows. However,
multi-dimensional (even 2D shallow water) simulations
are computationally expensive. This has led to the development of methods
to couple 1D channel simulations with 2D floodplain simulations.
River/floodplain coupling simulations started as quasi 2D models in which floodplains are represented by
storage cells, then coupled with an existing 1D river model \cite{cungeetal80, bladeetal1994}.
These approaches would not allow to simulate the fluid dynamics in the floodplains \cite{nietocoupled}.

In \cite{bladeetal2012}, the coupling of a 1D channel model with a full
2D floodplain model was achieved by including the 2D numerical fluxes
into the finite volume scheme for the 1D model, while
\cite{yongcanetal2012} utilized the theory of characteristics to couple 1D/2D models
through suitable matching conditions defined at the 2D/1D interface.
In \cite{solomonetal2012},
the 1D river model and a 2D non-inertia model were also coupled where the water
level differences between the flows in the two domains are used to calculate the
interacting discharges in the sub-domains.

Methods based on post-processing the separately computed solutions are presented in
\cite{moralesetal2013}, see also \cite{moralesthesis}. At 2D/1D interface, each model computes its own
solution from which the total water volume in a 1D cell and all its adjacent 2D cells is computed.
Then, the water height for 2D cells and the wetted Area for the associated 1D cell are found.
In \cite{moralesapplied2016}, the methods have been applied to Tiber River, Rome.
In \cite{monnierMarin2009super}, the 1D model including the coupling terms were classically derived
from the full 3D inviscid Euler's equations and an optimal control process applied to couple the models.
These models have been numerically treated with the finite volume method \cite{nietocoupled}
where the discrete exchange term, which lead to globally well-balanced scheme, is proposed.
This approach superposes a 2D grid over the 1D channel
grid and convergence is achieved using a Schwartz-like iterative algorithm.

A major difficulty in coupling 2D/1D shallow water models is the
computation of the channel flow lateral discharge.
In \cite{comparison2015}, this channel lateral discharge
was set to zero, while \cite{nicoleetal20142d} adopted an iterative technique
that uses the solution of successive Riemann problems to estimate the transverse velocity.
This difficulty in computing the channel lateral discharge remains
challenging to compute accurate fluxed between the channel and the
floodplain.

A further issue is the assumption that the flow will remain one dimensional
in the channel even in the region, where a flooding event is occurring.
In \cite{chineduandreas1} a method to compute different lateral discharges
at each channel boundary was purposed.
However, the free-surface and $x$-velocity component were still assumed to be laterally
constant across the channel and although the approach does improve the
accuracy of the method, it still does not recover the complete 2D flow structure during flooding.

Efficient methods that recover the 2D flow structure during flooding but revert back to
1D simulation if no flooding is occurring could solve both the issues
mentioned above.
Frontal coupling methods in which the floodplain extends into the
channel in parts of the domain, recover 2D flow structure but
they loose efficiency because they compute the 2D solutions at all times.
Moreover, most of the existing methods need to know the location of possible flooding
a-priori and can not take into account that flooding locations may vary with time.
So methods that can adapt to flooding regions a-posteriori are desirable.

The goal of this paper is to propose a method, the vertical coupling method (VCM),
that (i) recovers the 2D flow structure during flooding while reverting back to
1D simulation in the channel regions where no flooding occurs.
(ii) automatically detects flooding regions.
(iii) can be easily added to well established 1D channel and 2D floodplain
flow solvers.

The VCM is based on partitioning the flow in an overflowing channel
into two vertical layers where the flows in the lower
and upper layers are simulated using 1D and 2D models, respectively.
We then derive a coupled models for the two layers which represents the channel flow.
The upper layer model can be easily coupled to a 2D floodplain model.
The method is parameterized by prescribing a height function $\zwall$ along the
channel which determines at which water height the channel is in danger of
overflowing. Only when this water height is reached will the upper layer
model be activated, below that layer the channel flow is simulated using a
standard 1D model. Different choices for $\zwall$ will lead to different
method and special choices will lead to some standard lateral or frontal
coupling methods found in the literature. In this sense VCM is a superset
of some existing methods.

The rest of the paper is organised as follows:
The two layer channel model is presented in section \ref{vertsecmodels}.
We start by presenting the notation used to derive the VCM
in section \ref{vertsecbackground}, derive the lower and upper layer flow models
in sections \ref{vertseclower} and \ref{vertsecupper} respectively,
the complete two layer channel flow model is summarised in section
\ref{vertsecchanelmodelsummary}.
For the floodplain flow we use a standard shallow water model which is
summarized in section \ref{vertsec2dswes}.
A numerical algorithm for the coupled channel flow models is presented in
detail in section \ref{vsecnumerical}
and the properties of the method are considered in section \ref{vertnumsecprops}
where we prove that the method
is well-balanced, preserves no-numerical flooding and is mass conservative.
Numerical experiments are presented in section \ref{vsecnumResults} to
evaluate the performance of the method compared to simpler approaches found
in the literature. Finally, we give a summary in section \ref{vsecconc}.

\section{Mathematical Models for the Fluid Dynamics}\label{vertsecmodels}
This section presents the model equations used for flow in the channel and
fllodplain and derive, in detail,
the two layer models for the channel flow.
\subsection{Background}\label{vertsecbackground}
\begin{figure}[ht!]
  \subfigure[The full flow cross section at a fixed point $x$ showing the 2D bottom topography, $z_b(x,y)$ comprising
         of the channel and the floodplains.]
           {\label{figFlowCrossSection} \includegraphics[width=0.47\textwidth, height=0.38\textwidth]{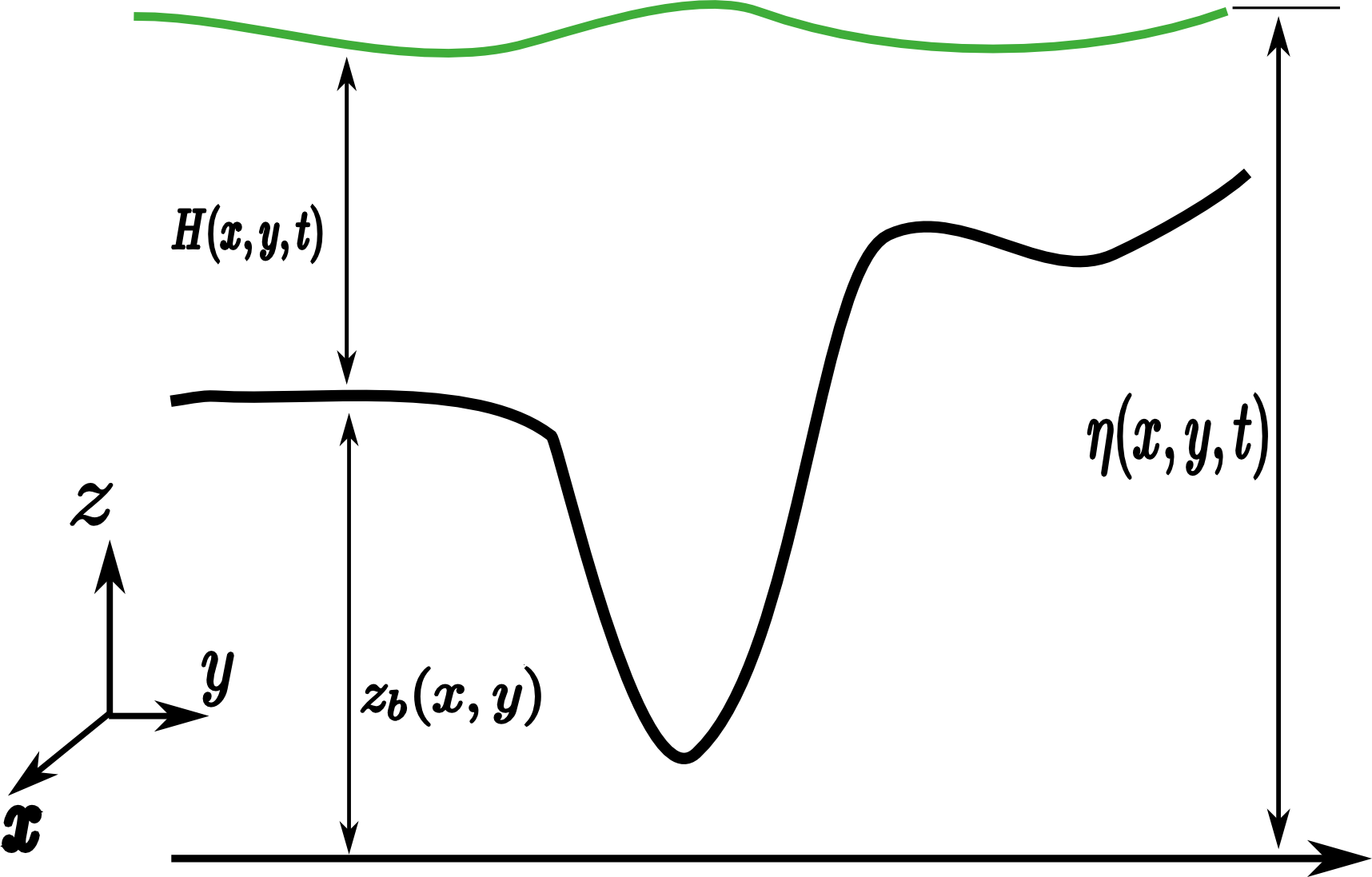} }
\hspace{0.02\textwidth}
  \subfigure[Vertical partitioning of the flow; water depths,  $h_1(x,y,t)$ and $h_2(x,y,t)$ in the lower and upper layers respectively.]
  {\label{fig-vcm-flow-partition}\includegraphics[width=0.47\textwidth, height=0.42\textwidth]{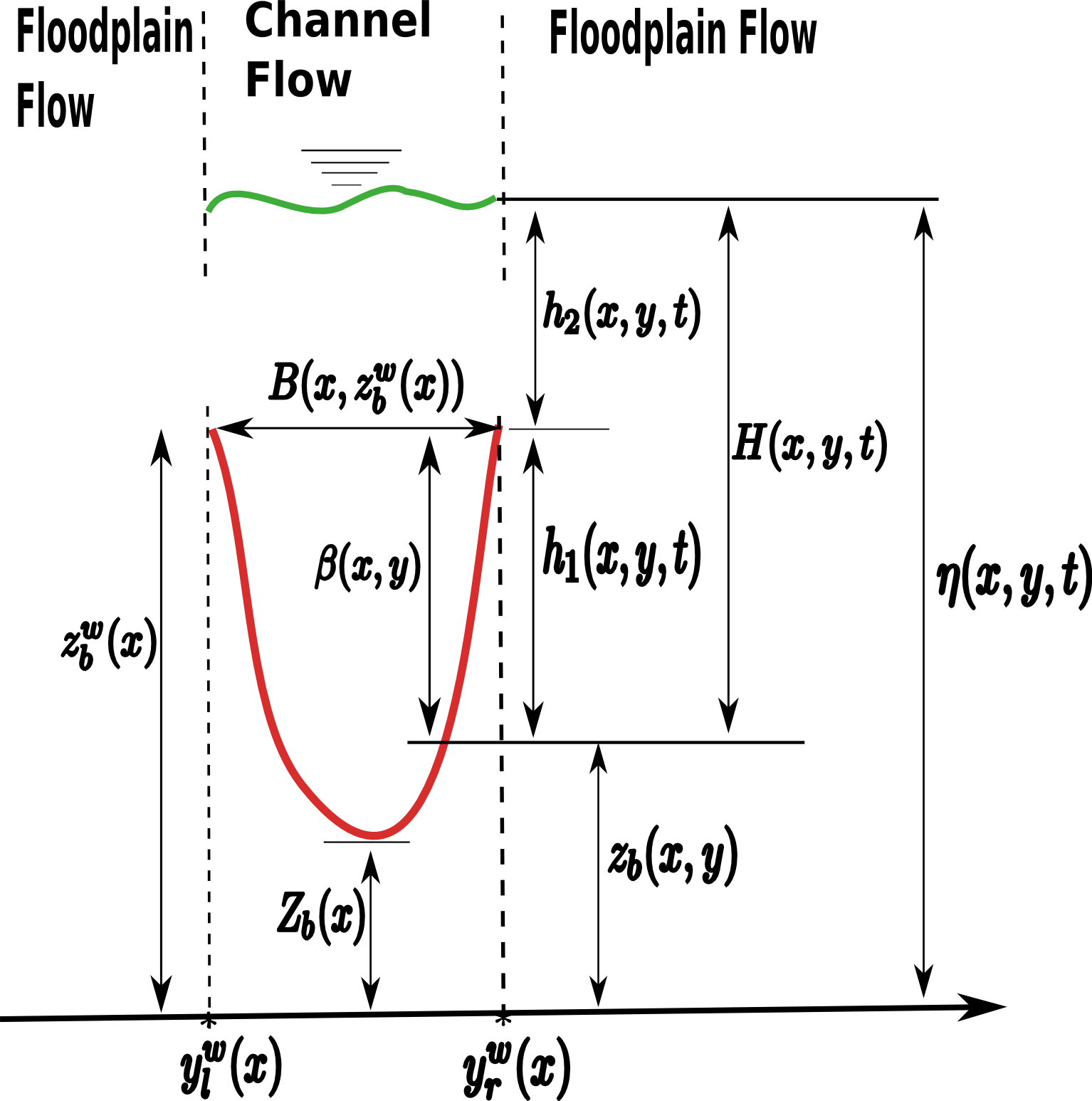} }
\caption{Illustration of flow cross section (left) and channel geometry (right).}
\end{figure}
Let $\Omega_H \subset \real^2$ be a fixed 2D horizontal domain with fixed bottom elevation, $\zb$, 
$\xy = (x,y)\in \Omega_H$. Let $H(\xy,t)\geq 0$ denote the depth of water at point, $\xy \in \Omega_H$ at
time, $t \ge 0$, so that
\begin{equation} \label{eq:waterHeight}
	 \eTa = z_b(\xy) + H(\xy,t)
\end{equation}
is the water free-surface elevation at point, $\xy$, at time $t$. Then, at time, $t$
the flow occupies the 3D domain, 
$\Omega_t$ defined by
\begin{equation} \label{moddomain3D}
	\Omega_t 
     	   =\{
			(\xy,z) \in \real^3 : \xy =(x,y) \in \Omega_H, \, 
			  z_b(\xy) \leq z \leq \eTa		
		   \}.
\end{equation}
A cross section of the flow domain, $\Omega_t$ at a fixed x is shown in 
figure \ref{figFlowCrossSection}.

\subsubsection{Channel Geometry:}
Figure \ref{fig-vcm-flow-partition} shows the flow cross section indicating the floodplain part and the
channel cross section. To simplify the derivation of the models we assume
its length lies along the $x$-axis (frontal direction) and the width along
the $y$-axis (lateral direction).
We define $Z_b(x)$ by
\begin{equation}\label{vcm-eqn-Zb}
  Z_b(x) = \min_{y}z_b(x,y).
\end{equation}
The functions $y_l(x,z)$ and $y_r(x,z)$ are the
$y$-coordinates of the left and right lateral wall boundaries respectively at the elevation, $z$.
The function, $B(x,z)$ gives the channel cross sectional lateral width,
i.e.,
	\begin{equation}\label{modwidthequation}
	  B(x,z) = \yrz - \ylz \quad \forall z \in \real,
   \end{equation}
and for convenience we set
\begin{equation}
	B(x,z) = 0, \quad \yrz = \ylz  \quad \mbox{ for all } z<Z_b(x),
\end{equation}
see figure \ref{fig-vcm-flow-partition}.

%
%
%

\subsubsection{The VCM Background}
Given the channel geometry above, the formulation of the VCM starts by choosing an elevation, $\zwall \geq Z_b(x)$
above which the channel is considered full, see figure
\ref{fig-vcm-flow-partition}.
This elevation is completely decided by the user; the only constraint 
is for it not to be less than $Z_b(x)$. 
So, $\zwall$ becomes the channel top and is referred to as the maximum wall elevation 
defined below.
\begin{mydefinition}[Maximum channel wall elevation, $z_b^w(x)$]\label{modzwalldef}
	The maximum channel wall elevation at cross section x, is the minimum elevation of the channel banks
  above which flooding is said to have occurred \cite{chineduthesis}. 
\end{mydefinition}
Once  $\zwall$  has been chosen, other quantities are derived.
The lateral width and the $y$-coordinates of the lateral wall boundaries at the top ($z=\zwall$) becomes
$B(x):=B(x,\zwall)$ and $y_{l,r}^w(x):=y_{l,r}(x,\zwall) $ respectively. 
For the VCM we then partition the flow into the channel flow within the region $\ylwall < y < \yrwall$,
and the flow in the floodplain the remaining region, namely $-\infty < y < \ylwall $ and $\yrwall < y < \infty$,
see figure \ref{fig-vcm-flow-partition}.
As noted above, the floodplain flows will be simulated using a standard 2D shallow water models
(given in section \ref{vertsec2dswes}), hence we concentrate on deriving
the models for the flow in the channel region.

\begin{myremark}
  The idea of the VCM method is to use a 1D channel model for regions in
  the domain where the channel is not full while using a two layer model
  otherwise, where only the lower layer is evolved using a 1D model
  while the upper layer uses a 2D model.
  As said above the choice of when to switch from a one layer to a two
  layer model depends on the choice of $\zwall$ and thus is up to the user.

  The closer the chosen $\zwall$ is to $Z_b(x)$,
  the smaller the channel flow region $\ylwall < y < \yrwall$, hence the larger
  the floodplain flow region (see figure \ref{fig-vcm-flow-partition}).
  In particular, choosing $\zwall=Z_b(x)$ in some part of the channel,
  will result in an approach where the 2D floodplain model is used in
  parts of the channel (a \emph{frontal} type coupling approach),
  while
  taking $\zwall=\infty$ will result in a one layer 1D model being used
  in that region leading to a \emph{horizontal} type coupling
  \cite{chineduthesis, chineduandreas1}.
  Thus, different choices of $\zwall$ lead to different coupling approaches.
\end{myremark}

Having identified the channel and floodplain flow regions, we now concentrate on the channel
region. First, we note that the following condition holds for fixed $x$:
\begin{align}\label{svmbanksGeneral}
  z_b(x,\ylz) = z_b(x,\yrz) = z  \quad  \forall z \in [ Z_b(x), z_b^w(x)   ].
\end{align}
We also extend the definition of the width functions, $B(x,z), \, y_{l,r}(x,z)$ to
the region above the top (that is where $z \ge \zwall$) as follows:
\begin{equation}\label{horBtop}
	    \begin{aligned}
          & B(x,z) = B(x, \zwall)  \quad \mbox{and} \quad y_{l,r}(x,z)=y_{l,r}^w(x)
        \end{aligned}          
          \quad  \forall  \quad z \ge \zwall,
\end{equation}
see figure \ref{fig-vcm-flow-partition}. This results in a channel with
straight vertical walls above the height where the channel is assumed to
be full.

For the scheme later on we need to assume that the channel geometry
is known such that we can reconstruct the water height given the wetted
cross section:
\begin{mydefinition}[Water height for given wetted cross section]
Let $A=A(x)$ be a given function describing the wetted area in the channel
defined by the bottom topography $z_b(x,\cdot)$.
We define the height function
$\Height(\xy;A):=\eta^*(x;A)-z_b(\xy)$ with
  $$ \eta^*(x;A) := \inf\{z \colon
           A(x) = \int_{y_l(x,z)}^{y_r(x,z)} z-z_b(x,y) \; dy \}~. $$
\end{mydefinition}
\begin{myremark}
With the above definition $\Height(\xy)-z_b(\xy)$ is independent of $y$
and represent a 1D flow with the wetter area given by $A$.
Due to the way we extended $y_{l,r}(x,z)$ for $z\ge\zwall$ it is easy to
see that the above definition provides a unique height function
$\Height$ even for large value of $A$, i.e., even in the case where the
given wetted cross section leads to a overfull channel.
\end{myremark}

\subsubsection{Important Quantities for the Channel Flow}
We now proceed to define other important quantities for the channel flow.
%
%
%

\begin{mydefinition}[Channel Depth]\label{vertbackdefbeta}
The channel depth, $\beta(\xy)$ is the laterally varying height between the channel bed and the chosen elevation, $\zwall$, that is
\begin{equation}\label{vertlowermodeqnbeta}
	\beta(\xy) = \zwall - \zb, \quad \quad \ylwall \leq y \leq \yrwall,
\end{equation}
see figure \ref{fig-vcm-flow-partition}.
\end{mydefinition}
We also introduce the critical area defined as follows.
%
\begin{mydefinition}[Critical Area]\label{vertbackdefcriticalArea}
	The critical area, $A_c(x)$ is the wetted cross sectional area of an exactly filled cross section. That is, the wetted cross         sectional area when the water level is exactly at the chosen  elevation, $\zwall$.
  It is defined by
	\begin{align}\label{vertbackeqncriticalarea}
     A_c(x)  &:=   \int_{y_l(x,\zwall)}^{y_r(x,\zwall)} \beta(\xy) dy.        
     \end{align}
\end{mydefinition}

\begin{myremark}
  If the channel is exactly full ($A=A_c$), then the water depth is exactly
  equal to the channel depth, $\beta$, and consequently
  %
\begin{equation}\label{vertnumeqnbetafromnotation}
  \Height(\xy;A) = \beta(\xy). 
\end{equation}
\end{myremark}  

\begin{mydefinition}[Channel Flow Lateral Boundaries]\label{vertbackdefywall}
Let $y_{l}^*(x, t)$ and $y_r^*(x, t)$  denote the $y$-coordinates of the 
left and right channel flow lateral boundaries, i.e.,
	\begin{align}\label{vertbackeqnywall}
	\begin{split}
		y_l^*(x,t) := \min \{ y : \eta(\xyt)>\zb, y \geq y_l(x,\zwall)  \},   \\
		y_r^*(x,t) := \max \{ y : \eta(\xyt)>\zb, y \leq y_r(x,\zwall)  \}.   
	\end{split}	
	\end{align}
Note that if there is water everywhere in the channel, then $y_{l,r}^*(x,t) = y_{l,r}^w(x)$. 	
\end{mydefinition}

\begin{myremark}[Channel Assumption]\label{vertchannelassumption}
 We assume that the channel never goes dry anywhere between the lateral wall boundaries,
  i.e., there are no islands in the channel:
   \begin{align*}
   	  \eta(\xyt) > \zb \quad \forall y \in \left( y_l^*(x,t), y_r^*(x,t) \right).
   \end{align*}
\end{myremark}

\begin{mydefinition}[Average Free Surface]
With the flow lateral boundaries known, we can define the laterally averaged  free surface elevation, namely
\begin{equation}\label{vertetaavg}
      \etabar = \frac{1}{  y_r^*(x,t)-y_l^*(x,t)  } \sint{y}{y_l^*(x,t)}{y_r^*(x,t)}{ \eTa}.
\end{equation}
\end{mydefinition}

Next, we define precisely what we mean by a full channel.
\begin{mydefinition}[Full Channel]\label{vertbackdeffullchannel}
We say a channel cross section at $x$ is full if
	\begin{equation}
		H(\xyt) \geq  \beta(\xy) \quad  \forall y \in \left( \ylwall, \yrwall) \right)~,
  \end{equation}
(recall \eqref{eq:waterHeight} for the definition of the water height and
definition \ref{vertbackdefbeta} for the channel depth $\beta$).

This means that if the channel is full the following statements hold:
\begin{itemize}
	\item given the total wetted area $A(x,t)$ in the cross section and the
        critical area $A_c$ as given in definition the following inequality
        holds (recall definition \ref{vertbackdefcriticalArea}):
          \begin{align}\label{vertbackeqntotalarea}
          	A(x,t) = \int_{\ylwall}^{\yrwall} H(\xyt) dy \geq \int_{\ylwall}^{\yrwall} \beta(\xy) dy = A_c(x),
          \end{align}
	\item 
    also $ \etabar \geq \zwall $ since by definition
				\begin{align*}
          \etabar
          = \avgsintb{y}{y_l^*(x,t)}{y_r^*(x,t)}{ H(\xyt) + \zb }
          \geq \avgsintb{y}{y_l^*(x,t)}{y_r^*(x,t)}{\beta(\xy) + \zb}
					 =  \avgsint{y}{y_l^*(x,t)}{y_r^*(x,t)}{ \zwall }
					 = \zwall.
        \end{align*}
\end{itemize}

\end{mydefinition}

\subsubsection{Vertically Partitioned Channel Flow}

The main idea of the VCM is to partition the channel flow into two layers,
based on the chosen channel top elevation, $\zwall$.
Since we want to approximate the flow with a 1D model if the channel is not full, we make
the following assumption to allow consistency with a 1D formulations.
\begin{mydefinition}[1D Consistency Assumption]\label{vertbackdefconsistencyrqd}
	If the channel is not full, then the lateral variation in free surface elevation, $\eTa$
	is negligible and the free surface can be taken to be its lateral average, $\etabar$.
\end{mydefinition}
\begin{myremark}
	This simply means that we assume the free surface to be laterally flat.
	This is exactly the 1D assumption and allows to apply 1D modelling
  whenever the channel is not full (a similar assumption is made for
  example for the 1D models in \cite{moralesetal2013,bladeetal2012,cungeetal80}).
\end{myremark}
%

Next, we define the following elevation:
\begin{mydefinition}[Time Dependent Interface]
The time dependent interface, $\etaone$ is the elevation defined as, 
\begin{equation}\label{vertbackeqneta1}
	\eta_1(x,t) = \min( \bar{\eta}(x,t), \zwall ).
\end{equation}
\end{mydefinition}
\begin{figure}[ht!]
	\subfigure[
             Full cross section; $\eta_1(x,t)$ equal to the channel wall elevation $\zwall$
             and lower layer water depth, $h_1(x,y,t)$
 	        equal to the the channel depth $\beta(x,y)$. Free surface $\eta(x,y,t)$ remains non flat laterally.
 	        Channel flow boundaries $y_{l,r}(x, \eta_1)$ equal to the channel top lateral walls, $y_{l,r}^w(x)$. ]
            {\label{fig-layers-full-case}\includegraphics[width=0.49\textwidth]{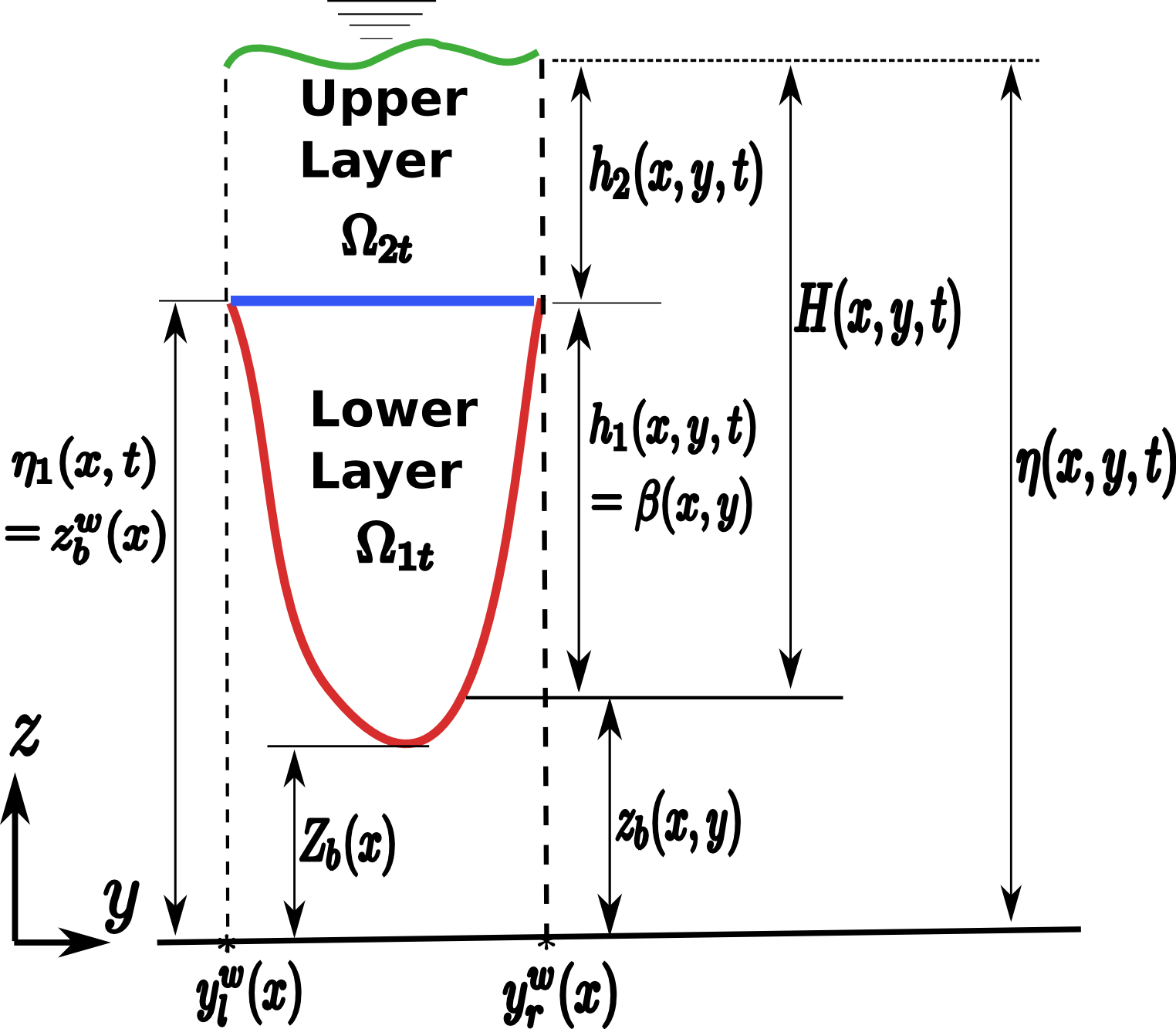} }
	\hspace{0.02\textwidth}
	\subfigure[
             Non full cross section; laterally flat free surface $\bar{\eta}(x,t)$,
             $\eta_1(x,t)$ equal to the laterally flat free surface $\bar{\eta}(x,t)$, which is less than the channel wall elevation $\zwall$,   $h_1(x,y,t)$ equal to the total water depth, $H(x,y,t)$ which is less than the 
            channel depth $\beta(x,y)$. The channel flow lateral interval, $|y_{r}(x, \eta_1) - y_{l}(x, \eta_1)|$ 
            is less than the total channel lateral interval $|y_{r}^w(x)- y_{l}^w(x)|$. ]
      {\label{fig-layers-nonfull-case}\includegraphics[width=0.48\textwidth]{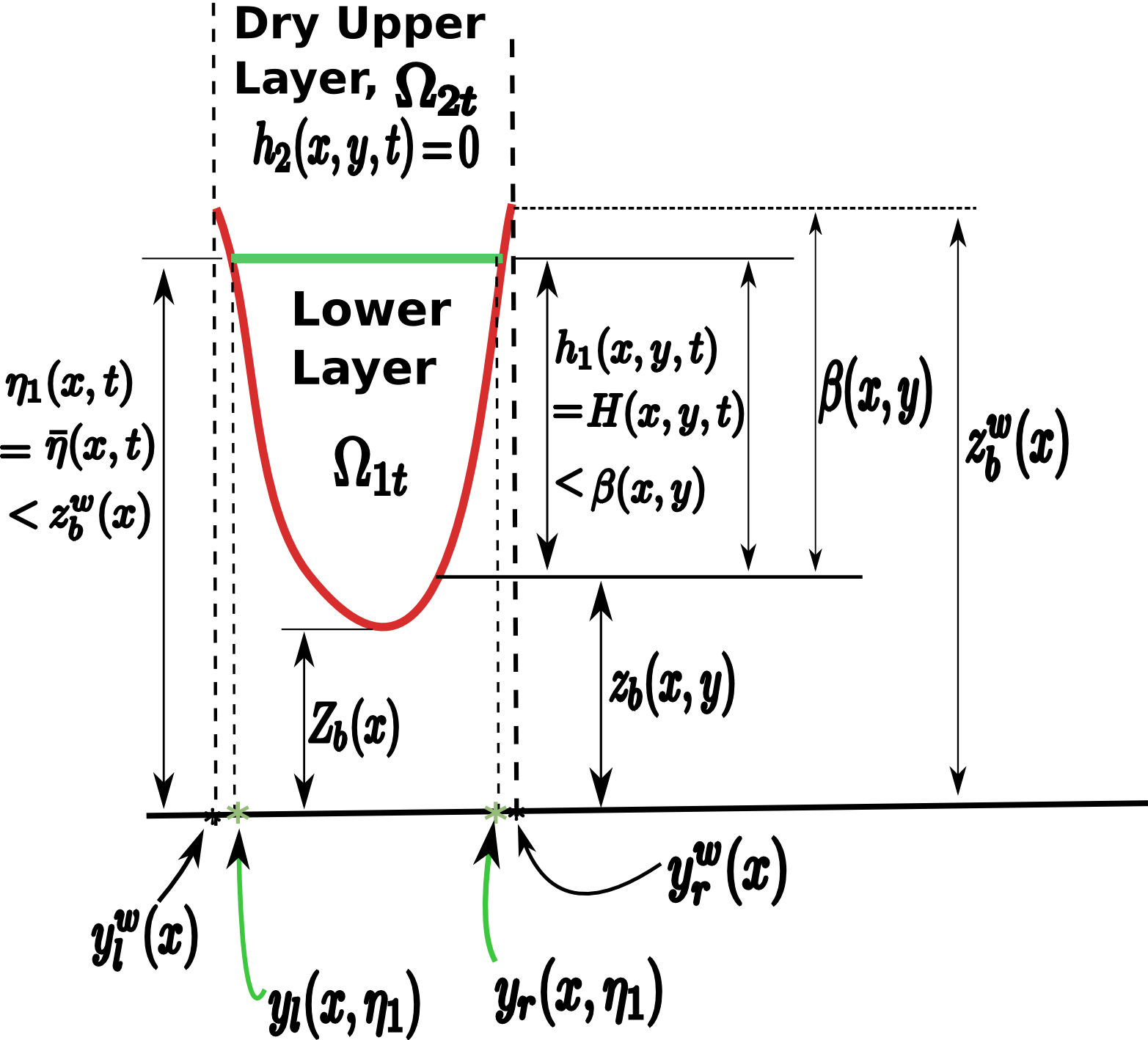} }
	\caption{The two layers in the VCM in the case of full channel (left) and the non full case (right)}
\end{figure}

\begin{myproposition}
  $y_{l,r}^*(x,t)=y_{l,r}(x, \eta_1(x,t) )$.
\end{myproposition}
\begin{proof}
Case 1 : If the channel is full (see figure \ref{fig-layers-full-case}), then
	by definition \ref{vertbackdeffullchannel} we have $\etabar \geq \zwall$.
  Hence $\eta_1 = \min(\etabar, \zwall) = \zwall$.
	Furthermore, by definition \ref{vertbackdefywall}, the lateral
	walls are $y_{l,r}^*(x, t) = y_{l,r}(x,\zwall)=y_{l,r}(x,\eta_1 )$.

Case 2 : If the channel is not full (see figure \ref{fig-layers-nonfull-case}), then
    $\etabar \leq \zwall, \mbox{ so }  \eta_1:=\min(\etabar, \zwall) = \etabar$.
    And definition \ref{vertbackdefconsistencyrqd} requires free surface to be constant and
    equal to the average, $\etabar$, hence the free surface at the lateral walls is $\etabar$. So, the lateral boundaries are
    $y_{l,r}^*(x,t) = y_{l,r}(x,\bar{\eta} ) = y_{l,r}(x,\etaone)$.
\end{proof}

We can now identify the channel flow domain, $\Omega_{ct}$  at time, $t$ as: 
\begin{equation}
	\Omega_{ct} = \{  (\xyz) \in \real^3: y_l(x, \eta_1(x,t) ) \leq y \leq  y_r( x, \eta_1(x,t) ), \zb \leq z \leq \eTa    \}.
\end{equation}

\begin{mydefinition}[Flow Partitions]\label{vertsecdefinelayers}
Given the total water depth, $H(\xyt)$, we partition the flow into the two
layers
\begin{align}
&	\Omega_{1t} = \{  (\xyz) \in \real^3 : \yleo \leq y \leq \yreo, \zb \leq z \leq \etaone          \},
	\\
&	\Omega_{2t} = \{ (\xyz) \in \real^3 :y_l(x, \eta_1(x,t) ) \leq y \leq  y_r( x, \eta_1(x,t) ), \etaone \leq z \leq \eTa      \}.
\end{align}
see figure \ref{fig-vcm-flow-partition}.
The water depths in the lower and upper layers are given by
\begin{align}
    h_1(\xyt) & = min( H(\xyt), \beta(\xy) ),  \quad  y_l(x, \eta_1(x,t) ) \leq y \leq  y_r( x, \eta_1(x,t)),   	\label{vertbackgeqnh1}
    \\
    h_2(\xyt) &= H(\xyt) - h_1(\xyt), \quad   y_l(x, \eta_1(x,t) ) \leq y \leq  y_r( x, \eta_1(x,t) ), \label{vertbackgeqnh2}
\end{align}
\end{mydefinition}
Note that $h_1(\xyt)+\zb = \etaone$ defines the top of the lower layer and
so by our 1D consistency assumptions $h_1$ is independent of $y$.



\subsubsection{Fundamental Equations:}\label{vertsecfundamentaleqns}
Under the assumption of hydrostatic pressure, the flow in the domain, $\Omega_{ct}$ is
governed by the incompressible free-surface Euler equations:
\begin{align}
	\begin{split} \label{vertcont}
	  \pdiff{x}{u(\xyzt)} + \pdiff{y}{v(\xyzt)} + \pdiff{z}{w(\xyzt)}  = 0,
	\end{split}
	\\
	\begin{split} \label{vertxmom}
	    \pdiff{t}{u(\xyzt)} 
	  +  \pdiffb{x}{u^2(\xyzt)} 
	  +  \pdiffb{y}{u(\xyzt)v(\xyzt)} 
	  \\
	  +  \pdiffb{z}{u(\xyzt)w(\xyzt)} 
	  =  -g\pdiff{x}{\eTa},
	\end{split}
	\\
	\begin{split} \label{vertymom}
	    \pdiff{t}{v(\xyzt)} 
	  +  \pdiffb{x}{u(\xyzt)v(\xyzt)} 
	  +  \pdiffb{y}{v^2(\xyzt)} 
	  \\
	  +  \pdiffb{z}{v(\xyzt)w(\xyzt)} 
	  =  -g\pdiff{y}{\eTa},
	\end{split}
\end{align}
where $(u,v,w)^T$ is the fluid velocity vector at point
$(\xyz)\in\Omega_{ct}$ at time $t$.
Furthermore, the following boundary conditions hold
\begin{flushleft}
\begin{equation}
	\left( 
		   u(\xyzt) \pdiff{x}{z_b(\xy)} 	+ v(\xyzt) \pdiff{y}{z_b(\xy)}	- 	   
		 w(\xyzt) 
	\right)\bigg|_{z=z_b(\xy)}	
	= 0,  \label{vertkinbed}
\end{equation}
\begin{equation}
	\left( 
		\pdiff{t}{\eTa} + u(\xyzt) \pdiff{x}{\eTa} 	+ v(\xyzt) \pdiff{y}{\eTa}	- w(\xyzt)
	\right)\bigg|_{z=\eTa}	 
	= 0,   \label{vertkinsurf}
\end{equation}
\end{flushleft}
%

\subsection{The Lower Layer Flow Model}\label{vertseclower}
The lower layer flow is assumed to be always one dimensional and thus the lateral variations
in the free-surface elevation do not have any impact on the lower layer flow.
Therefore, using the averaged free-surface, $\bar{\eta}(x,t)$, instead of
$\eta(\xyt)$. We thus replace the horizontal moment equation
\eqref{vertxmom} and the kinematic boundary condition \eqref{vertkinsurf} by
\begin{align}
&	    \pdiff{t}{u(\xyzt)} 
	  +  \pdiffb{x}{u^2(\xyzt)} 
	  +  \pdiffb{y}{u(\xyzt)v(\xyzt)} 
	  +  \pdiffb{z}{u(\xyzt)w(\xyzt)} 
	  =  -g\pdifft{x}{ \etabar },	 \label{vert1dxmom} \\
&	\left( 
		\pdifft{t}{\etabar} + u(\xyzt) \pdifft{x}{\etabar} 	- w(\xyzt)
	\right)\bigg|_{z=\etabar}
	= 0.
  \label{vertsvmkinsurf}
\end{align}
%
We define the area of wetted cross section, $A_1(x,t)$, the section averaged volumetric discharge,
$Q_1(\xt) $ and  section averaged velocity, $\uoned_1(\xt)$, for the lower layer, as follows:
\begin{align}
	Q_1(\xt) &= \dint{z}{y}{\yleo}{\yreo}{\zb}{\etaone}{ u(\xyzt)}, \label{vertQ1}
	\\
	A_1(\xt) &= \dint{z}{y}{\yleo}{\yreo}{\zb}{\etaone}{} = \sint{y}{\yleo}{\yreo}{h_1(\xyt)}, \label{vertA_1}
	\\
	\uoned_1(\xt) &= \frac{Q_1(\xt)}{A_1(\xt)}
		         = \frac{1}{A_1(\xt)}\dint{z}{y}{\yleo}{\yreo}{\zb}{\etaone}{ u(\xyzt)}\label{vertu1}.
\end{align}

Before we proceed, let us state the following important relations which are
easily proven (see \cite{chineduthesis} for details):
\begin{mylemma}\label{vertlowerlemmamin}
Let $A_c$ and $A$ be defined in \eqref{vertbackeqncriticalarea} and \eqref{vertbackeqntotalarea} respectively, and
 \begin{equation}\label{vcm-eqn-ac-star}
					 A_c^*(x,t) = \int_{\yleo}^{\yreo}\beta(\xy) dy,
	  \end{equation}
then
\begin{equation}	\min(A(x,t), A_c^*(x,t)) = \min(A(x,t), A_c(x) ). \end{equation}
\end{mylemma}
\begin{mytheorem}
Let  $A_c$, $A$,  $A_1$ and $A_c^*$ be defined in \eqref{vertbackeqncriticalarea}, \eqref{vertbackeqntotalarea},
\eqref{vertA_1} and \eqref{vcm-eqn-ac-star} respectively, then
	\begin{equation}\label{vertlowerlayereqnA1minAAc}
		A_1(x,t) = \min( A(x,t), A_c(x)).
	\end{equation}
\end{mytheorem}

\subsubsection{Mass Conservation Equation for Lower Layer:}

First we derive the model equations for the lower layer.
Integrating the equation \eqref{vertcont} over the lower layer cross section,
and
applying the Leibniz rule using the kinematic
boundary condition \eqref{vertkinbed}, we obtain
\begin{align}\label{vertlowermasseq}
	\begin{split}
	\pdifft{t}{A_1(\xt)} &+ \pdifft{x}{Q_1(\xt)} 
		=     -\sint{y}{\yleo}{\yreo}{ S(\xyt)  }, 
	 \end{split}
\end{align} 
where
\begin{equation}\label{vertS}
	S(\xyt) = \left[    w(\xyzt) - u(\xyzt) \pdifft{x}{\etaone} - \pdifft{t}{\etaone}  \right]\bigg|_{z=\etaone}, 
\end{equation}
\begin{myremark}
  Note that $ \sint{y}{\yleo}{\yreo}{ S(\xyt)  }  $
  is the averaged mass exchange term between the two layers. If the channel is not
  full, i.e., $\etaone=\bar{\eta}(\xt)$ then $S(\xyt)=0$ using the
  boundary condition \eqref{vertsvmkinsurf}
\end{myremark}

\subsubsection{Momentum Equation for Lower Layer:}
To derive the momentum equation for the lower layer, we integrate equation \eqref{vert1dxmom} over
the lower layer cross section
and using standard arguments arrive at
\begin{multline}\label{vertQequation}
	\pdifft{t}{Q_1(\xt)} + \pdiffb{x}{ \frac{Q_1^2(\xt) }{A_1(\xt)}  }
	 =
	   -gA_1(\xt)\pdifft{x}{\etabar} 
	  	- \sintb{y}{\yleo}{\yreo}
	  	  {
            u(\xyzt)|_{z=\etaone}S(\xyt)	  	  
	  	  },
\end{multline}
where the last integral on the right is the momentum exchange term between the two layers.


\subsection{Upper Layer Flow Model}\label{vertsecupper}
The upper layer flow is allowed to remain fully two-dimensional, so the Free-Surface Euler Equations,
\eqref{vertcont} - \eqref{vertymom} and \eqref{vertkinsurf} are applicable. However, the kinematic boundary 
condition on the bottom does not apply here because the bottom of the upper layer is $\etaone$ which is not a 
physical boundary that fluid particles cannot cross. Define the following quantities:
\begin{align}
  q_{2x}(\xyt) = \sint{z}{\etaone}{\eTa}{u(\xyzt)}, 
  \quad
  q_{2y}(\xyt) = \sint{z}{\etaone}{\eTa}{v(\xyzt)}.
\end{align}
So that the velocities are 
\begin{equation}
	u_2(\xyt) = \frac{q_{2x}(\xyt)}{h_2(\xyt)} \quad \mbox{ and }
	v_2(\xyt) = \frac{q_{2y}(\xyt)}{h_2(\xyt)},  
\end{equation}
where $\vec{q}_2(\xyt) = ( q_{2x}(\xyt), q_{2y}(\xyt) )^T $ is the upper layer $2D$ discharge vector and
$\vec{u}_2(\xyt) = ( u_{2}(\xyt), v_{2}(\xyt) )^T$ is the velocity vector
in the upper layer.

In the following, we derive the equations for the $2D$ quantities.
Integrating equation \eqref{vertcont} vertically over the upper layer,
we have
\begin{equation}\label{verth2equation}
	\pdifft{t}{h_2(\xyt)} + \pdifft{x}{q_{2x}(\xyt)} + \pdifft{y}{q_{2y}(\xyt)} = S(\xyt).
\end{equation}

Also, integrating equation \eqref{vertxmom} vertically over $\etaone \leq z \leq \eTa  $, applying the kinematic 
boundary condition (equation \eqref{vertkinsurf}) and simplifying, we have
\begin{multline}\label{vertq2xequation}
	\pdifft{t}{ q_{2x}(\xyt)  }
	+ \pdiffb{x}{  \frac{q^2_{2x}(\xyt)}{h_2(\xyt)}    + g/2h^2_2(\xyt)  } 
	+ \pdiffb{y}{  \frac{q_{2x}(\xyt)q_{2y}(\xyt)}{h_2(\xyt)}    }
	\\
	=
	  -gh_2(\xyt)\pdifft{x}{\zb}
	  -gh_2(\xyt)\pdifft{x}{h_1(\xyt)} 
	  + u(\xyzt)|_{z=\etaone}S(\xyt).
\end{multline}

Similarly, integrating equation \eqref{vertymom} vertically over $\etaone \leq z \leq \eTa  $, applying the kinematic 
boundary condition (equation \eqref{vertkinsurf}) and simplifying, we have
\begin{multline}\label{vertq2yequation}
	\pdifft{t}{ q_{2y}(\xyt)  }
	+ \pdiffb{x}{  \frac{q_{2x}(\xyt)q_{2y}(\xyt)}{h_2(\xyt)}    }
	+ \pdiffb{y}{  \frac{q^2_{2y}(\xyt)}{h_2(\xyt)}    + g/2h^2_2(\xyt)  } 
	\\
	=
	  -gh_2(\xyt)\pdifft{y}{\zb} 
	  + v(\xyzt)|_{z=\etaone}S(\xyt).
\end{multline}
The mass exchange $S(\xyt)$ between the layers is defined in equation \eqref{vertS}.

%

\subsection{Summary of Coupled Channel Flow Models}\label{vertsecchanelmodelsummary}
The two layer channel model derived so far consisting of the 1D
lower layer model
\begin{align}\label{vertmodeqnlayer1}
\begin{split}
\partial_t A_1 + \partial_x Q_1 &=  -\int_{\yleo}^{\yreo} S dy, \\
\partial_tQ_1 + \partial_x Q_1^2/A_1 & = - gA_1\partial_x\bar{\eta} - \int_{\yleo}^{\yreo} u_{\eta_1} Sdy,
\end{split}
\end{align}
and the 2D upper layer model,
\begin{align}\label{vertmodeqnlayer2}
\begin{split}
 \partial_th_2 + \nabla \cdot \vec{q}_2 &=  S,  \\
 \partial_t \vec{q}_2 + \nabla \cdot F^q(h_2, \vec{q}_2) &= -gh_2\nabla(z_b + h_1)  +  \vec{u}_{\eta_1}S,
\end{split}
\end{align}
%
where the fluxes are given by\\
$F^q(h_2, \vec{q}_2) = (F^x, F^y), F^x = (\frac{q_{2x}^2}{h_2} + \frac{g}{2}h^2_2, \frac{q_{2x} q_{2y}}{h_2}  )^T$,
$F^y= (\frac{q_{2x} q_{2y}}{h_2} , \frac{q_{2y}^2}{h_2} + \frac{g}{2}h^2_2,  )^T$,
$\vec{u}_{\eta_1}=(u(\xyt),v(\xyt))^T\bigg|_{z=\eta_1}$,
$h_1(\xyt) = \Height(\xy;A_1(x,t))$, $\eta_1 = h_1 + z_b $
and
$\etabar = \avgsint{y}{y_l^*(x,t)}{y_r^*(x,t)}{ \eTa}$.
%
%

%
%

Note that the models \eqref{vertmodeqnlayer1} and \eqref{vertmodeqnlayer2} are not closed
since the exchange term $S$ and interface velocity $\vec{u}_{\eta_1}$ are not known.
Following a similar idea as used in \cite{audusse2011multilayer} we will
solve the system using a two step approach where in the first step we solve
the equations without the exchange term
\begin{align}
\partial_t A_1 + \partial_x Q_1 = 0, \quad
\partial_tQ_1 + \partial_x \bigg( Q_1^2/A_1 \bigg)  = - gA_1\partial_x\bar{\eta} \label{vertnumeqn1dmodop1} 
\end{align}
and
\begin{align}
 \partial_th_2 + \nabla . \vec{q}_2 =  0 , \quad
 \partial_t \vec{q}_2 + \nabla .F^q(h_2, \vec{q}_2) = -gh_2\nabla(z_b + h_1) \label{vertnumeqn2dmodop1}
\end{align}
and then use this intermediate
solution to approximate the mass and momentum exchange between the layers.
But even then there is a difficulty in solving the 1D lower layer model since
the free-surface elevation term, $\bar{\eta}(x,t)$,  appearing in
\eqref{vertnumeqn1dmodop1} is different from  the actual top level, $\eta_1(x,t)$
of the lower layer flow, which is represented by $A_1$. This makes
directly deriving a well-balance scheme for the lower layer
model challenging and more importantly we cannot directly reuse an existing
1D channel solver for the solution. 
Since this is one of our major goals, we base our numerical scheme on a
reformulation of the above coupled channel model:

Integrating \eqref{vertmodeqnlayer2} over the cross section $y_l,y_r$ we arrive
at
\begin{align}\label{vertnumeqnA2upper}
\partial_t \begin{pmatrix}  A_2 \\ Q_2   \end{pmatrix}  =
  -\partial_x \begin{pmatrix} Q_2 \\ \frac{Q^2_2}{A_2}  \end{pmatrix}
    - \begin{pmatrix} -\bar{S} \\ gA_2 \partial_x
    \bar{\eta}-{u_{\eta_1}}\bar{S}\end{pmatrix}  +  \int_{\yleo}^{\yreo} u_{\eta_1} Sdy + \Phi(x,t), 
\end{align}
where
\begin{align}
\begin{split}
	&A_2(x,t) := \int_{y_l(x, \etaone}^{y_r(x, \etaone)}
	\int_{\etaone}^{\eTa} \, dz \, dy,
	\quad
	Q_2(x,t) := \int_{y_l(x, \etaone}^{y_r(x, \etaone)}
	 \int_{\etaone}^{\eTa} u(\xyzt) \, dz \, dy,
    \\
	&\Phi(x,t) := \bigg( \Phi^F(h_2, \myvec{q}_2) \bigg)_{y=y_l^w(x)}
				- \bigg( \Phi^F(h_2, \myvec{q}_2) \bigg)_{y=y_r^w(x)},
	\quad
	\Phi^F(h_2, \myvec{q}_2)  := \begin{pmatrix}
	          \myvec{q}_2 \cdot \myvec{n}  \\
	          \bigg(
	          		 \frac{q_{2x}^2}{h_2} + \frac{g}{2}h_2^2,
	          		 \frac{q_{2y}^2}{h_2}
	          \bigg)^T \cdot \myvec{n}
	\end{pmatrix}.
\end{split}	
\end{align}
assuming that the average horizontal velocities in both layers are very similar, i.e.,
$( \frac{Q_1}{A_1} - \frac{Q_2}{A_2})^2 \approx 0$ we can add the
equations for $A_1,Q_1$ and $A_2,Q_2$ together arriving at a model for the
1D full channel
\begin{align}\label{vertnumeqnAeqn}
\partial_t \begin{pmatrix} A \\ Q \end{pmatrix} 
					= -\partial_x \begin{pmatrix} Q \\ \frac{Q^2}{A}  \end{pmatrix}
                      - \begin{pmatrix} 0 \\ gA \partial_x \bar{\eta}	\end{pmatrix}  + \Phi(x,t).
\end{align}
So with this further approximation we end up with the following systems to
model the flow within the channel:
\begin{align}\label{verteqnAQfinal}
\begin{split}
  \partial_t A + \partial_x Q &= \Phi_A(x,t) \\
  \partial_t Q + \partial_x \frac{Q^2}{A} &=
         - gA \partial_x \bar{\eta} + \Phi_Q(x,t),
\end{split}
\end{align}
\begin{align}\begin{split}\label{verteqnh2q2final}                
  \partial_th_2 + \nabla . \vec{q}_2 &=  -S~, \\
  \partial_t \vec{q}_2 + \nabla .F^q(h_2, \vec{q}_2) &=
      -gh_2\nabla(z_b + h_1) - u_{\eta_1} S,
\end{split}      
\end{align}
where $h_1(\xyt) := \min( H(\xyt), \beta)$ and $\Phi = ( \Phi_A, \Phi_Q )^T$.
\begin{myremark}
The system \eqref{verteqnAQfinal} is a standard 1D approximation to the
channel flow with $\Phi$ providing the exchange terms between the channel
and the floodplain along the horizontal boundary. Thus a standard model can
be used to approximate \eqref{verteqnAQfinal}. The advantage of our
approach is that the second set of equations \eqref{verteqnh2q2final}
provides additional 2D information in the channel that can be used to
accurately compute the exchange terms $\Phi$.
\end{myremark}
\subsection{Floodplain Flow Model}\label{vertsec2dswes}  
We describe the flow in the floodplains using the 2D Shallow water equations, namely
\begin{align}\label{fv2dswewtfriceqnmodel}
  \partial_t \Pi + \nabla \cdot (F_1(\Pi),F_2(\Pi) =  S(\Pi, z_b ), 
\end{align}
where
\begin{align}\label{verteqnfloodmodel-fluxes-and-sources}
\begin{split}
&	\Pi = \begin{pmatrix}
		H \\ q_x \\ q_y
	\end{pmatrix}, 
\quad
%
  F_1(\Pi) = \begin{pmatrix}
   		q_x  \\ \frac{q_x^2}{H} + \frac{1}{2}gH^2 \\  \frac{q_xq_y}{H}
   \end{pmatrix}, \quad
   F_2(\Pi) = \begin{pmatrix}
   		q_y \\  \frac{q_xq_y}{H} \\ \frac{q_y^2}{H} + \frac{1}{2}gH^2
   \end{pmatrix},
\quad 
  %
  S(\Pi,z_b) = \begin{pmatrix}
   		                    0  \\  -gH\partial_x z_b(\xy) \\  -gH\partial_y z_b(\xy)
            \end{pmatrix}~.
\end{split}
\end{align}
%

\section{Numerical Algorithm for the Coupled Channel Flow Models}\label{vsecnumerical}

%
In this section, we describe in detail the numerical schemes and algorithms to solve
the two-layer models, \eqref{verteqnAQfinal} and \eqref{verteqnh2q2final}
for the channel flow. The 1D model \eqref{verteqnAQfinal} is solved on a 1D (quasi 2D) channel mesh (see figure \ref{vertnumfig1dmesh}), while the 2D upper layer model \eqref{verteqnh2q2final} is solved on a 2D channel mesh 
(see figure \ref{vertnumfig2dmesh}). The results are effectively combined to
derive the update values for the channel flow.
\subsection{Background}\label{vertnumsecbackground}
\begin{figure}[ht!]
  \begin{center}
  \subfigure[Top view of 1D (quasi 2D) channel mesh made of cells $K_i$. The 1D model is solved on each cell, $K_i$.]
   {\label{vertnumfig1dmesh} \includegraphics[width=0.44\textwidth]{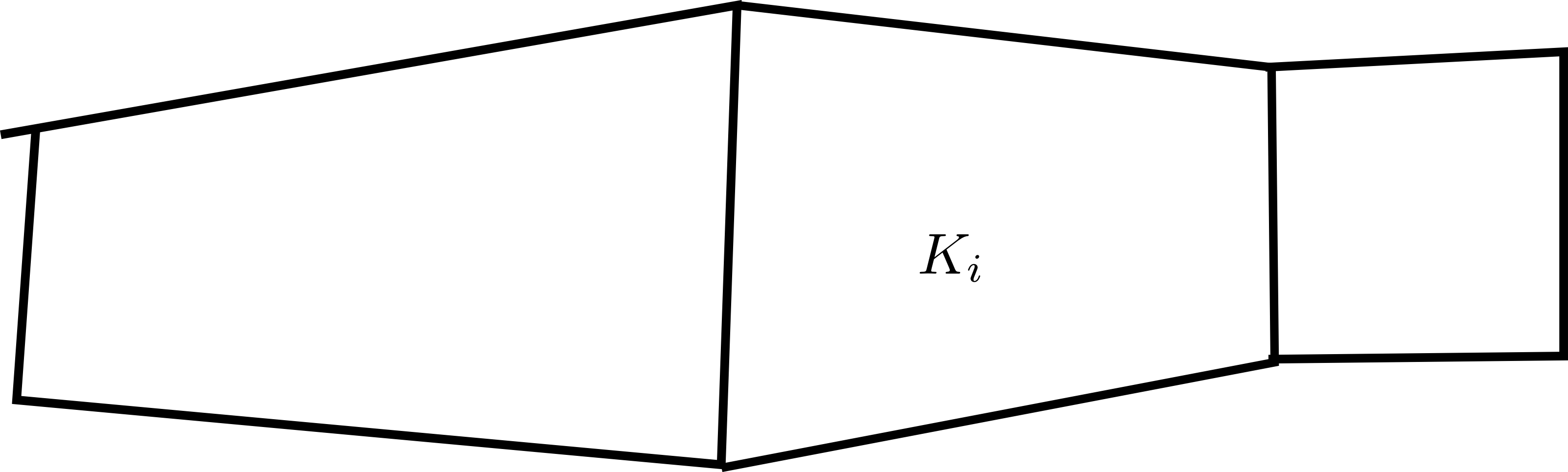}}
  \hspace*{0.03\textwidth}
   \subfigure[Top view of 2D mesh for the channel region only - no floodplains. Each 1D channel cell $K_i$ (in red) is further discretized to form a number, $N_y^1$ of 2D cells, $T_{ij}, j=1,2,..., N_y^1$.]
              {\label{vertnumfig2dmesh} \includegraphics[width=0.48\textwidth]{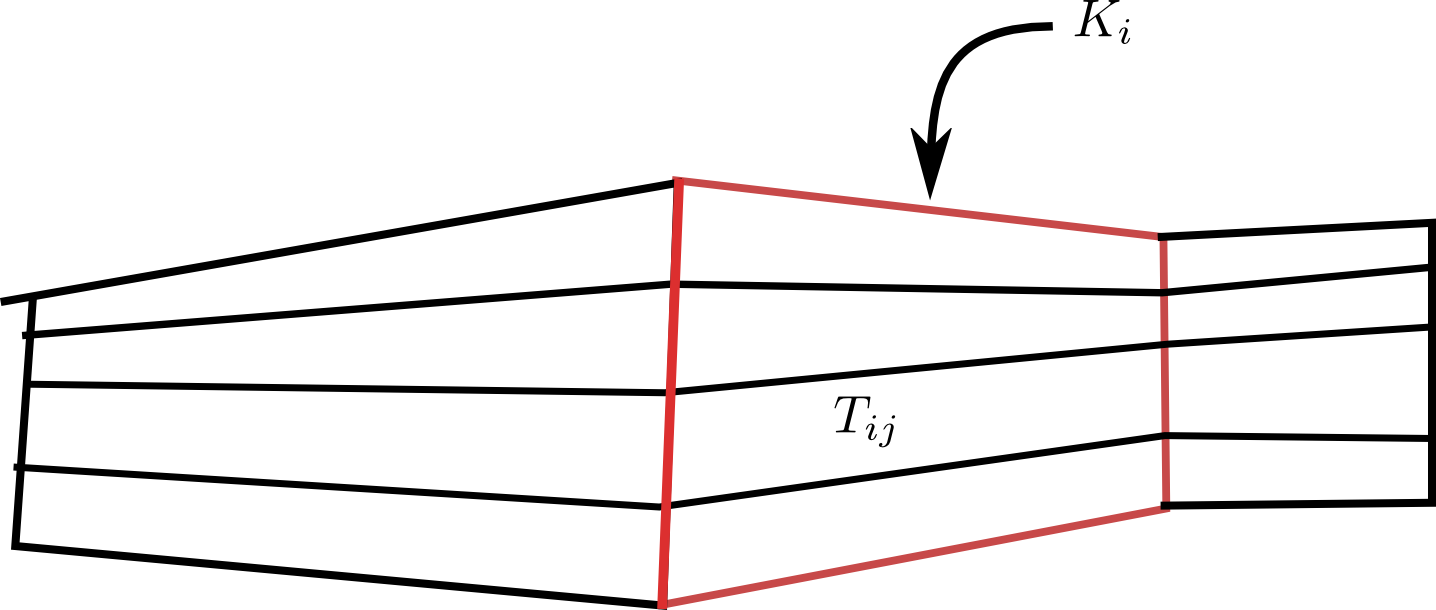}}   
  \end{center}
  \caption{Mesh for the two sub-models.}
\end{figure}
%

\begin{figure}[ht]
  \begin{center}
   \subfigure[1D channel cross section, $K_i$ depicting the discrete channel wall elevation, $\eta_i^\beta$, top width, $B_i$,
              depth $\beta_i$, critical area, $A_{c,i}$ and 1D bottom elevation $Z_{b,i}$.]
   {\label{vertnumfigcrosssec1d} \includegraphics[width=0.25\textwidth]{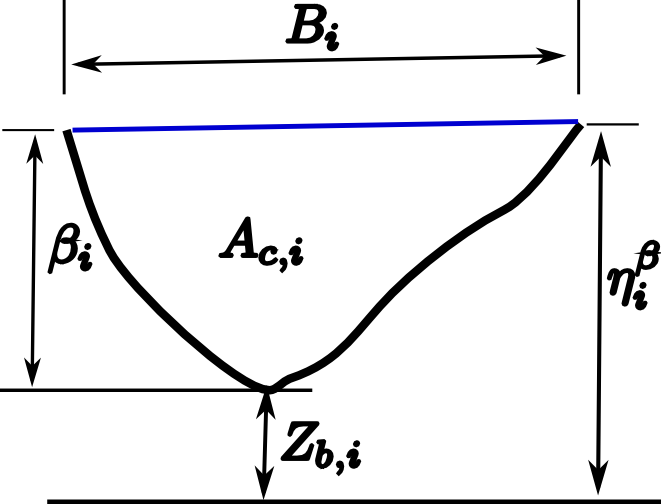}}
   \hspace*{0.02\textwidth}
   \subfigure[Full channel cross section viewed in $2D$ channel mesh. The laterally varying free-surface elevation, $\eta_{i,j}$ (in green)
   channel bottom elevation $z_{bij}$, channel depth $\beta_{ij}$ and water depth $H_{i,j}$ in the 2D channel cell $T_{ij} \in \Omega_h^2.$]
   {\label{vertnumfigcrosssec1d2dfull}
   \includegraphics[width=0.3\textwidth]{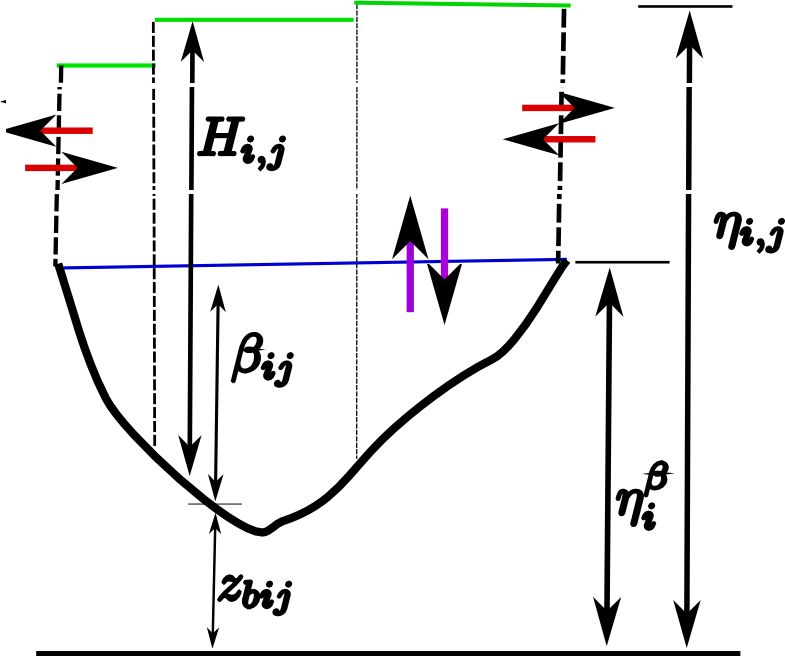}}
    \hspace*{0.02\textwidth}
   \subfigure[Non-full channel cross section viewed in $2D$ channel mesh cells. The laterally flat free-surface elevation, $\bar{\eta}_i$
   which is less than the channel elevation $\eta_i^\beta$.]
   {\label{vertnumfigcrosssec1d2dnotfull}
   \includegraphics[width=0.3\textwidth]{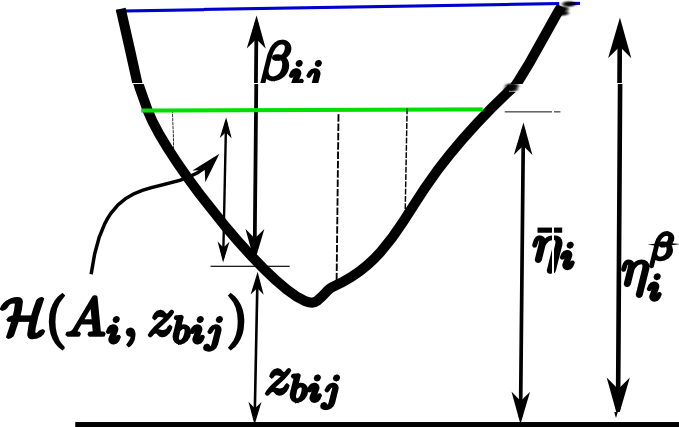}}
  \end{center}
  \caption{Discrete channel geometry in 1D channel cell (figure \ref{vertnumfigcrosssec1d}) and the
  	     cross sectional view of the channel flow  2D channel cells when full (figure \ref{vertnumfigcrosssec1d2dfull})
  and non-full (figure \ref{vertnumfigcrosssec1d2dnotfull}).}
\end{figure}

Let $\Omega_h^1$ be a 1D (quasi 2D) channel mesh with cells $K_i =
(x_{i-1/2}, x_{i+1/2} )$, see figure \ref{vertnumfig1dmesh}.
To consider the upper layer channel model, \eqref{vertmodeqnlayer2}, we further discretize each 1D channel cell $K_i$ into 2D cells,  $T_{ij}, j=1,...,N_y^1$, see figure \ref{vertnumfig2dmesh}.
This forms the 2D channel mesh $\Omega^2_h$, so that
\begin{align}\label{vertnumeqncells}
	\Omega_h^2 = \{  T_{ij} :  T_{ij} \subset K_i, j=1,...N_y^1; i=1,...,N_{1Dcells} \}.    
\end{align}
%
%
For the $1D$ cell, $K_i$, we define the following discrete channel quantities: 
\begin{align}\label{vertnumeqnTxdiscrete}
\begin{split}
	\eta_i^{\beta} = z_b^w(x_i), \quad	B_i = B(x_i, \eta_i^{\beta}), \quad	A_{c,i} = A_c(x_i), \quad
	Z_{bi} = Z_{b}(x_i), \quad \beta_i = \eta^\beta_i - Z_{b,i}
\end{split}	
\end{align}
as the channel wall elevation, channel top width, critical area, 1D bed
elevation, and channel depth, respectively, evaluated
at the center $x_i$ of each 1D cell, $K_i$ (see figure \ref{vertnumfigcrosssec1d}).
We also define the discrete quantities
\begin{align}
\begin{split}
	A_i^n       \approx  \frac{1}{\Delta x_i} \int_{x_i-1/2}^{x_i+1/2} A(x,t^n) dx, 
	\quad
	Q_i^n   \approx  \frac{1}{\Delta x_i} \int_{x_i-1/2}^{x_i+1/2} Q(x,t^n) dx
\end{split}
\end{align}
denoting approximations to the average total wetted cross section and section-averaged discharge in $K_i$.

For the associated $2D$ cells, $T_{ij}$ in $K_i$, we define two similar
sets of discrete quantities
\begin{align}
\begin{split}
	z_{bij}    = z_b(\vec{X}_{ij}), \quad
   \beta_{ij}  = \eta_i^{\beta} - z_{bij}
\end{split}
\end{align}
being the bed elevation and channel depth at the cell center, $ \vec{X}_{ij} $ of  $T_{ij}$,
see  figure \ref{vertnumfigcrosssec1d2dfull}.
Note the identity,
$Z_{bi} + \beta_i = z_{bij} + \beta_{ij} = \eta^\beta_i \quad \forall j=1,...,N_y^1$.
Furthermore, let
\begin{align}
\begin{split}
	h_{2,i,j}^n   \approx \frac{1}{|T_ij|} \int_{T_{ij}} h_2(\xy, t^n) dx dy, \quad
	\vec{q}_{2,i,j}^n \approx \frac{1}{|T_ij|} \int_{T_{ij}} \vec{q}_2(\xy, t^n) dx dy,
\end{split}
\end{align}
where $\vec{q}_2 = ( q_{2x} , q_{2y} )^T $ and $|T_{ij}|$ is the size(area) of $T_{ij}$.
%
%
%
Finally we also define full 2D cell averages given by
\begin{equation}
	H_{i,j}^n   \approx \frac{1}{|T_{ij}|} \int_{T_{ij}} H(\xy, t^n) dx dy, 
	\quad
	\vec{q}_{i,j}^n \approx \frac{1}{|T_{ij}|} \int_{T_{ij}} \vec{q}_(\xy, t^n) dx dy 
\end{equation}
which are cell averages for the full 2D data (sum of lower and upper layer) in 2D cells, $T_{ij}\in \Omega_h^2$.

%
%
%
%
Given a lateral distribution $\eta_{i,j}^n$ of the free-surface elevation in the 2D cells
$T_{ij}\subset K_i$ for $j=1,\dots,N_y^1$ we assume that we can
compute the wetted cross sectional area, $A_i^n$ in $K_i$ by
\begin{equation}
  A_i^n = \WArea( \eta_{i,1}^n, \eta_{i,2}^n,..., \eta_{i,N_y^1}^n)
\end{equation}
This requires information about the channel geometry,
for instance, if the channel has a rectangular cross-section,  with a laterally constant bottom 
elevation, $z_{bij}=Z_{bi}$ for
all $j$, then
$\WArea$ is given by
\begin{equation}\label{vertnumAcalpatotalarea}
 	\WArea( \eta_{i,1}^n, \eta_{i,2}^n,..., \eta_{i,N_y^1}^n ) 
			 := \sum_{j=1}^{N_y^1} ( \eta_{i,j}^n - Z_{bi} ) \Delta y_{ij}
			 = \sum_j^{N_y^1} H_{i,j}^n \Delta y_{ij}.
\end{equation}

%


Next, we formulate a discrete 1D consistency assumption
similar to definition \ref{vertbackdefconsistencyrqd}:
\begin{mydefinition}[Discrete Consistency requirement]\label{vertnumdefdiscreteconsistencyrqd}
 The solution,
 $(H, q_x, q_y)^n_{i,j}$ at  $t^n$  satisfies a discrete consistency requirement, if the following condition holds.
If $ \exists j^* \in \{ 1,2,..., N_y^1 \}$
  such that $\eta_{i,j^*}^n<\eta_i^{\beta}$, then
 \begin{align}\label{vertnumeqndiscreteconsistencyrqd}
	 \eta_{i,j}^n = \bar{\eta}^n_i < \eta_i^{\beta}   \quad  \forall j \in \{ 1,2,..., N_y^1 \}.
 \end{align}
\end{mydefinition}

We will later show that if the initial data satisfies this discrete consistency requirement, then 
the requirement is satisfied for all time by our scheme.

\subsection{Solution of the two-layer Channel Models}\label{vertnumsecoperator}
We now solve the two layer channel models,  \eqref{verteqnAQfinal},  \eqref{verteqnh2q2final} to evolve
the full 2D data, $(H, \myvec{q})_{ij}$ from the initial to the final time.
To this end we first initialize the value $(A_i^0,Q_i^0)_i$ for the full 1D
model and the values $(h_{2,i,j}i^0,\myvec{q}_{2,i,j}^0)_{ij}$ for the
2D upper layer model (in the region where the channel is full). To evolve
these values to the next time level we first update the 2D upper layer
values ignoring the exchange term between the two layers but taking the
horizontal fluxes with the floodplain into account. Next we update the values of the full 1D model,
using the channel-floodplain lateral fluxes computed by the 2D upper layer solver
to completely solve \eqref{verteqnAQfinal}. Finally, we correct the
solution in the upper layer by estimating the exchange term between the two layers.
Note that the first step only leads to computational cost in regions where the channel is
full (i.e. $h_{2,i,j}^n>0$). Otherwise this step can be skipped and only
the 1D step has to be carried out. Exchange terms with the floodplain are
also correctly taken into account in the first step.


\subsubsection{Initializing the two-layer model from given full 2D data}\label{vertnumsecdistri}
Given full 2D data, $(H, \myvec{q})^0_{i,j}$ at the initial time and $z_{bij}$ in each 2D channel cell
$T_{ij} \in \Omega_h^2$, we initialize the 1D model, \eqref{verteqnAQfinal}
and the 2D upper layer model, \eqref{verteqnh2q2final} as follows:
\begin{align}
  A_i^0  &= \WArea( \eta_{i,1}^0, \eta_{i,2}^0,..., \eta_{i,N_y^1}^0 ),
  \qquad
  Q_i^0  = \sum_{j=1}^{N_y^1} q_{x,ij}^0 \Delta y_{ij}.
\intertext{and}
       h_{1,i,j}^0  &:= \min(H_{i,j}^0, \beta_{ij} ), \qquad
           h^0_{2,i,j}  :=  H_{i,j}^0 - h_{1,i,j},  \qquad
   \myvec{q}_{2,i,j}^0  :=  h_{2,i,j}^0 \myvec{u}^0_{i,j},
  \label{vertnumeqndistriupper}
\end{align}
where $\myvec{u}_{i,j}^0 = \frac{ \myvec{q}_{i,j}^0}{H^0_{i,j}}$  and 
 $\Delta y_{i,j}$ is the average lateral width of $T_{i,j}$.

%
\subsection{Step 1: Solution of Upper Layer Model without Exchange Terms}\label{vertnumsec2dsolver}
We now propose a method to evolve the solution for the model, \eqref{vertnumeqn2dmodop1}.
Let us emphasis that any appropriate method for a 2D shallow water model can be used
here since \eqref{vertnumeqn2dmodop1} can be interpreted as a standard 2D
shallow water equation with topography, $\zb + h_1(\xyt)$ referred to as the apparent
topography in \citep{bouchutbook2004, rotating}. In our simulation we use an
apparent topography hydrostatic reconstruction scheme
\citep{audusseetal2004, bouchutbook2004, rotating}:

Let $T_j$ be a 2D cell in a 2D channel mesh and $T_k$ be its neighbour.
Let $e_{jk}$ be the edge between $T_j$ and $T_k$, and $\vec{n}_{jk}$ be
the unit vector normal outwards to $T_j$. 
Furthermore, let $|T_j|$ and $|e_{jk}|$ be the area of $T_j$ and length of $e_{jk}$ respectively and
let $\mathcal{E}_j$ be the set of all edges of $T_j$. 
Let $h_{1,j}^n, \beta_j, z_{b,j}$ and $\Pi_{2,j}^n = ( h_{2,j}^n, \myvec{q}_{2,j}^n )^T $ be the lower layer water depth,
channel depth, bottom elevation and upper layer cell average vector in 2D channel cell, $T_j$; while  
$ h_{1,k}^n, \beta_k, z_{b,k} $ and $ \Pi^n_{2,k} = ( h_{2,k}^n, \myvec{q}_{2,k}^n )^T$ are those in 2D cell, $T_k$. 
Note that the neighbour cell, $ T_k $ could actually be in the floodplain; the only requirement is that
the current cell, in which we compute the upper layer solution, must be in the channel.

Define the apparent topography:
\begin{align}
		\eta_{1,p}^n = h_{1,p}^n + z_{b,p}, \quad \mbox{ for } p = j,k.
\end{align}
Then reconstruct the apparent topography in the hydrostatic fashion:
\begin{align}
		\eta_{1,jk}^{n*} :=  \max( \eta_{1,j}^n, 	\eta_{1,k}^n ).
\end{align}
Next, define
\begin{align}\label{vert-num-inter-pi2}
	\tilde{h}_{2,p}^n :=  \max \bigg( \eta_{1,p}^n  + h_{2,p}^n - \eta_{1,jk}^{n*},  0  \bigg), \quad 
	\widetilde{ T_{\vec{n}_{jk}} \Pi }_{2,p} :=  \frac{ \tilde{h}_{2,p}^n }{h_{2,p}^n} T_{\vec{n}_{jk}} \Pi_{2,p}^n, 
	\quad \mbox{ for } p = j,k.
\end{align}

Then, the apparent topography hydrostatic reconstruction scheme \citep{bouchutbook2004, rotating} for the model reads:
\begin{align}\label{vertnumeqn2dsolver}
\begin{split}
\Pi_{2,j}^{n+1*}  =   \Pi_{2,j}^n - \frac{\Delta t}{|T_j|}  \sum_{e_{jk}\in \mathcal{E}_j} |e_{jk}|  \bigg(
								T_{\vec{n}_{jk}}^{-1} \phi(  \widetilde{T_{\vec{n}_{jk}}\Pi_{2,j}^n},
										\widetilde{T_{\vec{n}_{jk}}\Pi_{2,k}^n} )
							+  
								T_{\vec{n}_{jk}}^{-1}S^{hrm}( h_{2,j}^n, \tilde{h}_{2,j}^n) 
						\bigg)						
   + \Delta t S_b(\Pi_{2,j}^n),
\end{split}  
\end{align}	
where $S_b(.)$ is the friction term defined in \eqref{verteqnfloodmodel-fluxes-and-sources}, and
\begin{align}
	 S^{hrm}( a, b)  :=  \begin{pmatrix}   0 \\ \frac{g}{2}( a^2 - b^2   )  \\ 0  \end{pmatrix}, 
	 \quad
\phi(\Pi_L, \Pi_R ) 
	= \begin{cases} F_1(\Pi_L),   & \mbox{ if } s_L \geq 0,\\
	               F_1^* := \frac{s_RF_1(\Pi_L) -  s_LF_1(\Pi_R) + 	s_Ls_R(\Pi_R-\Pi_L) }{s_R-s_L}, & \mbox{ if } s_L \leq 0 \leq s_R, \\
	               F_1(\Pi_R),  & \mbox{ if }  s_R \leq 0.
	\end{cases}
\end{align}

\subsection{Step 2: Complete Solution of the Full 1D Model} \label{vertnumsecfull1dsolver}
We now solve the full 1D model, \eqref{verteqnAQfinal}. This involves two sub-steps, namely
\subsubsection{Step 2.1: Black-Box Stage}
Here, we use any available 1D solver to solve the full 1D model without the lateral flux term,
$\Phi$, namely
\begin{align}\label{verteqnAQ-no-phi}
\begin{split}
  \partial_t A + \partial_x Q = \Phi_A(x,t), \quad
  \partial_t Q + \partial_x \frac{Q^2}{A} =
         - gA \partial_x \bar{\eta} + \Phi_Q(x,t),
\end{split}         
\end{align}
with the data, $(A, Q)_i^n$.
In this paper, we adopt the scheme from \citep{moralesetal2013, moralesetallargedt}, see also \cite{chineduthesis, chineduandreas1}.
This gives the approximation, $( A_i^{\widetilde{n+1}}, Q_i^{\widetilde{n+1}})$.
\subsubsection{Step 2.1: Add the floodplain coupling term}
The discrete coupling term, $\Phi_i$ is the vector consisting of the first 
two components of the 2D lateral fluxes already computed by the upper layer 
solver in section \ref{vertnumsec2dsolver}. Therefore, the final update value for the full 1D model, \eqref{verteqnAQfinal} is then given by 
\begin{equation}
	(A_i^{n+1}, Q_i^{n+1} ) = ( A_i^{\widetilde{n+1}}, Q_i^{\widetilde{n+1}})
		                      + \Phi_i \Delta t, 
\end{equation}
where $\Delta t$ is the time step.

\subsection{Step 3: Final Update of the Upper Layer Model}
With the intermediate solution, $( h_{2}, \vec{q}_{2} )_{i,j}^{n+1*}$
known, we now completely solve the upper layer model including the exchange term,  \eqref{verteqnh2q2final}.
The approximate solution of \eqref{verteqnh2q2final} is the approximate solution of the system
\begin{align}
\begin{split}\label{vertnumeqn2dmodop2}
  \partial_th_{2,i,j} &=   S_{i,j},  \\
  \partial_t \vec{q}_{2,i,j}  &=  \vec{u}_{\eta_1,i,j} S_{i,j},
\end{split}  
\end{align}
with the initial data, $(h_{2}, \vec{q}_{2})_{i,j}^{n+1*}$, 
where
\begin{align*}
 S_{i,j} &\approx  \frac{1}{|T_{ij}|} \int_{T_{ij}} S dx \, dy, \quad
 \vec{u}_{\eta_{1,i,j}} \approx \frac{1}{T_{ij}} \int_{T_{ij}} \vec{u}_{\eta_1} \, dx \, dy.
\end{align*}

Then using forward Euler time discretization,
the approximate solution of \eqref{vertnumeqn2dmodop2}
is
\begin{align}
\begin{split}\label{vertnumeqnhq2new}
  h_{2,i,j}^{n+1}  &= h_{2,i,j}^{n+1*} + S_{i,j} \Delta t, \\
  \vec{q}_{2,i,j}^{n+1} &= \vec{q}_{2,i,j}^{n+1*} +  \vec{u}_{\eta_1,i,j}S_{i,j} \Delta t. 
\end{split}
\end{align}
%
%
To complete the description of the scheme, we need to define the exchange
terms and interface velocities $S_{i,j}$ and $\vec{u}_{\eta_1,i,j}$.
But there are no equations for these terms.
However they can be determined by requiring that the following conditions be satisfied:
(i) The operation from the intermediate solutions
$(n+1*)$ to the new update solutions $(n+1)$ must locally (and globally) conserve both mass and momentum.
(ii) The final update values must satisfy the discrete consistency requirement (definition \ref{vertnumdefdiscreteconsistencyrqd}).
(iii) the new update values must satisfy the non-negativity of water heights.
These conditions allow us to first obtain the heights $h_{2,i,j}^{n+1}$, then the exchange terms, $S_{i,j}=\frac{h_{2,i,j}^{n+1} - h_{2,i,j}^{n+1*}}{\Delta t}$, and 
finally the interface velocities $\vec{u}_{\eta_1,i,j}$ are calculated
so that we can compute $\vec{q}_{2,i,j}^{n+1}$.
%

\subsubsection{Step 3.1: Approximating the Lower Layer Flow}
By \eqref{vertlowerlayereqnA1minAAc} we have
\begin{equation}\label{eqn-new-A1}
	A_{1,i}^{n+1}  =  \min( A_i^{n+1}, A_{c,i} ).
\end{equation}
given $A^{n+1}_i$ computed above.
Note that we can not directly apply $A_{1,i}^{n+1} = A_i^{n+1}-\sum_j h_{2,ij}^{n+1} \Delta y$ since we do not yet know $h_{2,ij}^{n+1}$.
However, by mass/momentum conservation, the following must hold
\begin{align}\label{eqn-mass-mom-conserv} 
	\begin{pmatrix} A_{i} \\ Q_{i}   \end{pmatrix}^{n+1}
					 = \begin{pmatrix}  A_{1,i} \\ Q_{1,i}
   \end{pmatrix}^{n+1*} 
	+
	\sum_j \begin{pmatrix}
		 h_{2,i,j}^{n+1*} \\  q_{2x,i,j}^{n+1*}
   \end{pmatrix}\Delta y_{i,j}.
\end{align}
Hence the intermediate lower layer wetted area is given by
\begin{align}\label{eqnLowerLayerIntermSolver}
  \begin{pmatrix}  A_{1,i} \\ Q_{1,i}
   \end{pmatrix}^{n+1*} 
   =
	\begin{pmatrix} A_{i} \\ Q_{i}   \end{pmatrix}^{n+1}   
	-
	\sum_j \begin{pmatrix}
		 h_{2,i,j}^{n+1*} \\  q_{2x,i,j}^{n+1*}
   \end{pmatrix}\Delta y_{i,j}.
\end{align}
Using this lower layer update, $A_{1,i}^n$ and the intermediate solutions,
$A_{1,i}^{n+1*}$ and $h_{2,ij}^{n+1*}$ we can now compute the upper layer
heights in the next step.

\subsubsection{Step 3.2: Upper Layer Heights}
\label{sec:upperlayerheight}
Using the definition of $A_{1,i}^{n+1}$ from \eqref{eqn-new-A1} we have to
distinguish two cases $A^{n+1}_{1,i}=A_i^{n+1}$ and
$A^{n+1}_{1,i}=A_{c,i}<A_i^{n+1}$:
%
\subsubsection*{ Case 1: If $A_{1,i}^{n+1} = A_i^{n+1}$
(Lower Layer not full at $t^{n+1}$)} 
Recall, $ A_i^{n+1} = A_{1,i}^{n+1} + \sum_j h_{2,i,j}^{n+1} \Delta y_{ij} $, so we have $\sum_j h_{2,i,j}^{n+1} \Delta y_{ij}= A_i^{n+1} - A_{i,i}^{n+1} = 0 $.
Since $h_{2,i,j}^{n+1} \geq0$, we must have
\begin{equation}\label{vertnumeqnh2newcase1}
   h_{2,i,j}^{n+1} = 0.
\end{equation}
%
%
\subsubsection*{Case 2: If $A_{1,i}^{n+1} = A_{c,i}$ (Lower Layer full at $t^{n+1}$)}
Denoting $A_{2,i}^{n+1*} := \sum_j h_{2,i,j}^{n+1*}\Delta y_{ij} $,
$ A_{2,i}^{n+1} := \sum_j h_{2,i,j}^{n+1}\Delta y_{ij}$ and replacing $A_i^{n+1}$
by $A_{1,i}^{n+1} + \sum_j h_{2,i,j}^{n+1} \Delta y_{ij}  =
A_{c,i} + \sum_j h_{2,i,j}^{n+1} \Delta y_{ij}$,   then we can write the first equation in \eqref{eqn-mass-mom-conserv} in the following form:
\begin{align}\label{vertnumeqnA2newcase2}
	A_{2,i}^{n+1}  = A_{2,i}^{n+1*} + \bigg( A_{1,i}^{n+1*} - A_{c,i}  \bigg),
\end{align}
and consider two further cases.
%
\subsubsection*{Case 2a: If $ A_{1,i}^{n+1*} - A_{c,i}  \geq 0 $ (Lower Layer full at intermediate state)}
Then $A_{2,i}^{n+1}  \geq A_{2,i}^{n+1*} $  by an amount $A_{excess,i}^{n+1*} := A_{1,i}^{n+1*} - A_{c,i} \geq 0 $,
so we add the constant excess height, $h_{excess,i}^{n+1*} =
A_{excess,i}^{n+1*}/B_i$ to the intermediate solution $h_{2,i,j}^{n+1*}$
uniformly over all 2D cells, namely
\begin{equation}\label{vertnumeqnh2newcase2a}
	 h_{2,i,j}^{n+1} = h_{2,i,j}^{n+1*} + \frac{A_{excess,i}^n}{B_i}.
\end{equation}
%
\subsubsection*{ Case 2b: If $  A_{1,i}^{n+1*} - A_{c,i}  < 0 $  (Lower Layer not full at intermediate state)}
Then  $A_{2,i}^{n+1}  < A_{2,i}^{n+1*} $ by the amount, $A_{gap,i}^{n+1*} = A_{c,i} - A_{1,i}^{n+1*}  > 0$.
Hence we remove $A_{gap,i}^{n+1*}$ from  $A_{2,i}^{n+1*}$ using the following algorithm.
\begin{itemize}
	\item initialize $h_{2,i,j}^{n+1} = h_{2,i,j}^{n+1*} $ for all $T_{ij}$.
	\item $h_{gap,i}^{n+1*} = \frac{A_{gap,i}^{n+1*}}{B_i}, \quad TOL=10^{-12}$.
	\item while( $h_{gap,i}^{n+1*} > TOL$ )
		\begin{enumerate}[i]
			\item for all $T_{ij}$
				\begin{enumerate}
					\item $h_t = h_{2,i,j}^{n+1}$,
					\item Reduce upper layer height by the gap height :
						 \begin{equation} h_{2,i,j}^{n+1} = max( 0, h_{2,i,j}^{n+1} - h_{gap,i}^{n+1*} ),  \label{vertnumeqnh2newcase2b}  \end{equation} 
					\item Remove area of reduced height from total gap area:\\
							$A_{gap,i}^{n+1*} = A_{gap,i}^{n+1*} - |  h_{2,i,j}^{n+1} - h_t |\Delta y_{i,j} $, 
				\end{enumerate}
			\item $h_{gap,i}^{n+1*} = \frac{A_{gap,i}^{n+1*}}{B_i}$.
    \end{enumerate}
\end{itemize}
Equations \eqref{vertnumeqnh2newcase1}, \eqref{vertnumeqnh2newcase2a} and
\eqref{vertnumeqnh2newcase2b} compute $h_{2,i,j}^{n+1}$ in $T_{ij} \subset K_i$ for all 
cells.

\subsubsection{Step 3.3: Upper Layer Discharge and Lower Layer Discharge}
Having computed $h_{2,i,j}^{n+1}$, we compute the exchange terms, $S_{i,j}$ 
using the first equation in \eqref{vertnumeqnhq2new}, hence
\begin{align}\label{vert-num-eqn-sij}
	S_{i,j}    = \frac{ h_{2,i,j}^{n+1} - h_{2,i,j}^{n+1*}  }{\Delta t}, 
\end{align}
and compute the interface velocity, $\vec{u}_{\eta_1,i,j}$ following
\citep{audusse2011multilayer} namely
\begin{align}\label{vert-num-eqn-inter-velo}
	\vec{u}_{\eta_1,i,j} = \begin{cases} \vec{u}^{n+1*}_{2,i,j}
        = \frac{\vec{q}_{2y,i,j}^{n+1*}}{h_{2,i,j}^{n+1*}}, & \mbox{ if }  S_{i,j} \leq 0, \\
        (u^{n+1*}_{1,i} , 0 )^T
        = \bigg( \frac{Q_{1,i}^{n+1*}}{A_{1,i}^{n+1*}}, 0 \bigg)^T,   & \mbox{ if } S_{i,j} > 0.
        \end{cases} 	                             
\end{align} 
Thus the interface velocity is equal to the upper layer velocity if mass is
flowing from the upper to the lower layer and otherwise equal to the lower
layer velocity.

Using the above definitions of the exchange terms and interface velocity,
we define the upper layer discharge $\myvec{q}_{2,i,j}^{n+1}$ by using the second equation in
\eqref{vertnumeqnhq2new}.

This completes the upper layer solution $(h_{2,ij}^{n+1},
\myvec{q}_{2,ij}^{n+1} )$.

Implicitly we can also compute the lower layer solution $(A_{1,i}^{n+1}, Q_{1,i}^{n+1})$
by using \eqref{eqn-new-A1}
and local momentum conservation
 \begin{align}\label{vertnumeqnMomConservQ}
	 Q_i^{n+1} = Q_{1,i}^{n+1} + \sum_j q_{2x,i,j}^{n+1}\Delta y_{ij}
\end{align}
which leads to the lower layer discharge given by
 \begin{align}\label{vertnumeqn1DQnew}
	 Q_{1,i}^{n+1} = Q_{i}^{n+1} - \sum_j q_{2x,i,j}^{n+1}\Delta y_{ij}
\end{align}

\subsection{Obtaining the full $2D$ data $H_{i,j}^{n+1}$,
$q_{x,i,j}^{n+1}$ and $q_{y,i,j}^{n+1}$ for $T_{i,j} \subset K_i$}\label{vertnumsecupdatefull2d}

This final step is only required for postprocessing and for studying the
properties of the scheme in the following section.

The full water depth, $H_{i,j}^{n+1}$ is given by
\begin{equation}\label{vertnumeqnHnew}
	H_{i,j}^{n+1} = h_{1,i,j}^{n+1} + h_{2,i,j}^{n+1},
\end{equation}
where
\begin{equation}\label{vertnumeqnh1new}
  h_{1,i,j}^{n+1} = \Height( \vec{X}_{ij}; A_{1,i}^{n+1} ).
\end{equation}
%
The full 2D  $x$-discharge is computed as
\begin{align}\label{vertnumeqnqxnew}
	q_{x,i,j}^{n+1} &:= q_{1x,i,j}^{n+1} + q_{2x,i,j}^{n+1} \notag \\
		            &= h_{1,i,j}^{n+1}\frac{Q_{1,i}^{n+1} }{ A_{1,i}^{n+1} } + q_{2x,i,j}^{n+1}.
\end{align}
%
Finally, the $y$-discharge is computed as follows :
\begin{align}\label{vertnumeqnqynew}
	q_{y,i,j}^{n+1}  = H^{n+1}_{i,j} v_{i,j}^{n+1},
\end{align}
where
\begin{align}
	 &v_{i,j}^{n+1} := v_{2,i,j}^{n+1}
   = \begin{cases} \frac{  q_{2y,i,j}^{n+1*} + v_{\eta_1,i,j}S_{i,j}\Delta t }{ h_{2,i,j}^{n+1}}, & \mbox{ if } h^{n+1}_{2,i,j} > 0, \\
	 	                                             0,  & \mbox{ else }.
	                                 \end{cases} \label{vertnumeqnvnew}
\end{align}

This completes the numerical algorithm for the channel flow,
for the floodplain flow we use a hydrostatic reconstruction
scheme as described for the upper layer model, see \citep{chineduthesis} for details.
At the channel/floodplain interface, the lateral fluxes are effortlessly
computed since in the 2D channel cell we have complete 2D value $( H_{i,j}, \vec{q}_{i,j})^{n}$
at least if the channel is not full. In this important case,
we have an approximation of the vertical discharge which is not the case for most existing methods
where only a 1D model is used in the channel which does not provide values
for the lateral discharge.
If the channel is not full but the
surrounding floodplain is not dry, the channel 1D flow data
$(A_i^n,Q_i^n)$ is used to compute the coupling fluxes.

\section{Properties of the Scheme}\label{vertnumsecprops}
We now consider some important properties of the scheme proposed above.
%
\begin{mydefinition}[Consistency of Distribution Operation]\label{vertnumrmkdistri}
We require a distribution operation such as the one in section \ref{vertnumsecdistri} to be
\textbf{consistent} in the sense  that whenever
the channel is not full ($H_{i,j} \leq \beta_{ij}$ for all $T_{ij} \subset K_i$),
then the following conditions must hold:
$ h_{2i,j},  q_{2x,i,j} $ and $ q_{2y,i,j} $ each equal zero $\forall
  T_{ij} \subset K_i;$.
\end{mydefinition}
This ensure that we do not solve a two-layer problem when the channel is not full.

\begin{mydefinition}[No-Numerical Flooding Properties]\label{vertnumrmkassemble}
	In order to satisfy the no-numerical flooding property, the full 2D data  $(H_{i,j}, q_{xi,j}, q_{yi,j})^{n}$
	must satisfy the following condition for all $n$:
	\begin{enumerate}[i]
		\item Either the lower layer is full or the upper layer is empty, i.e.,
        \begin{equation} \label{vertnumeqnnonumfloodi}
								( \beta_{ij} - h_{1i,j}^{n} )h_{2i,j}^{n} = 0 \quad
                \forall T_{ij} \subset K_i;
				\end{equation}
		\item If $h_{2,i,j}^{n}=0$, then
              (a) $q_{2x,i,j}^{n} = 0$
              (b) $u_{i,j}^{n} =  \frac{Q_{1,i}^{n}}{A_{1,i}^{n} }$  (which is constant laterally)
					     and (c) $q_{y,i,j}^{n}=0$; 
		\item If $ \exists j^*\in\{1, 2, ..., N_y^1 \}$ such that $\eta_{i,j^*}^{n} < \eta_i^{\beta}$, then
		         $ \eta_{i,j}^{n} = \eta_{i,j^*}^{n} \forall j  $.
	\end{enumerate}
\end{mydefinition}
Condition $(i)$ means that the lower layer is either full $(h_{1,i,j}^{n+1} = \beta_{i,j})$ or upper layer is dry $h_{2,i,j}^{n+1}=0$. In other
words, there should not be gap between the two layers. A consequence of
this property is also that no overflowing of the channel can occure unless
it is full.
Condition ($ii$) states that if the upper layer is dry,
the upper layer velocities and discharges must vanish, and the full layer flow velocity, $u$ must
be laterally uniform.
Condition $(iii)$ means that if the channel is not full,
then the free-surface must be flat, i.e., not vary in the lateral direction.
Thus the discrete consistency requirement (definition
\ref{vertnumdefdiscreteconsistencyrqd}) is satisfied at all time steps.

In the following we will show that if the solution at time $t^n$ satisfies
the above \emph{no numerical flooding properties} then this is also true at time $t^{n+1}$.
In addition, we also prove that the scheme is well-balanced and conserves mass
under suitable conditions on the intermediate solutions.

\begin{mytheorem}[Consistency of Distribution Operation]\label{vertnumtheoremdistri}
    The distribution operation proposed in section \ref{vertnumsecdistri} is consistent with the problem in the 
    sense of definition \ref{vertnumrmkdistri}.
\end{mytheorem}    
\begin{proof}
We need to prove that the ascertion in definition \ref{vertnumrmkdistri} is true.
Let the channel not be full, that is, $H_{i,j} \leq \beta_{ij}$. Then, we have
	$h_{2,i,j} = H_{i,j} - \min( H_{i,j}, \beta_{ij})  = H_{i,j}  - H_{i,j} = 0,$
	and   
	$ q_{2x,i,j} = q_{2y,i,j} = 0 \quad  $  (see \eqref{vertnumeqndistriupper}).
\end{proof}


\begin{mytheorem}[No-Numerical Flooding Property]\label{vertnumtheoremnonumflood}
    The vertical coupling scheme as derived in \eqref{vertnumeqnHnew}, \eqref{vertnumeqnqxnew} and
    \eqref{vertnumeqnqynew}, preserves the no-numerical flooding property in the sense of
    definition \ref{vertnumrmkassemble}.
\end{mytheorem}  
\begin{proof}
  \begin{enumerate}[i]
	\item We prove \eqref{vertnumeqnnonumfloodi} on case-by-case bases
   (see section \ref{sec:upperlayerheight}).\\
    Case 1:  $A_{1,i}^{n+1} = A^{n+1}_{i}$, then
		     $$ h_{2,i,j}^{n+1} = 0, \quad \mbox{ (see } \eqref{vertnumeqnh2newcase1}).$$
		Cases 2a and 2b : $A_{1,i}^{n+1} = A_{c,i} $, then
         $$ h_{1,i,j}^{n+1} := \Height(\xy_{i,j} ; A_{1,i}^{n+1}) = \Height( \xy ; A_{c,i}) = \beta_{ij}
		     \quad (\mbox{see } \eqref{vertnumeqnbetafromnotation}) $$
		Therefore, we have 
		$h_{2,i,j}^{n+1} = 0$ or  $ h_{1,i,j}^{n+1} = \beta_{ij} $ in either case. 
		Hence, $$( \beta_{ij} - h_{1i,j}^{n+1} )h_{2i,j}^{n+1} = 0 \quad
      \forall T_{ij} \subset K_i $$ as claimed.
	\item Let $h_{2,i,j}^{n+1} = 0$, then \\
		(a) \begin{align*} 
			& S_{i,j} = -\frac{ h_{2,i,j}^{n+1*}} {\Delta t} < 0  
			 \quad \mbox{ (see \eqref{vert-num-eqn-sij}) }  
			 \\
			  & => u_{\eta_1,i,j} = u_{2,i,j}^{n+1*} \quad \quad \mbox{ (see \eqref{vert-num-eqn-inter-velo}) } 
			  \\
              & => u_{\eta_1,i,j}S_{i,j}\Delta t = - q_{2x,i,j}^{n+1*} => q_{2x,i,j}^{n+1}   = 0.  \quad \quad
             \mbox{(see \eqref{vertnumeqnhq2new})}.
       \end{align*}
    (b) Since $ q_{2x,i,j}^{n+1}   = 0  $, then by \eqref{vertnumeqnqxnew}, we have
		    $$q_{x,i,j}^{n+1} = h_{1,i,j}^{n+1}\frac{Q_{1,i}^{n+1}}{A_{1,i}^{n+1} } 
        = \frac{Q_{1,i}^{n+1}}{A_{1,i}^{n+1} } H_{i,j}^{n+1} 
          \quad \mbox{ (because } h_{2,i,j}^{n+1}=0).$$
	    	Hence,
	    	 $$ u_{i,j}^{n+1} := \frac{q_{x,i,j}^{n+1}}{H_{i,j}^{n+1}} = \frac{Q_{1,i}^{n+1}}{A_{1,i}^{n+1} } \quad \mbox{ (which is constant in } j).$$
		(c) $$v_{i,j}^{n+1}=0 
		     \quad
		     \mbox{(see} \eqref{vertnumeqnvnew})
		     \quad
		      => q_{y,i,j}^{n+1} = 0.$$
	 \item   
	 	Assume that for $j^* \in \{1,2, ...,N_y^1 \}  $ we have $\eta^{n+1}_{i,j*} < \eta_i^{\beta}$,
      then this corresponds to case 1 in section \ref{sec:upperlayerheight} because cases 2a and 2b satisfy 
	 	$$A_{1,i}^{n+1}=A_{c,i} => \eta_{i,j}^{n+1} \geq \eta^\beta_i \forall j=1,...,N_y^1.$$
	 	Since, $\eta^{n+1}_{i,j*} < \eta_i^{\beta}$ corresponds to case 1, then it satisfies 
	 	$$A_{1,i}^{n+1}=A_i^{n+1} \leq A_{c,i}  => h_{2,ij}^{n+1} = 0 \forall j.$$
	 	So, $$\eta_{i,j}^{n+1} := H_{i,j}^{n+1} + z_{bij} =  h_{1,ij}^{n+1} +
	 	 z_{bij} =  \eta^{n+1}_{1,i}  \forall j=1,...,N_y^1.$$
	 	 That is $\eta_{i,j}^{n+1}$ is independent of $j$, hence also equal to
	 	 $\eta_{i,j^*}^{n+1}$.
	 	Therefore, $$\eta_{i,j}^{n+1} =  \eta^{n+1}_{1,i} = \eta_{i,j^*}^{n+1} \forall j ={1,..,N_y^1}.$$
	 	Hence the free surface is laterally flat. Therefore, the  discrete consistency requirement (definition \ref{vertnumdefdiscreteconsistencyrqd})
	 	holds at $t^{n+1}$, and the scheme satisfies the no-numerical flooding property.	
 \end{enumerate}
\end{proof}

\begin{mytheorem}[Well-balanced property]\label{vertnumthmwellbalance}
  If the numerical schemes used to compute the intermediate solutions, $(A_{1,i}^{n+1*}, Q_{1,i}^{n+1*})^T$ and
  $(h_{2,i,j}^{n+1*}, q_{2x,i,j}^{n+1*}, q_{2y,i,j}^{n+1*})^T$ are well-balanced,  then the vertical
  coupling method, \eqref{vertnumeqnHnew}, \eqref{vertnumeqnqxnew} and
    \eqref{vertnumeqnqynew}, is well-balanced. This is true for any kind of well-balance, not only 
    for lake at rest.
\end{mytheorem}

\begin{proof}
	Let the intermediate solutions be well-balanced, i.e.,
  \begin{alignat*}{2}
    (A_{1,i}^{n+1*}, Q_{1,i}^{n+1*})^T &=  (A_{1,i}^n, Q_{1,i}^n )^T,
    \quad&\quad
	  (h_{2,i,j}^{n+1*}, q_{2x,i,j}^{n+1*}, q_{2y,i,j}^{n+1*}    )^T 
    &= (h_{2,i,j}^n, q_{2x,i,j}^n, q_{2y,i,j}^n    )^T.
  \end{alignat*}
We need to show that $(H, \myvec{q})_{i,j}^{n+1} = (H, \myvec{q})_{i,j}^{n}$.\\	
\textbf{First}, we show that $ A^{n+1}_{1,i} = A^n_{1,i}$.
	Recall that by \eqref{eqn-mass-mom-conserv}
	 \begin{align*}
	      A_{1,i}^{n+1} &:= \min( A_i^{n+1}, A_{c,i}  )
	                      = \min( A_{1,i}^{n+1*}+A_{2,i}^{n+1*}, A_{c,i})
	                    = \min( A_{1,i}^n + A_{2,i}^n, A_{c,i})
                       = A^n_{i,1},
\end{align*}
Hence,
$
h_{1,i,j}^{n+1} 
  :=\Height( \xy ; A_{1,i}^{n+1} )
                   = \Height( \xy ; A_{1,i}^n ) 
                   =: h_{1,i,j}^n.                                
$\\
	 
\textbf{Secondly}, we show that $A^{n+1}_{2,i} = A^{n}_{2,i}$
using \eqref{eqn-mass-mom-conserv} and that
$A_{1,i}^{n+1*} = A_{1,i}^n$:
\begin{align*}
	 	A_{2,i}^{n+1} &= A_{i}^{n+1} - A_{1,i}^{n+1}
	 	               = (A_{1,i}^{n+1*} + A_{2,i}^{n+1*} ) - A_{1,i}^n
	 	              = A_{2,i}^{n+1*}
	 	              = A_{2,i}^n.
\end{align*}
Hence, 
$
	A_i^{n+1} = A_i^n.
$

\textbf{Thirdly}, we use the above results to show that $h^{n+1}_{2,i,j} = h^{n}_{2,i,j}$. 	
\\
Case 1: $ A_{1,i}^{n+1} = A_{i}^{n+1} < A_{c,i}$. We have already shown that
  $A_{i}^{n+1} = A_{i}^n$ so
	$A_{i}^n < A_{c,i} => h^{n}_{2,i,j} = 0 $ and $A_{i}^{n+1} < A_{c,i}
  => h^{n+1}_{2,i,j} = 0 $.
     Hence,  $h_{2,i,j}^{n+1} = 0 = h^{n}_{2,i,j}. $ in case 1.
\\
Cases 2a and 2b: $ A_{1,i}^{n+1} =  A_{c,i}$, that is
	 $A_{1,i}^{n+1} =A_{1,i}^{n} = A_{1,i}^{n+1*} = A_{c,i}$
   using $A_{1,i}^{n+1} =A_{1,i}^{n} = A_{1,i}^{n+1*} )$.
	 But $$A_{1,i}^{n+1*} = A_{c,i} => A_{excess,i}^{n+1*}=A_{gap}^{n+1*}=0 => h_{excess,i}^{n+1*} = h_{gap}^{n+1*}= 0. $$
	 Hence by \eqref{vertnumeqnh2newcase2a} and \eqref{vertnumeqnh2newcase2b}, we have
	  $$h^{n+1}_{2,i,j} = h^{n+1*}_{2,i,j} = h^{n}_{2,i,j}.$$

\textbf{Fourtly}, we show that $\myvec{q}_{2,i,j}^{n+1} = \myvec{q}_{2,i,j}^{n} $
and $Q_{1,i}^{n+1} = Q_{1,i}^n $. \\
Since $ h_{2,i,j}^{n+1} = h_{2,i,j}^{n+1*} $, then we have
using \eqref{vert-num-eqn-sij}, \eqref{vertnumeqnhq2new} that
$ S_{i,j} = 0$ and
$\myvec{q}_{2x,i,j}^{n+1} = \myvec{q}_{2x,i,j}^{n+1*} = \vec{q}_{2x,i,j}^n$.

Also,
\begin{align*}
	    Q_{1,i}^{n+1} &:=  Q_i^{n+1} - Q_{2,i}^{n+1} 
	                       \quad \quad (\mbox{see } \eqref{vertnumeqn1DQnew}) \\
	                &= Q_{1,i}^{n+1*} + Q_{2,i}^{n+1*} - Q_{2,i}^{n+1}
	                \quad \quad (\mbox{see } \eqref{eqn-mass-mom-conserv})  \\
	                &= Q_{1,i}^{n+1*}
	                 = Q_{1,i}^{n} 
	                 \quad \quad (\mbox{since } 
	                 Q_{1,i}^{n+1*} = Q_{1,i}^{n}
	                  \mbox{ and }
	                 Q_{2,i}^{n+1*} = Q_{2,i}^{n}).
\end{align*}	
 Finally, we use the above results to show that $(H, \vec{q})_{i,j}^{n+1} =  (H, \vec{q})_{i,j}^n$.    By \eqref{vertnumeqnHnew}-\eqref{vertnumeqnqynew}, we have     
 \begin{align*}	 	                         
  H_{i,j}^{n+1} &:= h_{1,i,j}^{n+1} + h_{2,i,j}^{n+1}
      = h_{1,i,j}^n+h_{2,i,j}^n =: H_{i,j}^{n}, 
      \\
 q_{x,i,j}^{n+1}    &:=  h_{1,i,j}^{n+1}\frac{Q_{1,i}^{n+1}}{A_{1,i}^{n+1}} 
 							+  q_{2x,i,j}^{n+1}
            =  h_{1,i,j}^n\frac{Q_{1,i}^n}{A_{1,i}^n} + 
	    q_{2x,i,j}^n =:  q_{x,i,j}^n,  
	    \\
 q_{y,i,j}^{n+1} &:= H_{i,j}^{n+1}\frac{q_{2y,i,j}^{n+1}}{h_{2,i,j}^{n+1} }
      					 = H_{i,j}^n\frac{q_{2y,i,j}^{n}}{h_{2,i,j}^{n} } 
      					 =:  q_{y,i,j}^n.
\end{align*}
Hence, we have shown that 
$(H, \vec{q})_{i,j}^{n+1} =  (H, \vec{q})_{i,j}^n$, so the method is well-balanced.
\end{proof}

Next, we prove that the proposed schemes for the intermediate solutions are mass conservative and also
well-balanced for lake at rest.
Since the hydrostatic reconstruction method is mass conservative, \citep{audusseMarie2005well, audusseetal2004}
and the $2D$ model, \eqref{vertnumeqn2dmodop1} does not introduce
any source term to the height equation, then the scheme \eqref{vertnumeqn2dsolver} is mass conservative.
We state and prove a theorem below to show that this scheme, like the standard hydrostatic reconstruction \citep{audusseetal2004} 
method, preserves well-balance of lake at rest.
\begin{mytheorem}
	The upper layer scheme, \eqref{vertnumeqn2dsolver} is well-balanced with respect to lake at rest.
\end{mytheorem}
\begin{proof}
 The proof follows the same lines given in \citep{audusseMarie2005well,
  audusseetal2004}, for details see \cite{chineduthesis}.
\end{proof}

\begin{mytheorem}\label{vertthm1dwellbal}
	If the underlying 1D solver for \eqref{vertnumeqnAeqn} is well-balanced then the intermediate lower layer scheme, 
\eqref{eqnLowerLayerIntermSolver} is also well-balanced.
\end{mytheorem}
\begin{proof}
	Let scheme for \eqref{vertnumeqnAeqn} be well-balanced, then $ (A_{i}^{n+1}, Q_{i}^{n+1})^T = (A_{i}^n, Q_{i}^n)^T $. Since
	the $2D$ solver, \eqref{vertnumeqn2dsolver}, is also well-balanced, 
	so $ (h_{2,i,j}^{n+1*}, \vec{q}_{2,i,j}^{n+1*})^T = ( h_{2,i,j}^n, \vec{q}_{2,i,j}^n )^T $. Therefore, the 
	lower layer scheme, \eqref{eqnLowerLayerIntermSolver} becomes
	\begin{align*} 
	      \begin{pmatrix} A_{1,i} \\ Q_{1,i}   \end{pmatrix}^{n+1*} 
	      &:= 
	      \begin{pmatrix} A_i \\ Q_i \end{pmatrix}^{n+1} 
	      -
	      \begin{pmatrix} A_{2,i} \\ Q_{2,i} \end{pmatrix}^{n+1*} 
	    =
	    \begin{pmatrix} A_i \\ Q_i \end{pmatrix}^{n} 
	    - \begin{pmatrix} A_{2,i} \\ Q_{2,i} \end{pmatrix}^{n}\\
	&= \begin{pmatrix} A_{1,i} \\ Q_{1,i} \end{pmatrix}^{n} 
    \end{align*}
	So, the $1D$ scheme is well-balanced.
\end{proof}

\begin{mytheorem}[Conservation]\label{vertnumthmmassconsrvation}
	If the intermediate solutions are mass conservative, then the vertical coupling solution is also
	mass conservative.
\end{mytheorem}

\begin{proof}
	Let the intermediate solutions be mass conservative, then
\begin{align}\label{vertnumeqnmassconsv}
			\sum_i A_{1,i}^{n+1*} = 
			\sum_i A_{1,i}^n \mbox{ and } \sum_i A_{2,i}^{n+1*} = \sum_i A_{2,i}^n.  
\end{align}
By \eqref{eqn-mass-mom-conserv}, we have
\begin{align*}
	A_{1,i}^{n+1} +A_{2,i}^{n+1} = A_i^{n+1} = 	A_{1,i}^{n+1*} +A_{2,i}^{n+1*}. 
\end{align*}
Hence, 
    \begin{align*}			
        \sum_i( A_{1,i}^{n+1} + A_{2,i}^{n+1} ) &= \sum_i( A_{1,i}^{n+1*} 
        									+ A_{2,i}^{n+1*} ) 
           = \sum_i ( A_{1,i}^{n} + A_{2,i}^{n} ) 
   	              	   \quad \quad \mbox{(by \eqref{vertnumeqnmassconsv} ).} \notag 
		\end{align*}
\end{proof}

\begin{mytheorem}\label{vertthm1dmassconsv}
	If the underlying solver for \eqref{vertnumeqnAeqn} is mass conservative, then the intermediate lower layer scheme, \eqref{eqnLowerLayerIntermSolver} is also mass-conservative.
\end{mytheorem}
\begin{proof}
	Let scheme for \eqref{vertnumeqnAeqn} be mass conservative, then $ \sum_i A_{i}^{n+1*} = \sum_i A_{i}^n $. Since
	the $2D$ solver, \eqref{vertnumeqn2dsolver}, is also mass conservative, 
	so $ \sum_i A_{2,i}^{n+1*} = \sum_i A_{2,i}^n $. Therefore,  
	\eqref{eqnLowerLayerIntermSolver} gives
	\begin{align*} 
	     \sum_i A_{1,i}^{n+1*} &:= \sum_i A_i^{n+1} - \sum_i A_{2,i}^{n+1*} \\
	     		               &=  \sum_i A_i^n - \sum_i A_{2,i}^n  \quad \mbox{(by conservation of $A_i^{n+1}$ and $A_{2,i}^{n+1*}$)}  \\
	     		               &=\sum_i \bigg( A_i^n - A_{2,i}^n \bigg)  =   \sum_i A_{1,i}^n   \quad  \mbox{(by definition)}.
    \end{align*}
	So, the $1D$ scheme is mass conservative.
\end{proof}

\begin{mytheorem}
  The vertical coupling method (VCM), described in \eqref{vertnumeqnHnew}- \eqref{vertnumeqnvnew} is well-balanced with respect to lake at rest
	and is mass conservative.
\end{mytheorem}
\begin{proof}
 Since the $1D$ solver \citep{moralesetal2013}  is well-balanced for lake at rest
 and mass conservative,
 then by theorems \ref{vertthm1dwellbal} and \ref{vertthm1dmassconsv}
 the same holds for the intermediate lower layer scheme.
 Since the intermediate upper layer scheme is also
 well-balanced and conservative
 the results follows using theorems
 \ref{vertnumthmwellbalance}, \ref{vertnumthmmassconsrvation}.
\end{proof}

This completes the theoretical aspect of the VCM. In the next section, we
present some numerical experiments to evaluate its performance compared to
other coupling methods.

\section{Numerical Results}\label{vsecnumResults}
We consider two simple test cases to verify the performance of the proposed method.
We use the full 2D simulation results
as the reference solution, and compare the results of the VCM with those of the Horizontal Coupling
Method (HCM) \cite{chineduandreas1} and the Flux-Based Method (FBM) \cite{bladeetal2012}.
Recall that in addition to the usual data required for the simulation,
e..g, domain size and intitial conditions, the VCM method also requires a
choice of the function $\zwall$ used to determin when a channel is
considered full and flodding might occure.

\subsection{Test Case 1 : Dam-Break Flow into a Flat Floodplain}\label{numsectest1}

\begin{figure}[ht!] 
  \subfigure[Top view of Channel and Floodplain for Test 1]{%
    \label{fvmnumfigwhycouplinggeometry}
	  \includegraphics[width=0.49\textwidth]{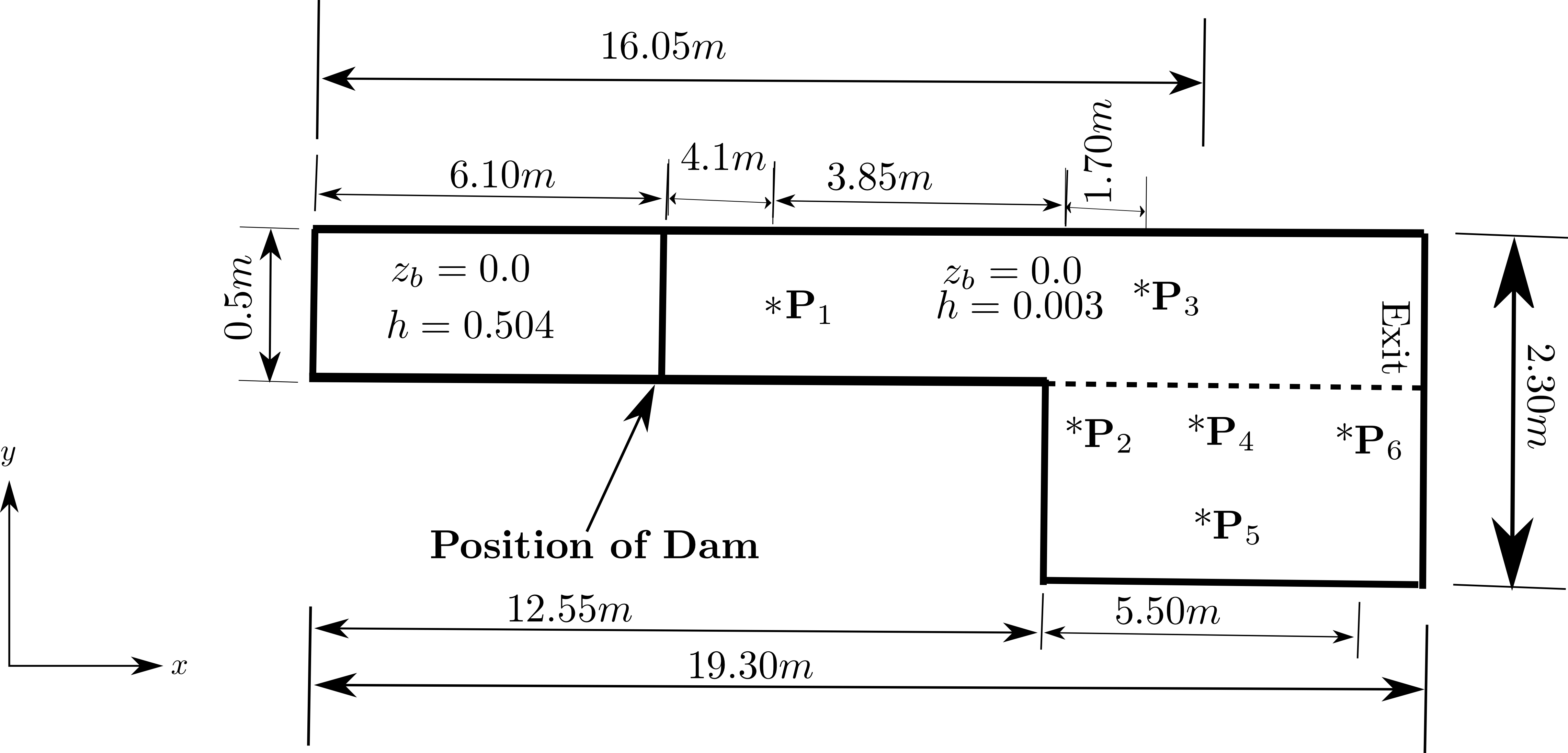}%
  }
  \subfigure[Top view of Channel and Floodplain for Test 2]{%
    \label{numfiggeometrytest2}%
	  \includegraphics[width=0.49\textwidth]{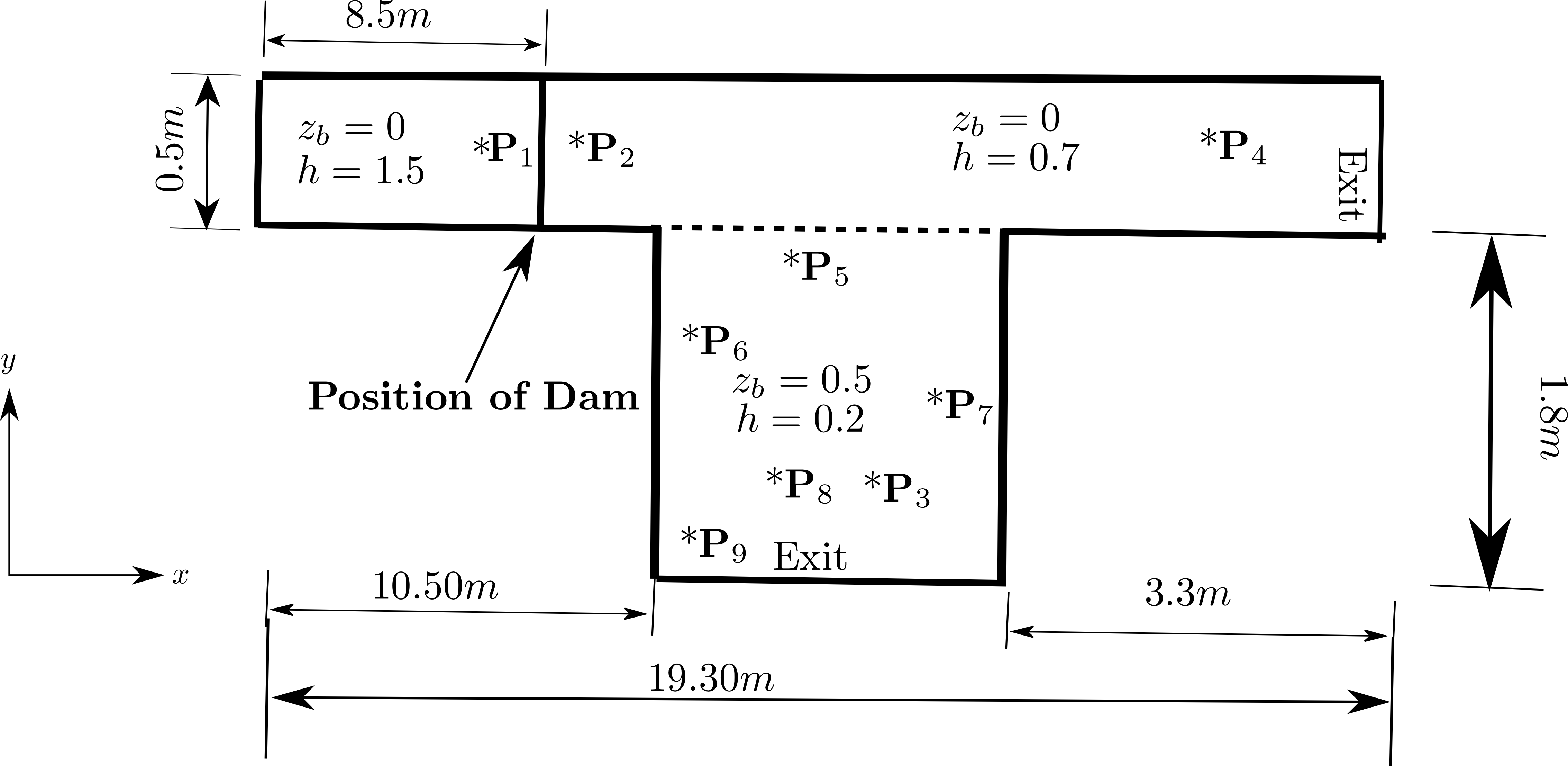}%
  }
\end{figure} 
We consider a dam break flow in a 19.3 meter long, 0.5 meter constant width
flat channel with adjacent flat floodplain
\citep{viseu3numerical, moralesetal2013, chineduthesis},
see figure \ref{fvmnumfigwhycouplinggeometry}.
The labels $P_1, P_2, \dots, P_6$ are chosen probe points in the flow
domain. To compare our simulations
with the literature we add bottom friction terms to the channel and the floodplain
models using a manning coefficient of $0.009$s/m$^{1/3}$. The boundaries
are all closed walls except the right side as indicated in the figure.
The wall elevation within the channel is $2.5$ meters. For $x\geq 12.5$
there is no channel wall.
The initial flow condition is given by
\begin{align*}
&	H(x,y,0) = \begin{cases}
						0.504, & \mbox{ at the reservoir, that is }  0\leq x \leq 6.10 \mbox{ and } 1.8\leq y \leq 2.3, \\
						0.003, & \mbox{ elsewhere},
		       \end{cases}\\
&	u(x,y,0) = v(x,y,0)	= 0 \quad \mbox{ everywhere}.	       
\end{align*}

To apply the VCM to this problem, we consider the following $\zwall$:
\begin{equation}
  \zwall  = \begin{cases}
  						   \tanh( (10.0-x) )+ 1.0, & \mbox{ if } x < 14.0, \\
                           0.0,  & \mbox{ if } x \geq 14.0.
              \end{cases}
\end{equation}

\subsubsection{Result of test 1:}
The four simulation methods were all run with a grid of $68 \times 90$ cells in the floodplain. For the channel region,
the full 2D simulation used $193\times25$ 2D cells, the VCM used 193 1D cells and an upper layer 2D grid of $193\times8$ cells
(that is, 8 2D upper layer cells per one 1D cell), while the HCM and FBM
used 193 1D cells each.
Each simulation was run for ten seconds with a CFL number of 0.95.
%

\begin{figure}[ht!]
	\includegraphics[width=\textwidth]{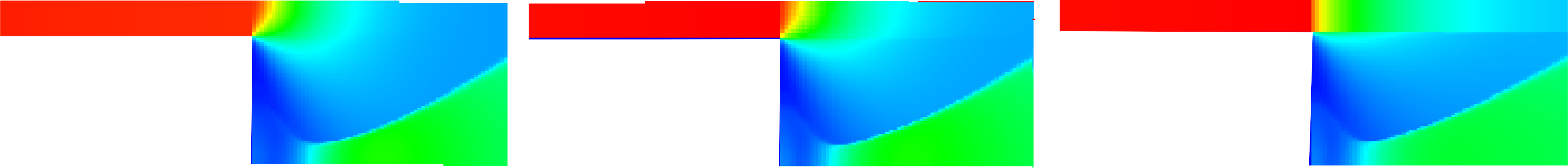}
  \caption{Comparison of free surface elevation for the different methods:
  Test 1. From left to right: full 2d, VCM, HCM} 
  \label{numfigtest1etaview}
\end{figure}

\begin{figure}[ht!]
  \begin{center}
    {\includegraphics[width=0.3\textwidth]{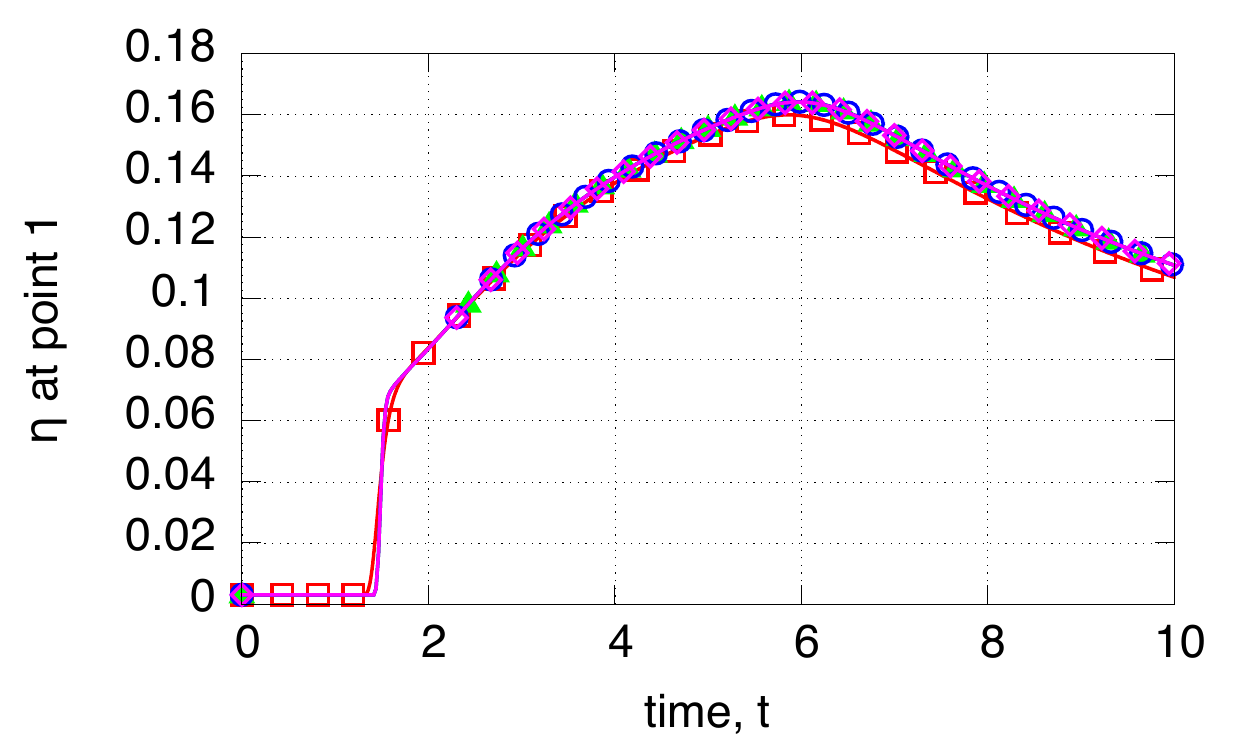}}
    {\includegraphics[width=0.3\textwidth]{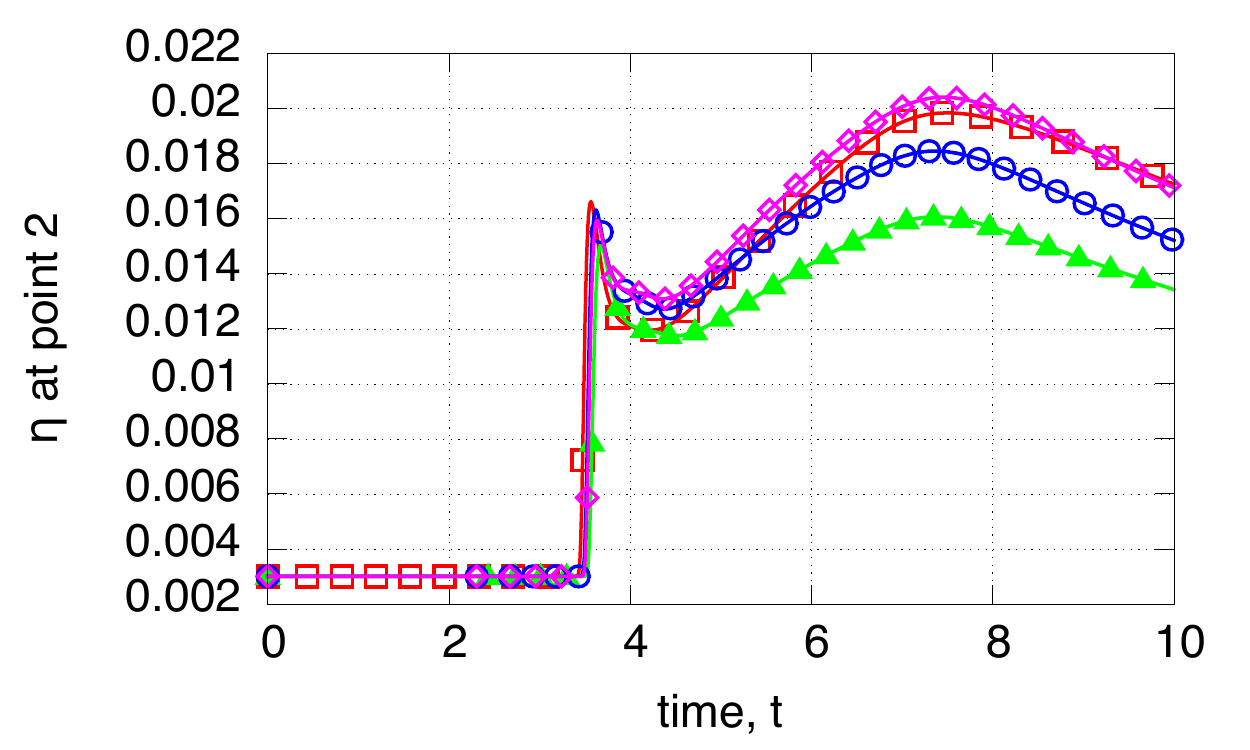}}
	  {\includegraphics[width=0.3\textwidth]{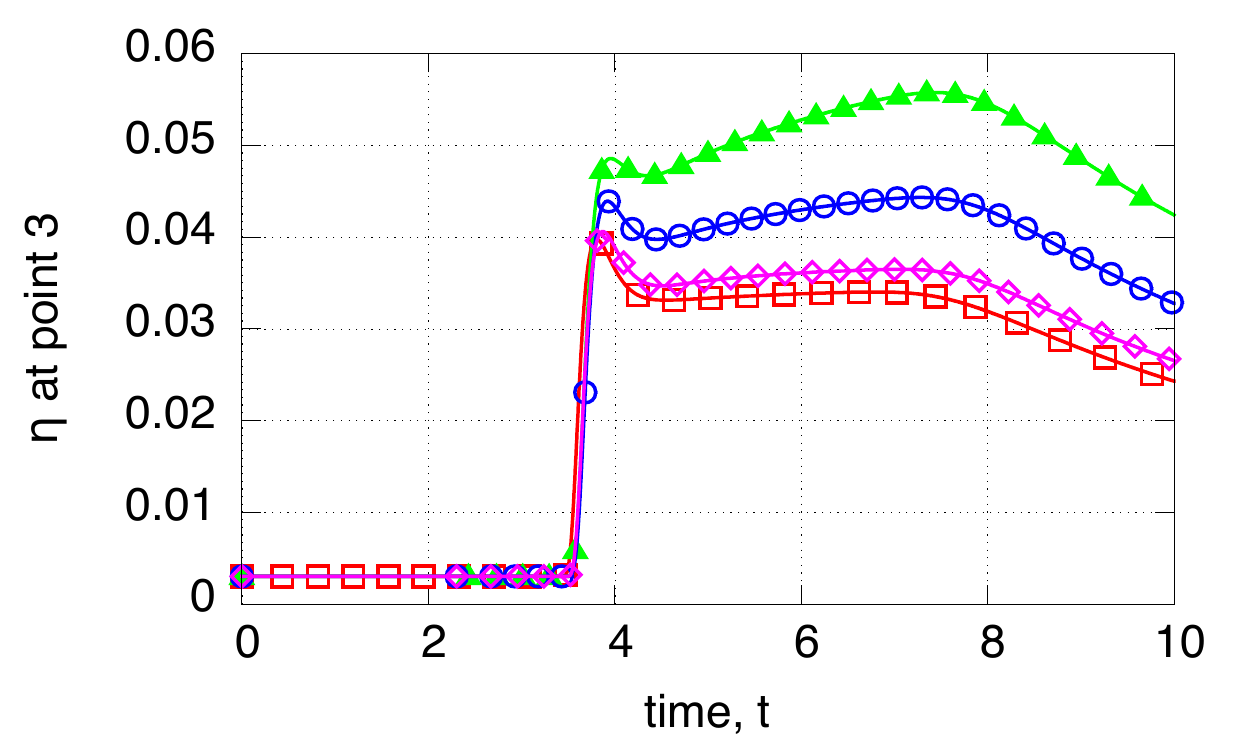}}

    {\includegraphics[width=0.3\textwidth]{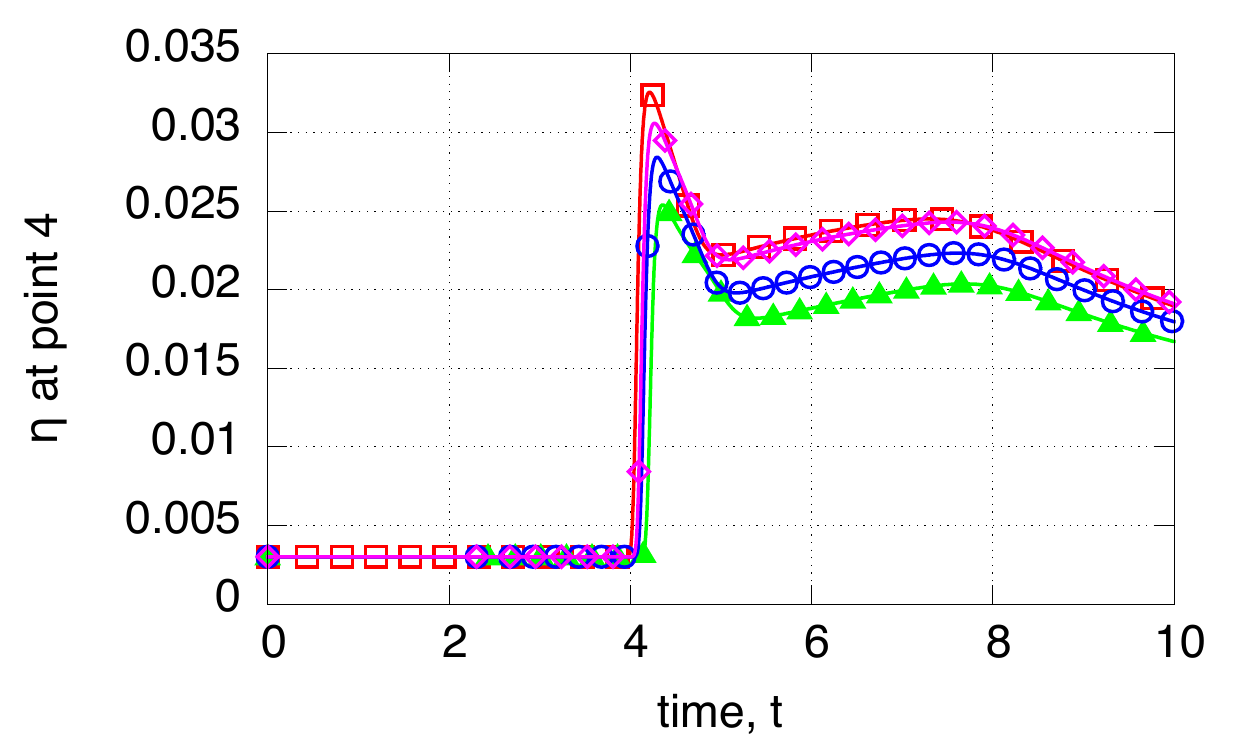}}
	  {\includegraphics[width=0.3\textwidth]{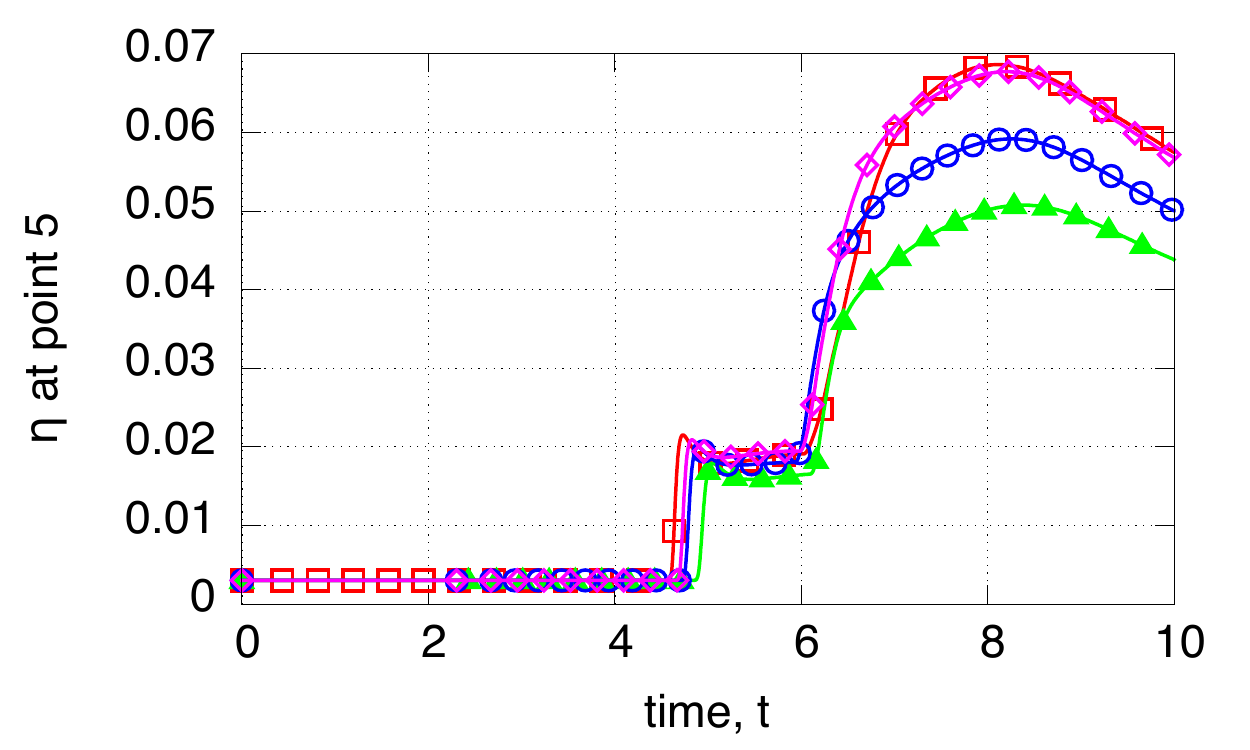}}
	  {\includegraphics[width=0.3\textwidth]{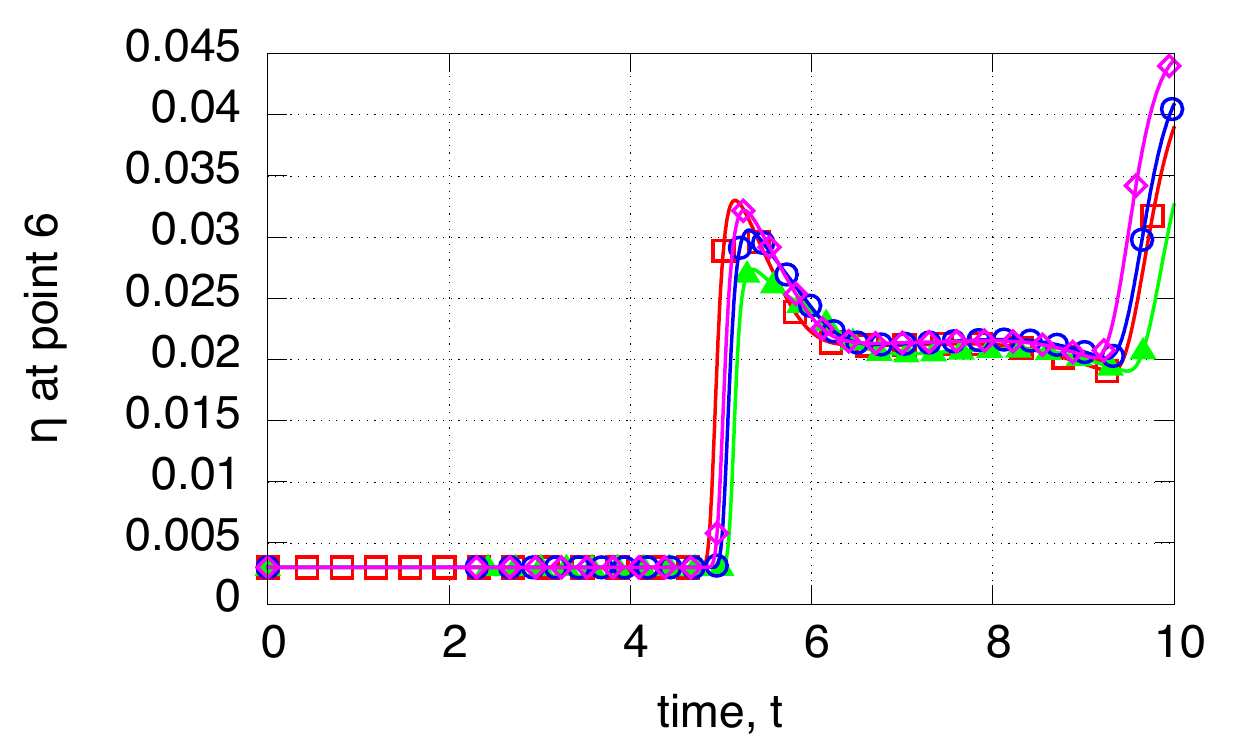}}

	  {\includegraphics[width=0.25\textwidth]{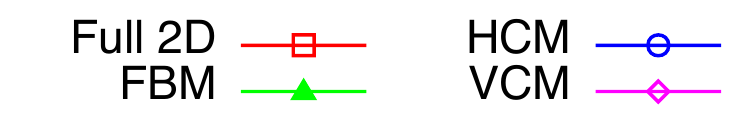}}
  \end{center}
  \caption{Comparison of free surface elevation at probe points : Test 1}
  \label{numfigtest1etaprobepoints}
\end{figure}

The free surface elevation at the last time step is shown in figure \ref{numfigtest1etaview} 
for the full 2D, VCM and HCM. The figure shows that the coupling methods, VCM and HCM, approximates
the full flow field with good accuracy, however, one can also see that the VCM computes better approximation
than the HCM.  Especially the 2D flow structure in the right part of the channel
where the upper layer is active, is correctly captured by the VCM, unlike the HCM.
The accuracy of the vertical coupling method is further  illustrated in figure \ref{numfigtest1etaprobepoints} 
which displays the time evolution of the free surface elevation, at the selected probe points. It can be seen that the VCM captures the
full 2D results better than both the HCM and FBM, at all the probe points. 
Again, from figure \ref{numfigtest1etaview}, one can see that the VCM recovers 2D flow structure within the channel
unlike the HCM and FBM. 
%
%
For the this test case, the vertical coupling method results in about 48\% gain
in computational time over the full 2D simulation. The fastest method (FMB)
leads to a gain of about 55\%.
Hence in this test case, the VCM truly has shown good accuracy improving on the simple flux
based coupling method while retaining most of the gain in efficiency.


\subsection{Test case 2 : Channel Flow into Elevated 2D Floodplain}\label{numsectest2}
The second test case was suggested in \cite{chineduthesis}, see also \cite{chineduandreas1}.
The set up consists of a channel of length 19.3 metre, width 0.5 metres, zero bottom topography
and connected to a 0.5 metre high floodplain which is located in the
region, $(x,y) \in [10.5, 16.0] \times [0, 1.8]$,  see figure \ref{numfiggeometrytest2}. 
The manning coefficient for both channel and floodplain is taken as 0.009$s/m^{1/3}$ the boundaries
are only open at the sides indicated "exit" in figure \ref{numfiggeometrytest2}, others are closed.
Nine probe points are considered. 
From the bottom topographies for the channel and floodplain given above,
the channel wall elevation is equal to $3.0$ meters within the channel
except in a breach region for $10.5\leq x \leq 16.0$ where the wall height
drops to $0.5$ meters.
The initial condition is the following.
\begin{align}
 & H(x,y,0)  = \begin{cases} 1.5, & \mbox{ if } x \leq 8.5, \quad y \geq 1.8, \\
                            0.7, & \mbox{ if } x > 8.5, \quad y \geq 1.8, \\
                            0.2, & \mbox{ if } 10.5\leq x \leq 16.0, \quad 0 \leq y \leq 1.8, \\
                            0.0, & \mbox{ else}.
               \end{cases} \\
& u(x,y,0) = v(x,y,0) = 0.               
\end{align}

To apply the vertical coupling method to this problem, we consider the following smoother version 
of the wall elevation:
\begin{equation}\label{numeqntest2etabeta}
    \zwall  = \begin{cases} \tanh( 0.5(4.5-x) )+1.5, & \mbox{ if } x < 10.0, \\
                                    0.5,    & \mbox{ if }  10.0 \leq x \leq 16.5,  \\
                                    \tanh(x-19.2)+1.5, & \mbox{ elsewhere}.
       \end{cases}
\end{equation}

\subsubsection{Result of test 2:}
For all four methods we simulated this problem with $55\times90$ grid cells in the
floodplain, while the grid resolutions in the channel are exactly the same as in
the first test case. 
The simulation was run for ten seconds and CFL of $0.95$.

\begin{figure}[p] 
	\includegraphics[width=\linewidth]{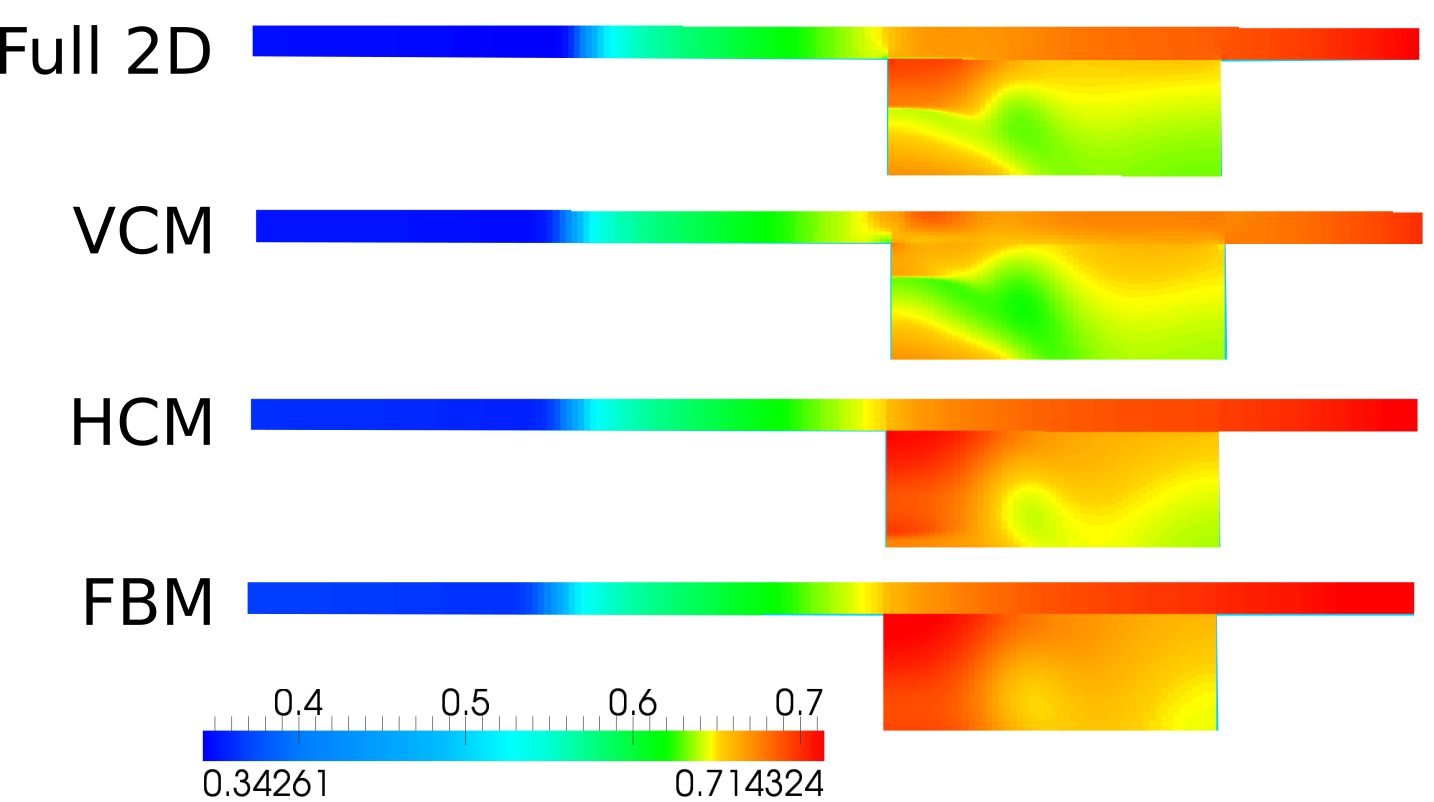}
  \caption{Comparison of free surface elevation for the different methods: Test 2 } 
  \label{numfigtest2etaview}
\end{figure}

\begin{figure}[p] 
	\includegraphics[width=\linewidth]{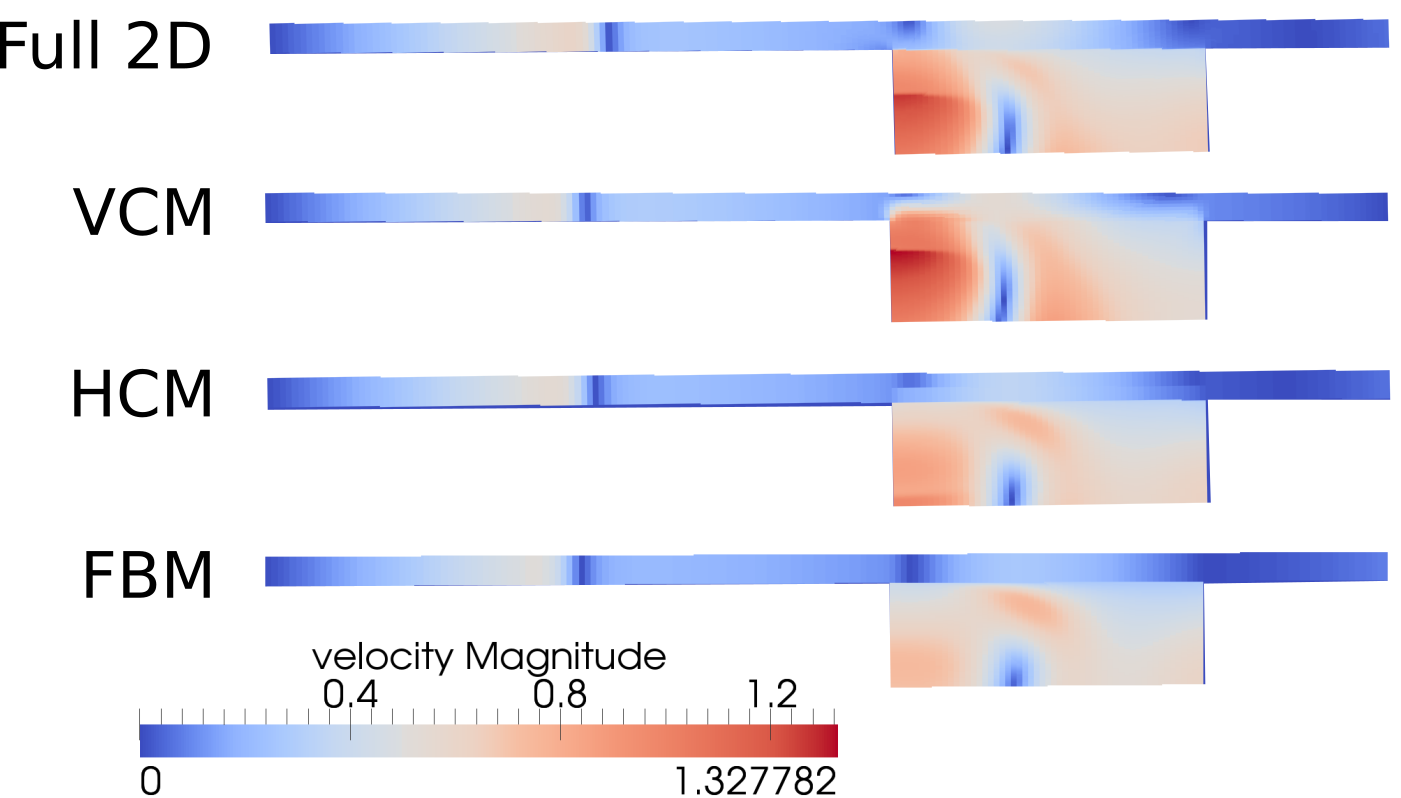}
  \caption{Comparison of velocity magnitude for the different methods: Test 2} 
  \label{numfigtest2velomagview}
\end{figure}

\begin{figure}[p] 
  \begin{center}
	\subfigure{\includegraphics[width=0.3\textwidth]{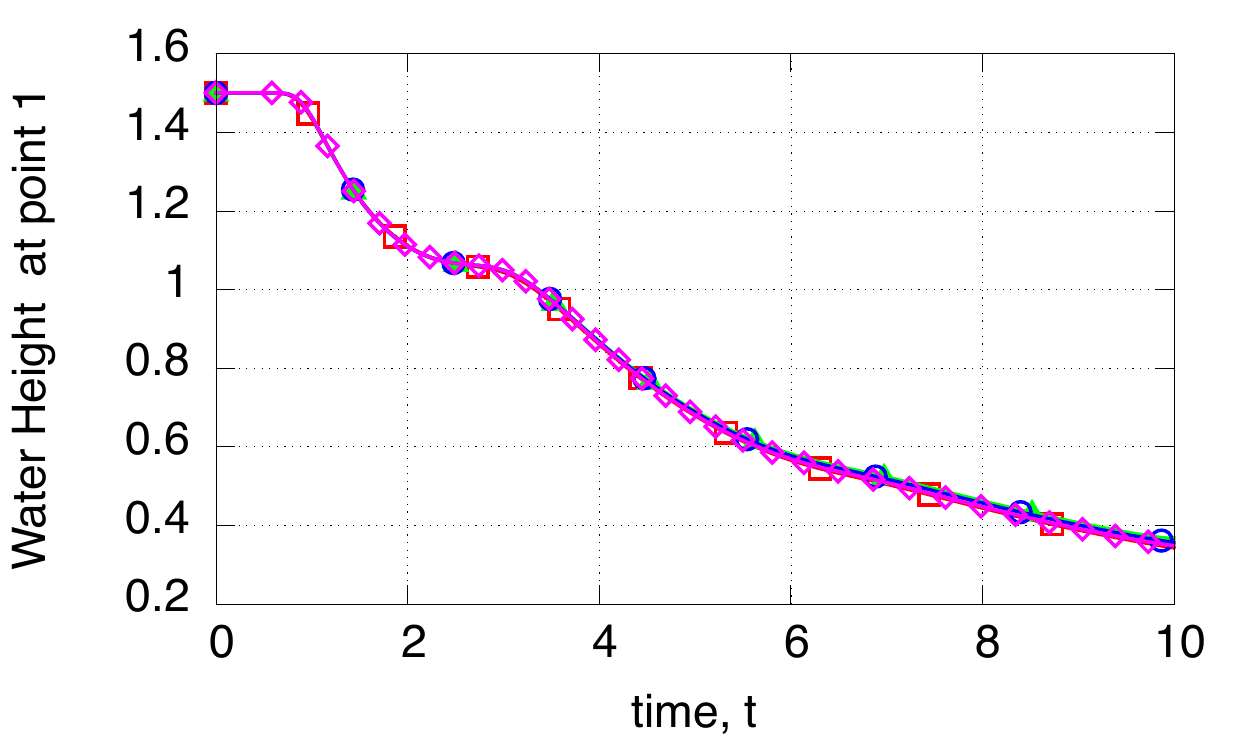}}
 	\subfigure{\includegraphics[width=0.3\textwidth]{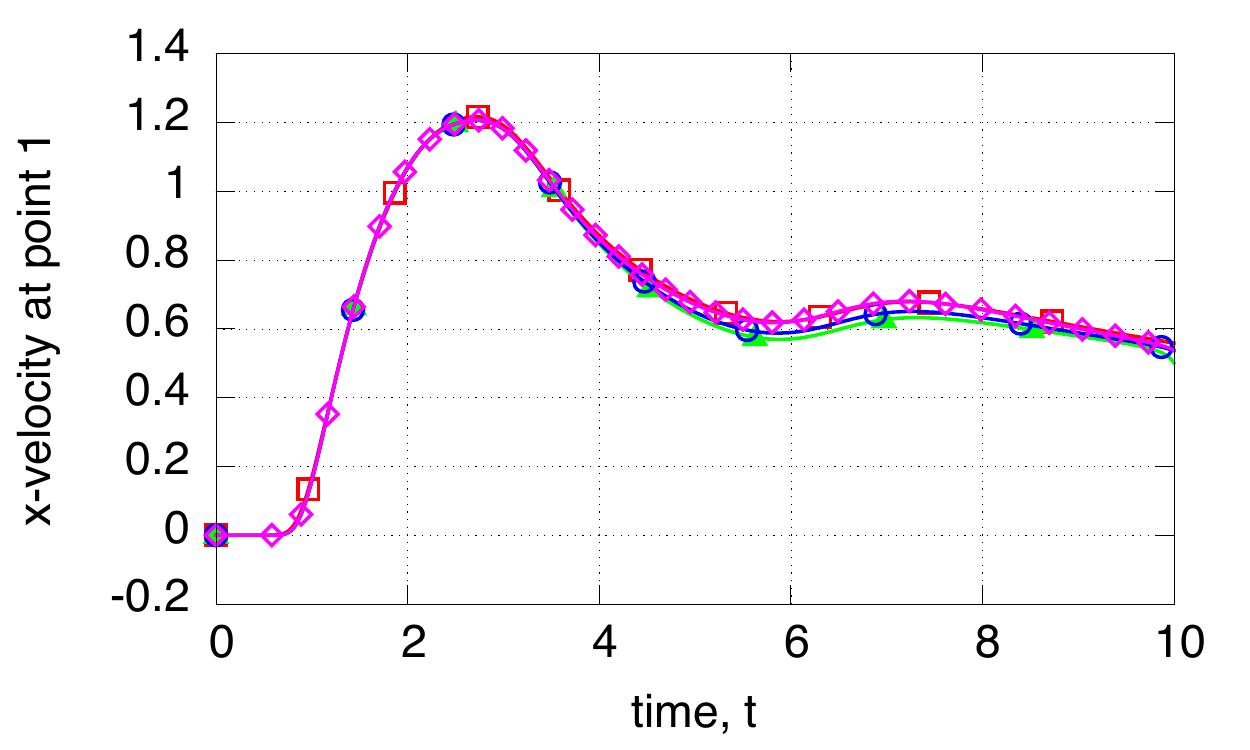}}
 	\subfigure{\includegraphics[width=0.3\textwidth]{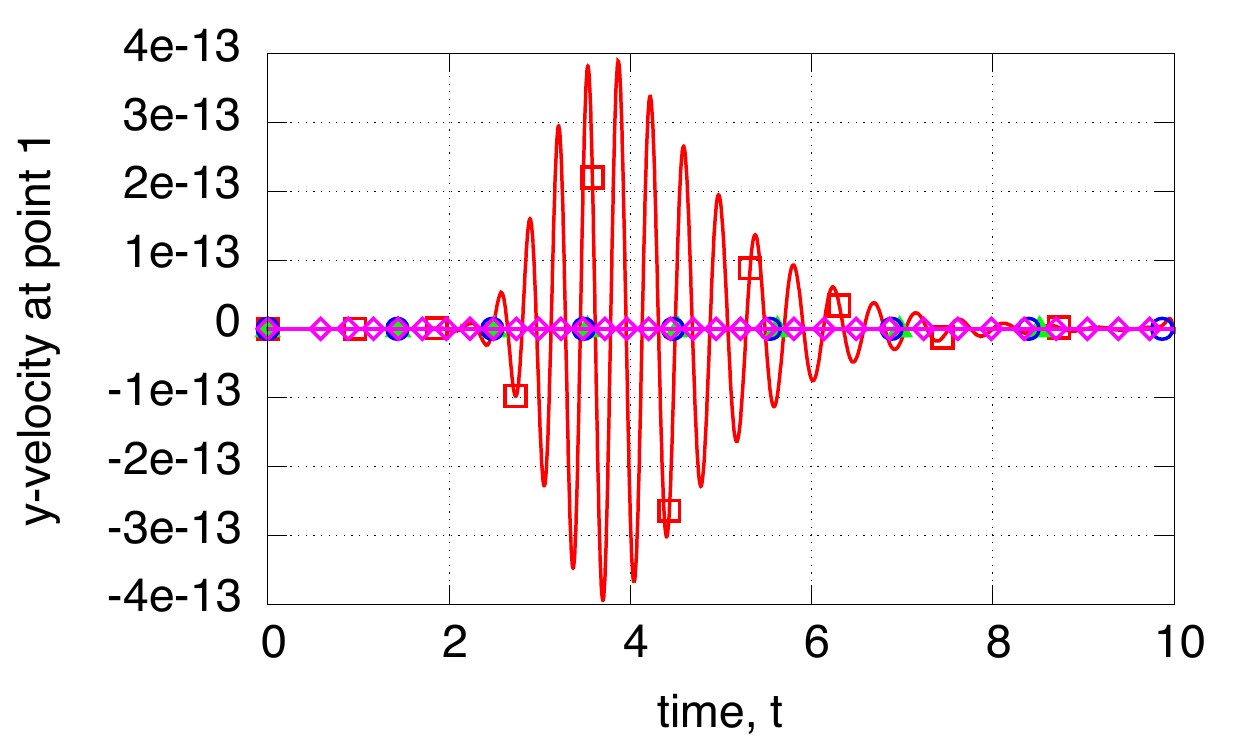}}

	\subfigure{\includegraphics[width=0.3\textwidth]{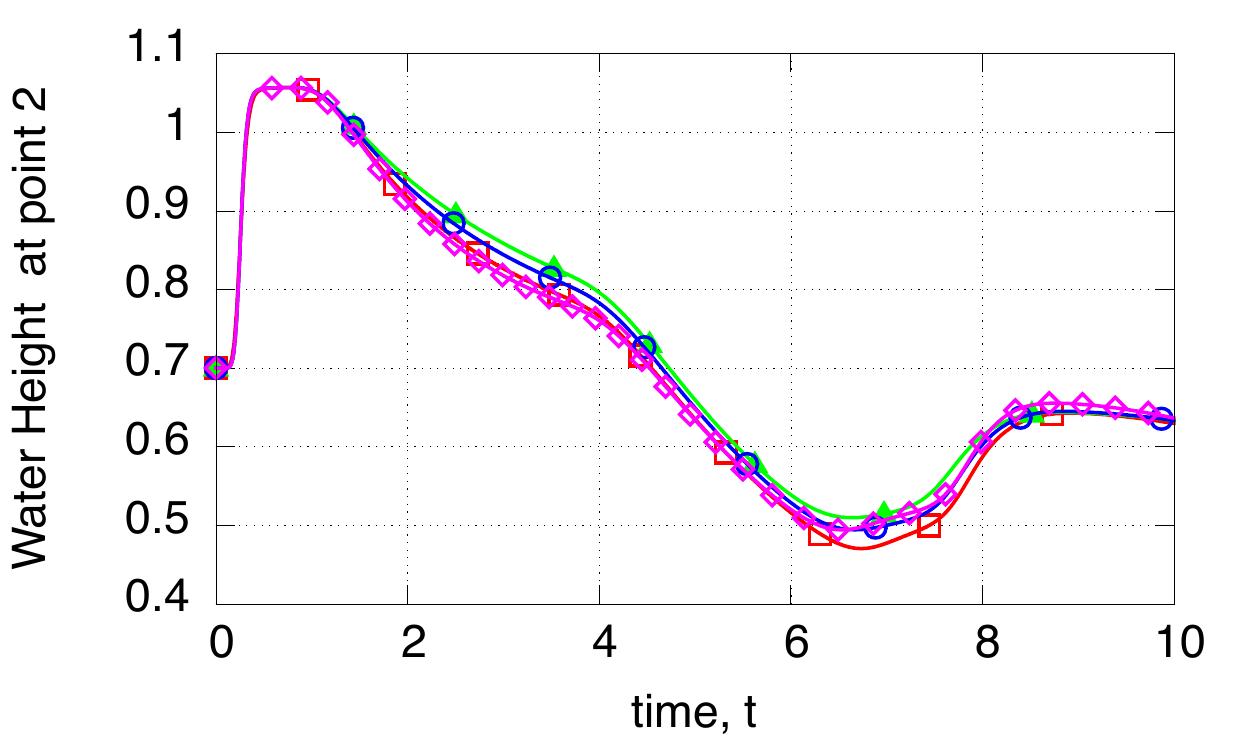}}
	\subfigure{\includegraphics[width=0.3\textwidth]{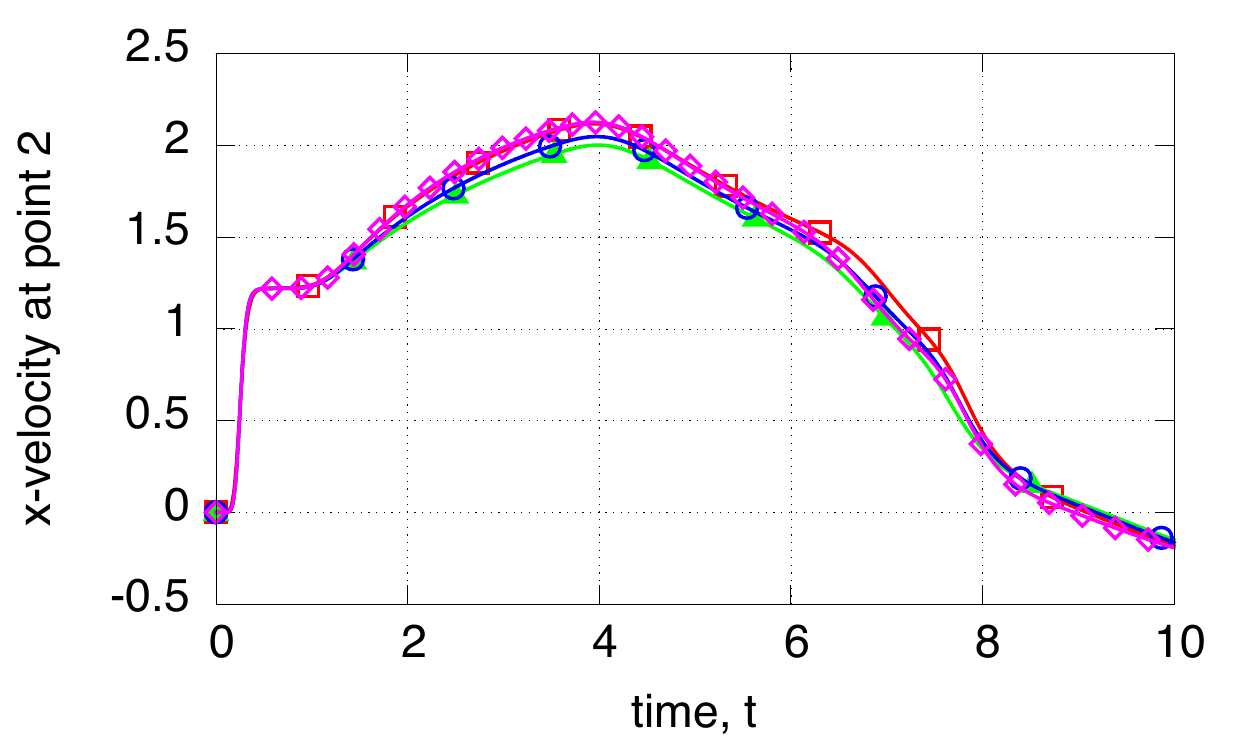}}
	\subfigure{\includegraphics[width=0.3\textwidth]{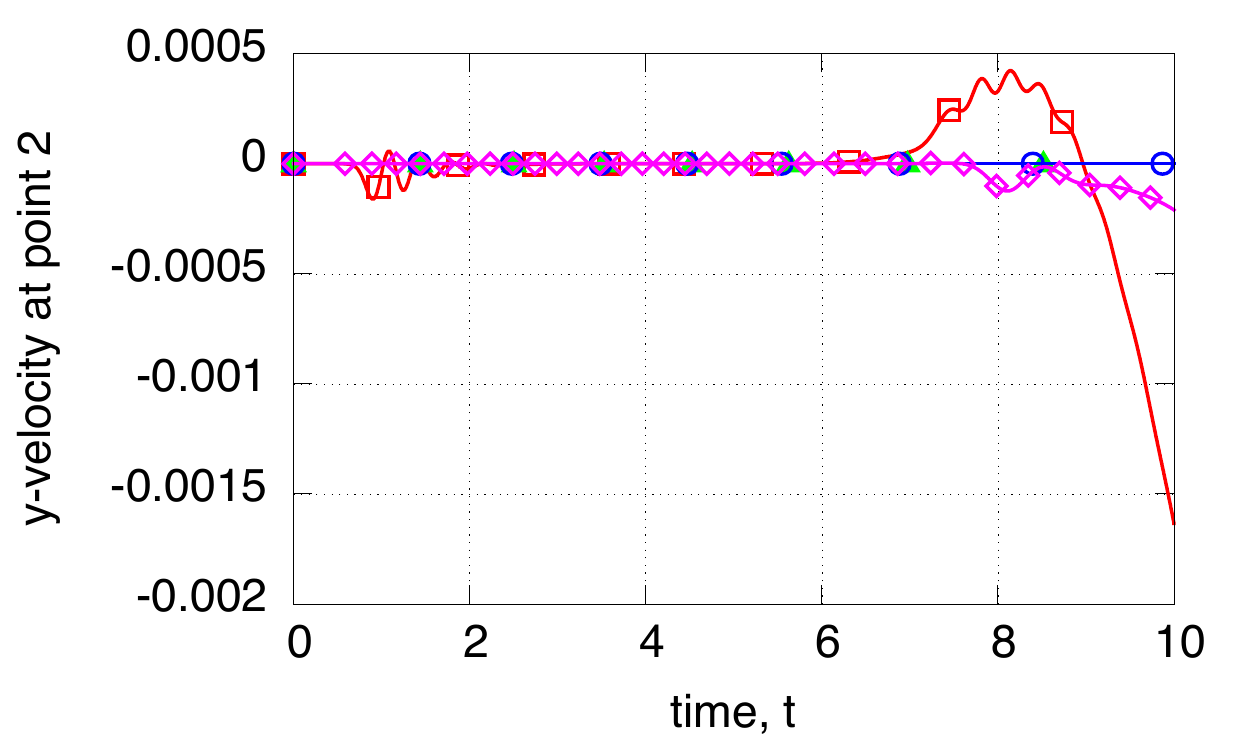}}

	\subfigure{\includegraphics[width=0.3\textwidth]{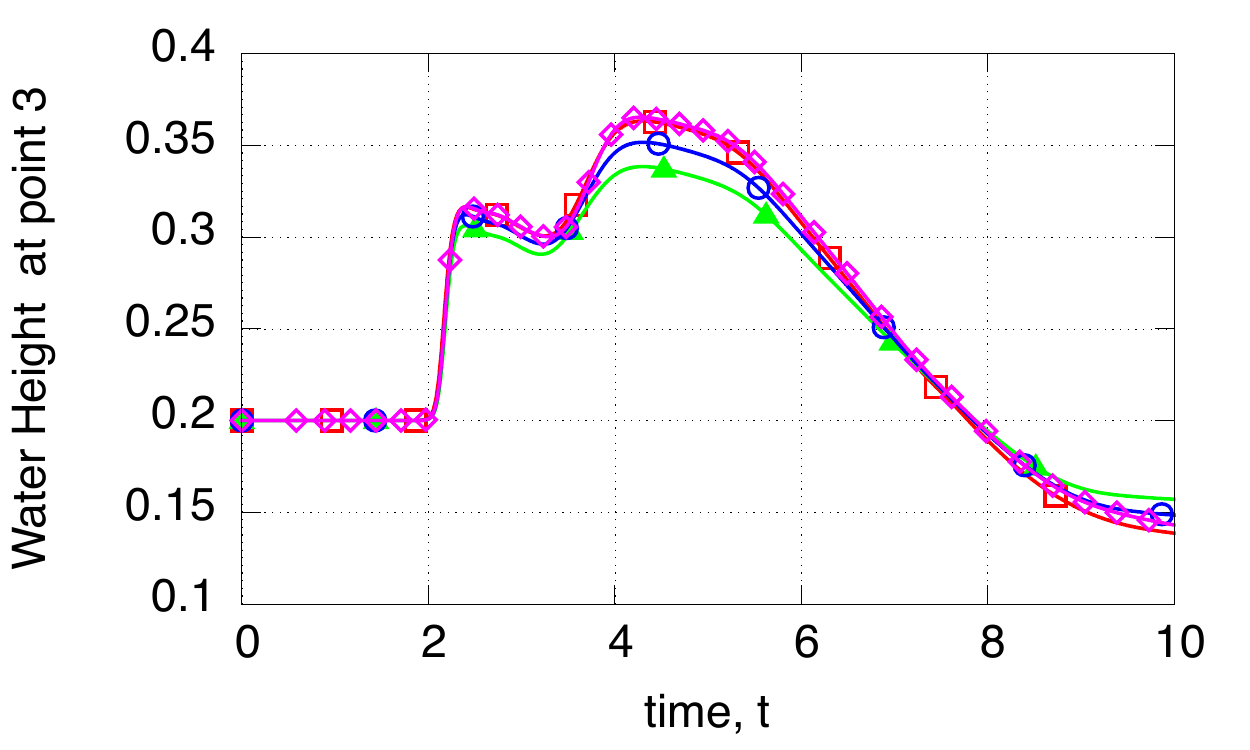}}
	\subfigure{\includegraphics[width=0.3\textwidth]{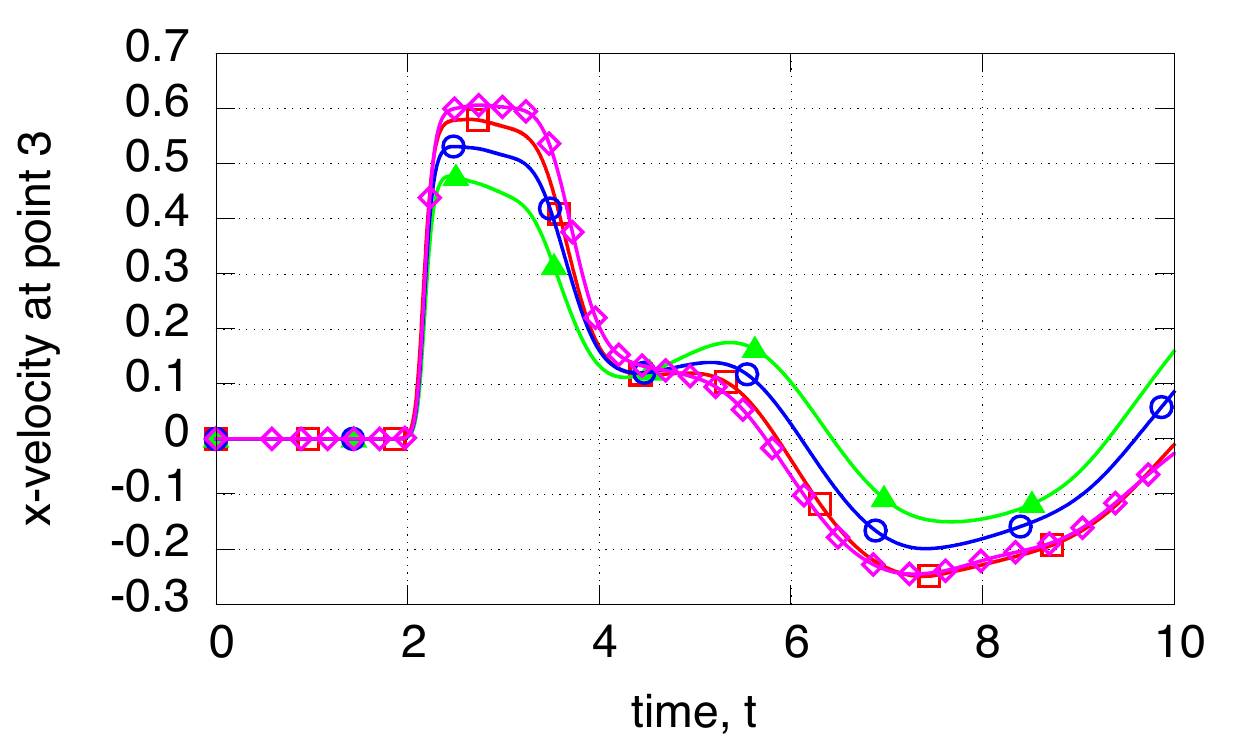}}
	\subfigure{\includegraphics[width=0.3\textwidth]{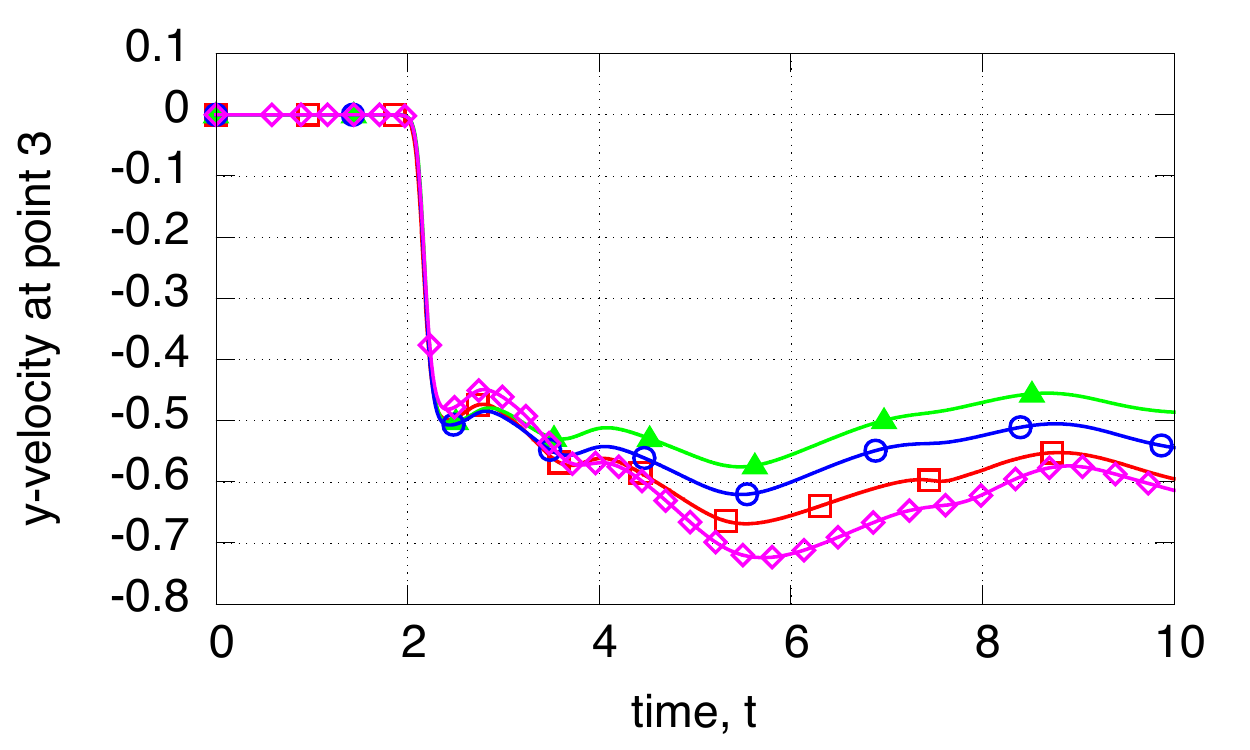}}

	\subfigure{\includegraphics[width=0.3\textwidth]{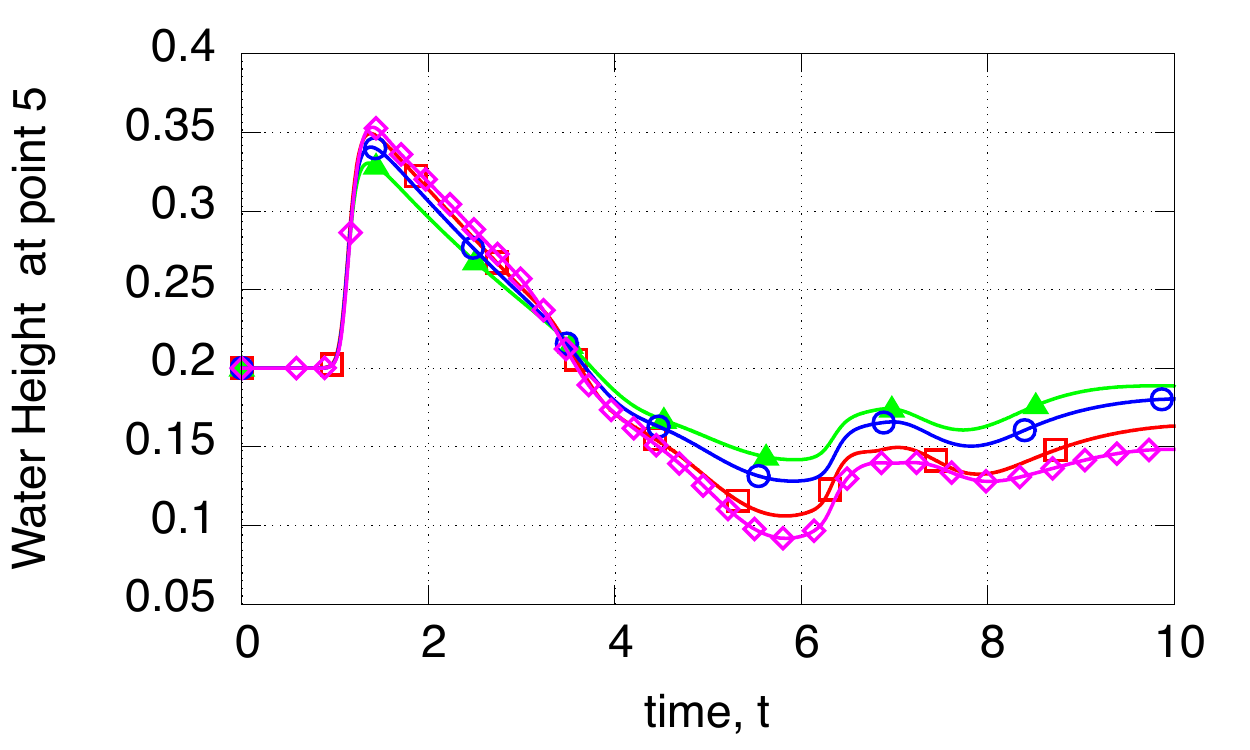}}
	\subfigure{\includegraphics[width=0.3\textwidth]{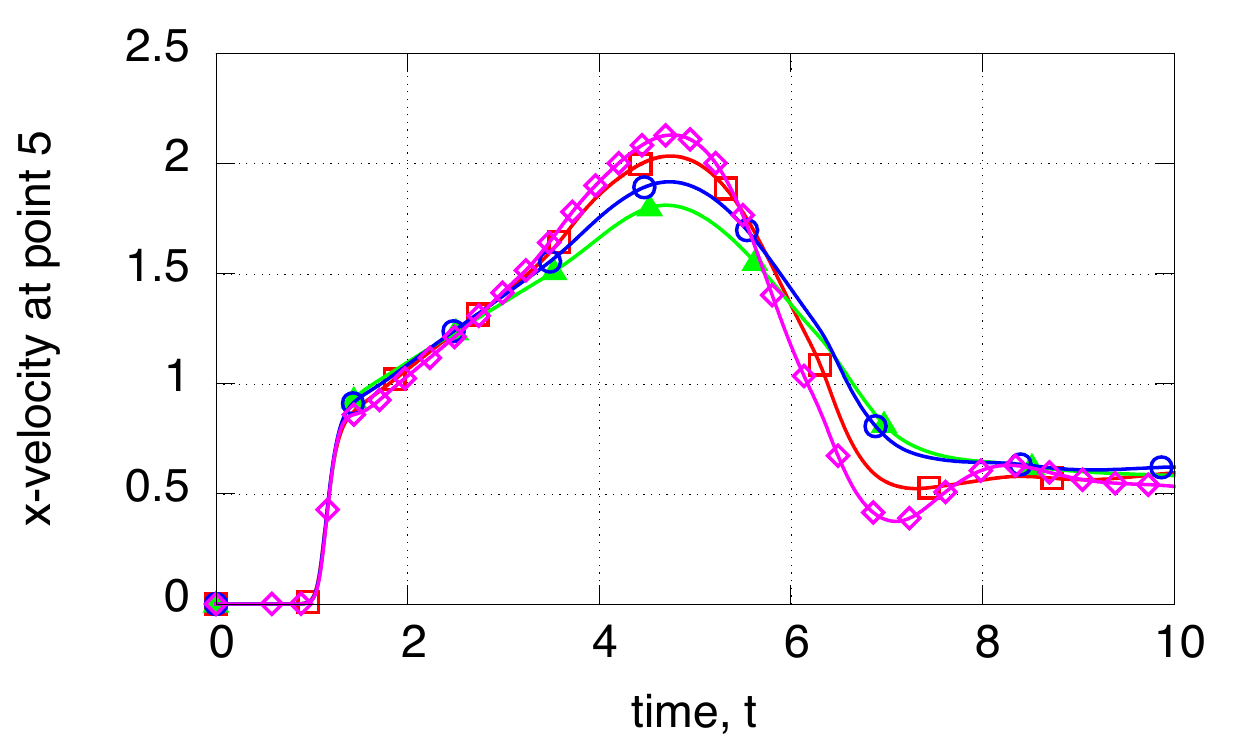}}
	\subfigure{\includegraphics[width=0.3\textwidth]{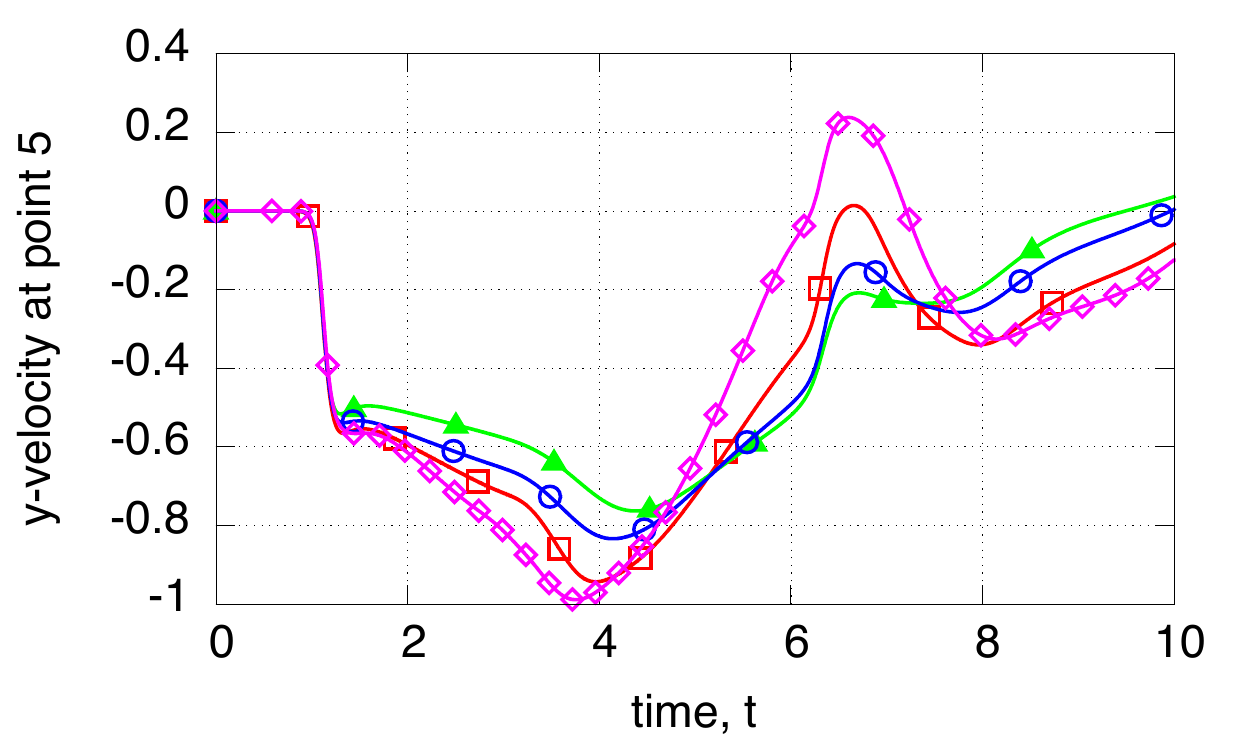}}


	\subfigure{\includegraphics[width=0.3\textwidth]{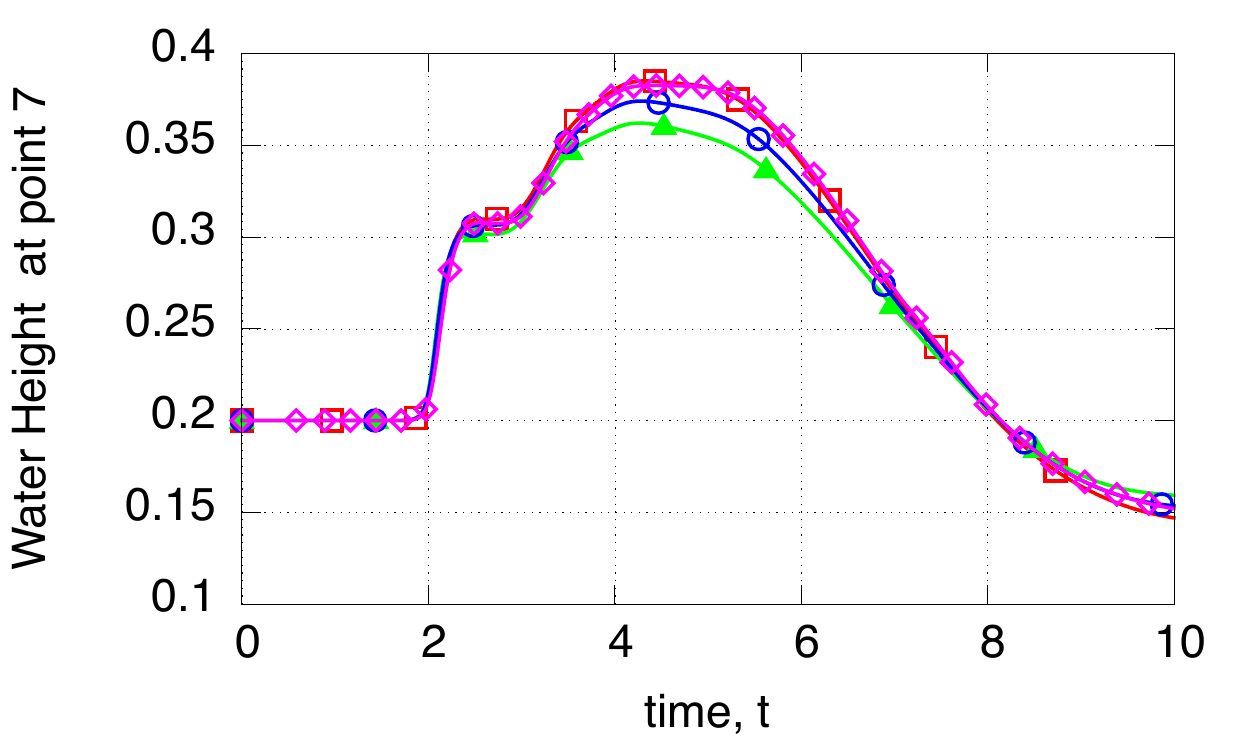}}
	\subfigure{\includegraphics[width=0.3\textwidth]{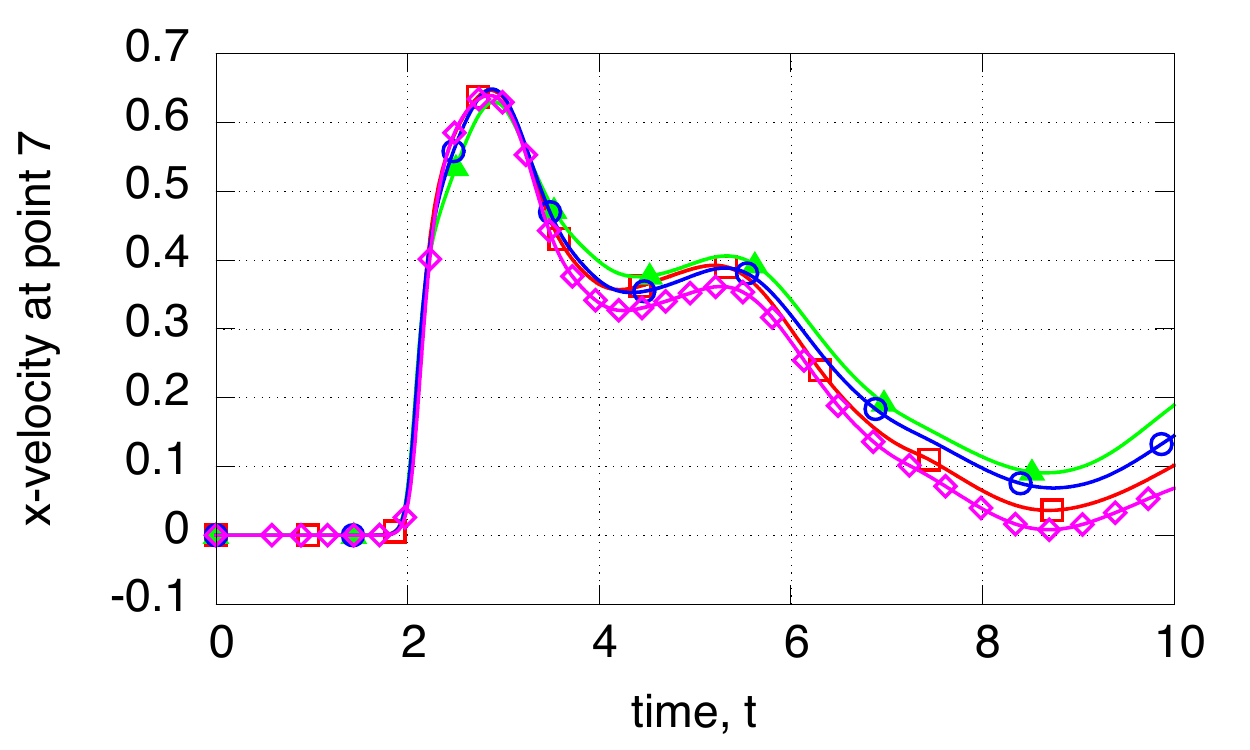}}
	\subfigure{\includegraphics[width=0.3\textwidth]{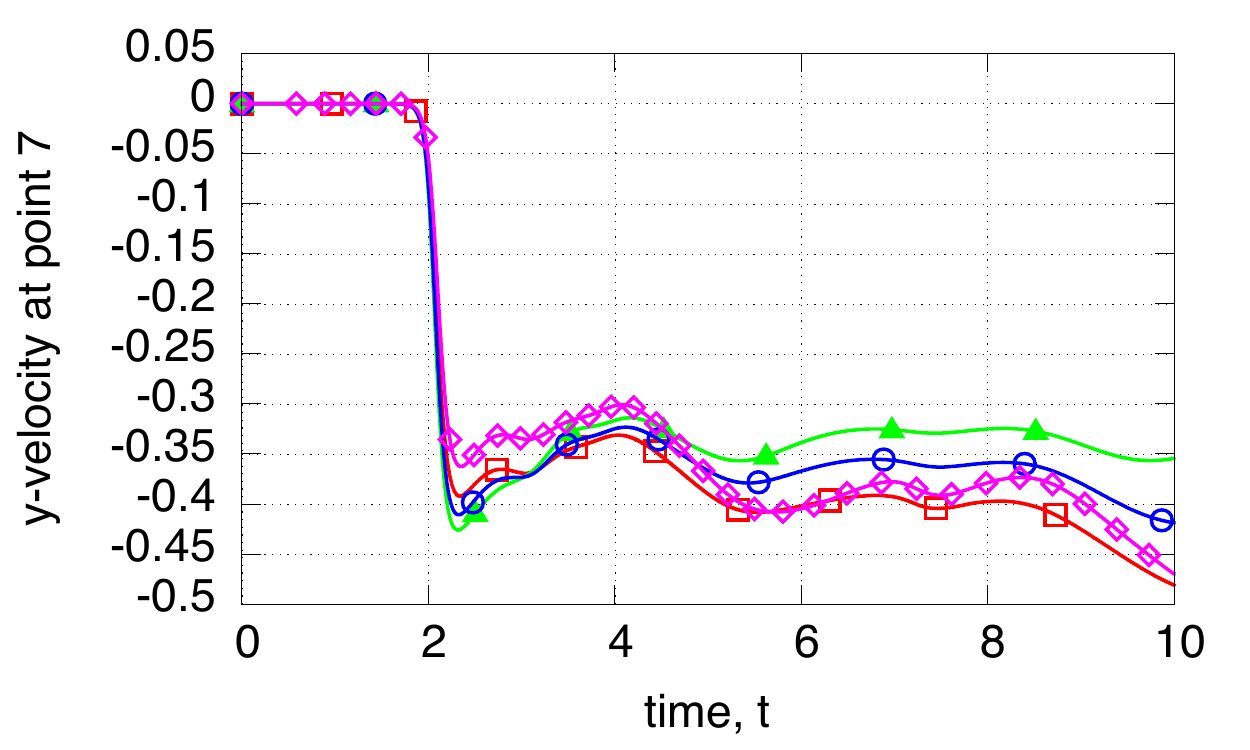}}

	\subfigure{\includegraphics[width=0.3\textwidth]{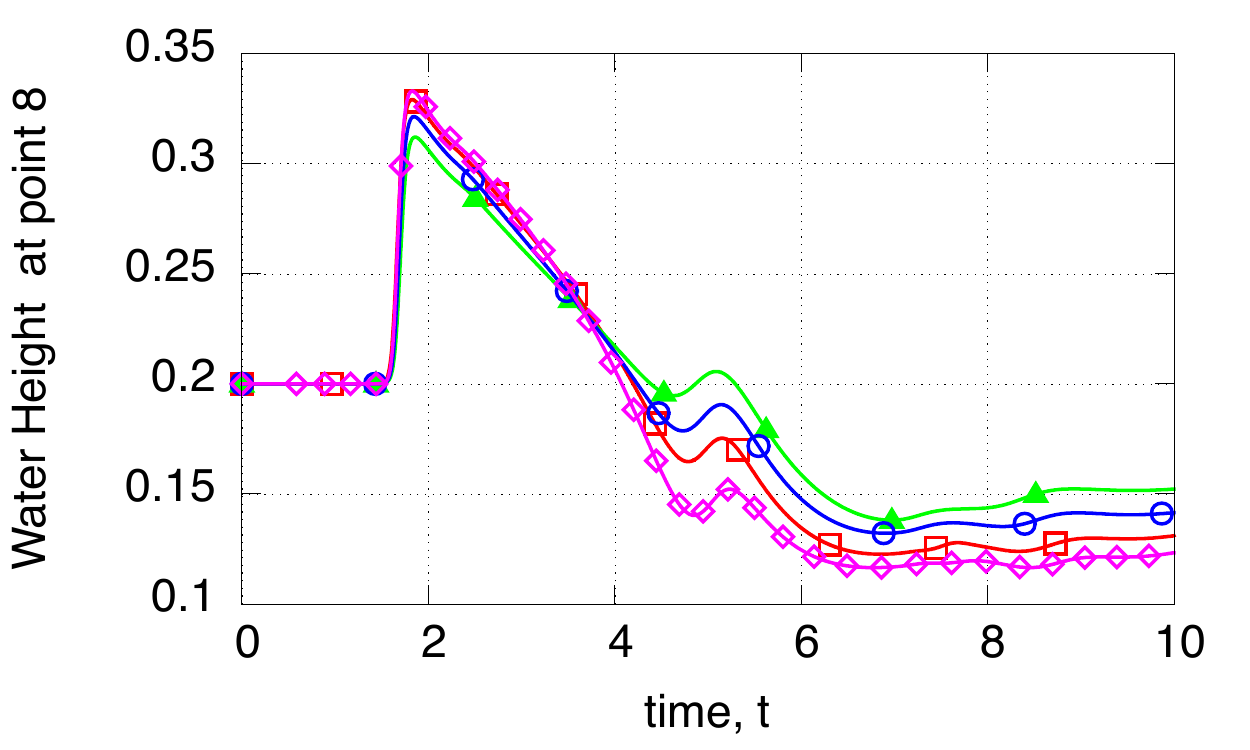}}
	\subfigure{\includegraphics[width=0.3\textwidth]{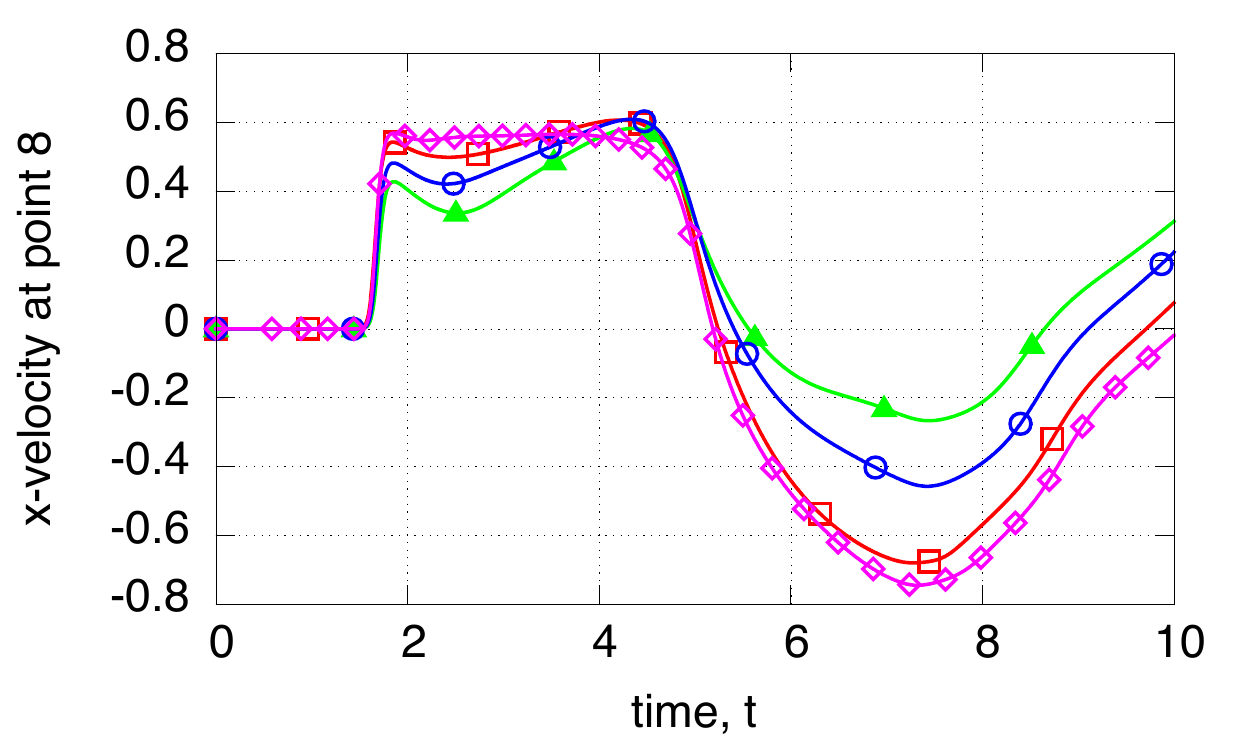}}
	\subfigure{\includegraphics[width=0.3\textwidth]{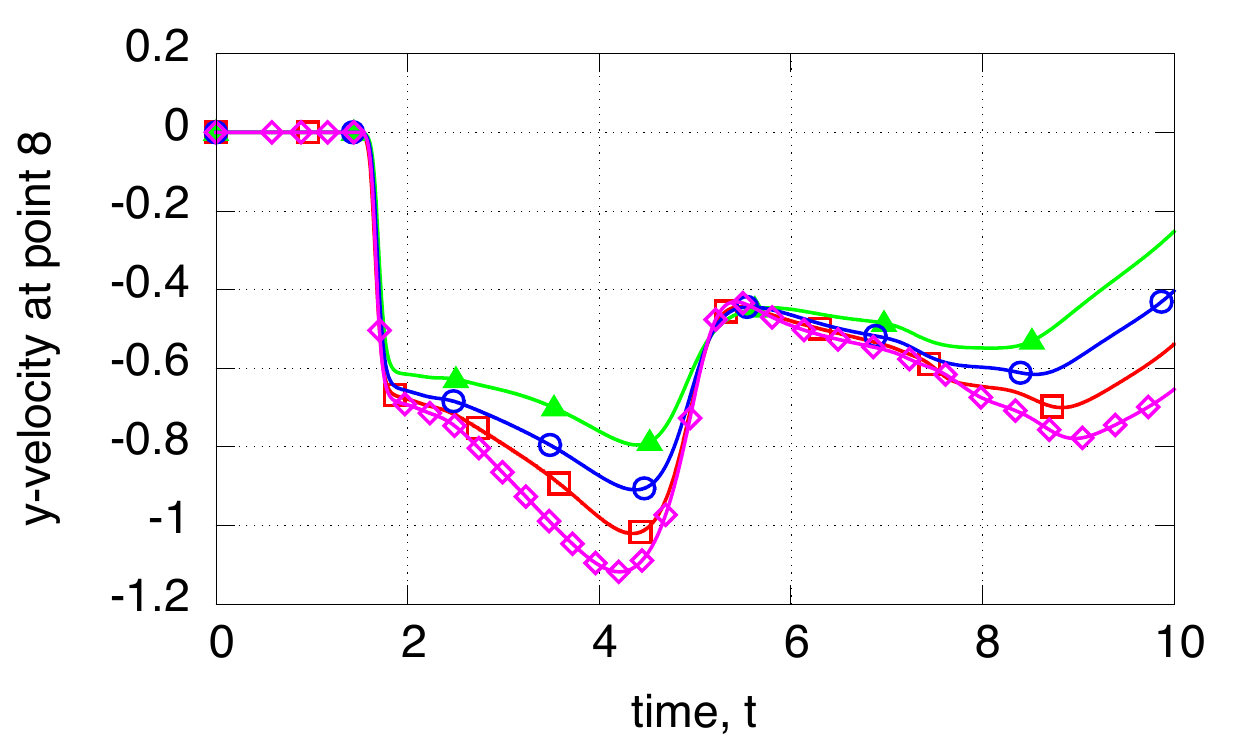}}


	\subfigure{\includegraphics[width=0.25\textwidth]{key}}
  \end{center}
  \caption{Comparison of water height at velocities at probe points : Test 2} 
  \label{numfigtest2etaprobepoints}
\end{figure}

The free surface elevation and velocity magnitude are shown in
figures \ref{numfigtest2etaview} and \ref{numfigtest2velomagview} respectively.
We can see that that the VCM computes a better approximation of the full 2D solution
than the HCM which in turn, is more accurate than the FBM. Moreover, the VCM
reproduces a non laterally constant free surface elevation and velocity
within the channel, which can not be achieved by either HCM or FBM.
To further understand the results, the free surface elevation,
the $x$-component and $y$-component velocity are plotted for selected probe points
in figures \ref{numfigtest2etaprobepoints}.  It can be seen that the vertical coupling
method is more accurate than the other methods for all three flow quantities at all
probe points and almost throughout the duration of simulation.
Again, the vertical coupling method really captures the flow structure of the
full 2D simulation.

Finally, figures \ref{numfigtest1etaview} and \ref{numfigtest2velomagview} show that the VCM, unlike the other methods,
recovers the 2D flow structure within the channel at the flooding regions. Again, the VCM continues to compute
1D solutions at non flooding regions; this demonstrates the self-adaptive nature of the method.

In this example VCM was about 75\% slower then the computationally much
simpler FBM method and about 35\% more efficient then the full 2D
simulation.

For this problem the VCM clearly improved the accuracy of the simulation
both within the channel and in the floodplain compared to simpler coupling
method but at a increased computational cost. But note that in this examples
the percentage of the domain where a 1D assumption for the flow is valid is
quite small so that it is not surprising that the gain in computational
efficiency between full 2D and VCM is not so large.
For a simulation of a very large network of rivers where the 2D region might
be very small compared to the entire computational domain, the difference
in efficiency of one coupling method over another will be far less significant,
so that accuracy becomes the deciding factor.
In this regard, the VCM is the best method of the three considered here.

\section{Conclusion}\label{vsecconc}
In this paper we investigated a new method for coupling 2D and 1D shallow water
models. We focused on the need to efficiently recover 2D channel flow structure
during a flooding event to achieve high accuracy while maintaining the efficiency of
the 1D channel model as much as possible. To this end we presented a
vertical coupling approach which adds a second 2D layer inside the channel
in regions where overflow can be expected to occur. This makes coupling the
channel flow to the floodplain straightforward since 2D information is
always available also within the channel when required.
Any standard 1D channel flow  and 2D floodplain solver can be used as building
blocks for the new VCM method and only a slight modification of standard 2D
solvers is required for the evolution of the second layer.
We proved that the resulting method retains many properties of the
1D and 2D solvers used, e.g., mass conservation and well balancing. In
addition we also studied a \emph{no-numerical flooding property}.
Our numerical results show that the VCM, in most cases, outperforms the
other methods studied in this paper and does
accurately recover 2D flow structures also within the channel when flooding
occurs.

More detailed investigations with more complex channel geometries will be
the focus of our further research.
As mention in this paper, VCM actually consists of a very large family of
methods. It depends on the choice of $\zwall$ which is used to determine
when a channel is considered full and the 2D layer model is used.
For example, taking $\zwall=\infty$ everywhere leads to a horizontal flux
coupling method similar to the FBM method. On the other hand taking
$\zwall=0$ in regions in danger of flooding and $\zwall$ very large
away from these regions using a smooth transition,
results in a frontal type coupling method where a standard 2D solver is used in the
region where $\zwall=0$ and on the boundary of this region the VCM will
lead to a blending type approach between the 2D and the 1D regions. The
width of the blending region will depend on how $\zwall$ changes from $0$
to $\infty$.
Further tests are required to understand the influence of different choices
of $\zwall$ and are the focus of ongoing work.

\section*{Acknowledgement} 
We are grateful to the Petroleum Technology Development Fund (PTDF), Nigeria for funding
this study and to the Centre for Scientific Computing, University of Warwick for providing
the computing resources.


\bibliographystyle{abbrv}
\bibliography{mybib}            
\end{document}